\documentclass{article}

\usepackage{amsfonts}

\usepackage{amsthm}

\usepackage{amssymb}

\usepackage{amsmath}

\theoremstyle{definition} \newtheorem{propa}{Definition}[subsection]

\theoremstyle{definition} \newtheorem{coarse}[propa]{Definition}

\theoremstyle{plain} \newtheorem{coarseinv}[propa]{Proposition}

\theoremstyle{remark} \newtheorem{cinvrem}[propa]{Remark}

\theoremstyle{definition} \newtheorem{coarseembed}[propa]{Definition}

\theoremstyle{remark} \newtheorem{cerem}[propa]{Remarks}

\theoremstyle{plain} \newtheorem{embedhilbert}[propa]{Theorem}

\theoremstyle{definition} \newtheorem{bg}{Definition}[subsection]

\theoremstyle{definition} \newtheorem{equivdef}[bg]{Definition}

\theoremstyle{definition} \newtheorem{partition}[bg]{Definition}

\theoremstyle{plain} \newtheorem{equivalences}[bg]{Theorem}

\theoremstyle{remark} \newtheorem{finite}{Example}[subsection]

\theoremstyle{remark} \newtheorem{trees}[finite]{Example}

\theoremstyle{definition} \newtheorem{equi}{Definition}[subsection]

\theoremstyle{plain} \newtheorem{subequi}[equi]{Lemma}

\theoremstyle{plain} \newtheorem{gluing}[equi]{Theorem}

\theoremstyle{plain} \newtheorem{fibreing}[equi]{Corollary}

\theoremstyle{plain} \newtheorem{products}[equi]{Corollary}

\theoremstyle{plain} \newtheorem{covers}[equi]{Proposition}

\theoremstyle{definition} \newtheorem{fad}[equi]{Definition}

\theoremstyle{plain} \newtheorem{fad2}[equi]{Corollary}

\theoremstyle{plain} \newtheorem{unions}[equi]{Corrolary}

\theoremstyle{plain} \newtheorem{unions2}[equi]{Corollary}

\theoremstyle{plain} \newtheorem{unions3}[equi]{Corollary}

\theoremstyle{definition} \newtheorem{wmetric}{Definition}[subsection]

\theoremstyle{remark} \newtheorem{integers}[wmetric]{Example}

\theoremstyle{plain} \newtheorem{wmce}[wmetric]{Proposition}

\theoremstyle{plain} \newtheorem{subgp}[wmetric]{Lemma}

\theoremstyle{plain} \newtheorem{dlim}[wmetric]{Proposition}

\theoremstyle{plain} \newtheorem{ext}[wmetric]{Theorem}

\theoremstyle{plain} \newtheorem{produ}[wmetric]{Corollary}

\theoremstyle{plain} \newtheorem{wreath}[wmetric]{Corrolary}

\theoremstyle{definition} \newtheorem{bassserre}[wmetric]{Definition and lemma}

\theoremstyle{definition} \newtheorem{metrictree}[wmetric]{Definition and lemma}

\theoremstyle{plain} \newtheorem{fprod}[wmetric]{Theorem}

\theoremstyle{plain} \newtheorem{hnn}[wmetric]{Corollary}

\theoremstyle{remark} \newtheorem{cecprop}[wmetric]{Remark}

\theoremstyle{definition} \newtheorem{kernels}{Definitions}[subsection]

\theoremstyle{remark} \newtheorem{kerex}[kernels]{Examples}

\theoremstyle{plain} \newtheorem{kernelrem}[kernels]{Lemma}

\theoremstyle{plain} \newtheorem{kerthe}[kernels]{Theorem}

\theoremstyle{remark} \newtheorem{mokerex}[kernels]{Examples}

\theoremstyle{plain} \newtheorem{krecom}[kernels]{Lemma}

\theoremstyle{plain} \newtheorem{recom}[kernels]{Lemma}

\theoremstyle{plain} \newtheorem{eh2}[kernels]{Theorem}

\theoremstyle{remark} \newtheorem{eh2rem}[kernels]{Remarks}

\theoremstyle{plain} \newtheorem{ceunions}[kernels]{Proposition}

\theoremstyle{plain} \newtheorem{schur}{Lemma}[subsection]

\theoremstyle{plain} \newtheorem{expneg}[schur]{Lemma}

\theoremstyle{plain} \newtheorem{schoen}[schur]{Lemma}

\theoremstyle{plain} \newtheorem{pto2}[schur]{Proposition}

\theoremstyle{plain} \newtheorem{geneh2}[schur]{Proposition}

\theoremstyle{definition} \newtheorem{mazur}[schur]{Definition}

\theoremstyle{plain} \newtheorem{robgeneh2}[schur]{Proposition}

\theoremstyle{plain} \newtheorem{embcor}[schur]{Corollary}

\theoremstyle{plain} \newtheorem{embthe}[schur]{Theorem}

\theoremstyle{definition} \newtheorem{finprop}{Definition}[subsection]

\theoremstyle{definition} \newtheorem{roealg}[finprop]{Definition}

\theoremstyle{remark} \newtheorem{roealgrem}[finprop]{Remark}

\theoremstyle{definition} \newtheorem{leftrep}[finprop]{Definition}

\theoremstyle{remark} \newtheorem{leftreprem}[finprop]{Remark}

\theoremstyle{definition} \newtheorem{gpalg}[finprop]{Definition}

\theoremstyle{definition} \newtheorem{redgp}[finprop]{Definition}

\theoremstyle{definition} \newtheorem{posfun}[finprop]{Definition}

\theoremstyle{plain} \newtheorem{cstarfacts}[finprop]{Lemma}

\theoremstyle{definition} \newtheorem{cdd}[finprop]{Definition}

\theoremstyle{plain} \newtheorem{cstarfacts2}[finprop]{Lemma}

\theoremstyle{plain} \newtheorem{redroe}[finprop]{Lemma}

\theoremstyle{definition} \newtheorem{c*algdefns}{Definitions}[subsection]

\theoremstyle{plain} \newtheorem{ucplem}[c*algdefns]{Lemma}

\theoremstyle{plain} \newtheorem{stinespring}[c*algdefns]{Theorem}

\theoremstyle{plain} \newtheorem{stinec}[c*algdefns]{Corollary}

\theoremstyle{plain} \newtheorem{ucplem2}[c*algdefns]{Lemma}

\theoremstyle{plain} \newtheorem{ucplem3}[c*algdefns]{Lemma}

\theoremstyle{definition} \newtheorem{exact}[c*algdefns]{Definition}

\theoremstyle{plain} \newtheorem{nuclabel}[c*algdefns]{Proposition}

\theoremstyle{plain} \newtheorem{nuclexact}[c*algdefns]{Theorem}

\theoremstyle{plain} \newtheorem{roenucl}[c*algdefns]{Proposition}

\theoremstyle{plain} \newtheorem{diseq}{Lemma}[subsection]

\theoremstyle{definition} \newtheorem{nondbg}{Definition}[subsection]

\theoremstyle{definition} \newtheorem{nonda}[nondbg]{Definition}

\theoremstyle{remark} \newtheorem{nondarem}[nondbg]{Remarks}

\theoremstyle{plain} \newtheorem{nondacda}[nondbg]{Lemma}

\theoremstyle{plain} \newtheorem{cinvcor}[nondbg]{Corollaries}

\theoremstyle{plain} \newtheorem{nonda2}[nondbg]{Proposition}

\theoremstyle{definition} \newtheorem{amenable}{Definition (and theorem)}[subsection]

\theoremstyle{plain} \newtheorem{amenimpa}[amenable]{Lemma}

\theoremstyle{plain} \newtheorem{kerfun}[amenable]{Lemma}

\theoremstyle{definition} \newtheorem{atmen}[amenable]{Definition}

\theoremstyle{plain} \newtheorem{nuclamen}{Theorem}[subsection]

\theoremstyle{definition} \newtheorem{amenact}[nuclamen]{Definition}

\theoremstyle{definition} \newtheorem{stcech}[nuclamen]{Definition (and theorem)}

\theoremstyle{plain} \newtheorem{amenact1}[nuclamen]{Lemma}

\theoremstyle{plain} \newtheorem{amenact2}[nuclamen]{Lemma}

\theoremstyle{plain} \newtheorem{amenact3}[nuclamen]{Lemma}

\theoremstyle{plain} \newtheorem{amenacta}[nuclamen]{Theorem}

\theoremstyle{definition} \newtheorem{finquot}{Definition}[subsection]

\theoremstyle{definition} \newtheorem{boax}[finquot]{Definition}

\theoremstyle{plain} \newtheorem{boxamen}[finquot]{Proposition}

\theoremstyle{plain} \newtheorem{boxatmen}[finquot]{Proposition}

\theoremstyle{definition} \newtheorem{transl}[finquot]{Definition}

\theoremstyle{definition} \newtheorem{warpm}[finquot]{Definition}

\theoremstyle{plain} \newtheorem{wmet}[finquot]{Proposition}

\theoremstyle{plain} \newtheorem{wcbgp}[finquot]{Proposition}

\theoremstyle{plain} \newtheorem{wamenimpa}[finquot]{Proposition}

\theoremstyle{definition} \newtheorem{cone}[finquot]{Definition}

\theoremstyle{definition} \newtheorem{warpc}[finquot]{Definition}

\theoremstyle{plain} \newtheorem{wcamenimpa}[finquot]{Corollary}

\theoremstyle{plain} \newtheorem{wcaimpamen}[finquot]{Proposition}

\theoremstyle{definition} \newtheorem{adjlap}{Definition}[subsection]

\theoremstyle{plain} \newtheorem{laplac}[adjlap]{Lemma}

\theoremstyle{definition} \newtheorem{expander}[adjlap]{Definition}

\theoremstyle{definition} \newtheorem{exp2}[adjlap]{Definition}

\theoremstyle{plain} \newtheorem{expembed}[adjlap]{Proposition}

\theoremstyle{definition} \newtheorem{propT}{Definition}[subsection]

\theoremstyle{plain} \newtheorem{replem}[propT]{Lemma}

\theoremstyle{plain} \newtheorem{atment}[propT]{Lemma}

\theoremstyle{plain} \newtheorem{proptexp}[propT]{Proposition}

\theoremstyle{definition} \newtheorem{quanta}{Definitions}[subsection]

\theoremstyle{plain} \newtheorem{qaprops}[quanta]{Proposition}

\theoremstyle{plain} \newtheorem{aeqamen}[quanta]{Theorem}

\theoremstyle{plain} \newtheorem{suppinf}[quanta]{Theorem}

\theoremstyle{plain} \newtheorem{cenota}[quanta]{Proposition}

\title{Some notes on Property A}

\author{Rufus Willett}

\begin{document}

\maketitle

\section*{Introduction}

The \emph{coarse Baum-Connes} conjecture states that a certain \emph{coarse assembly 
map}
\begin{displaymath}
\mu:KX_*(X)\to{}K_*(C^*(X))
\end{displaymath}
is an isomorphism (see for example \cite{roe3} or the piece by N. Higson and J. Roe in \cite{frr}).  
It has many important consequences including one of the main motivations for this piece:
a descent technique connects it to injectivity of the (`ordinary') Baum-Connes
assembly map
\begin{displaymath}
\mu:K_*(\underline{E}G)\to{}K_*(C^*_r(G))
\end{displaymath}
when $X$ is the underlying coarse space of a countable discrete group $G$.  Injectivity of
this map is
called the \emph{strong Novikov conjecture} and (as one might guess) implies the 
`usual' Novikov conjecture on higher signatures in the case $G$ is the fundamental group of a
compact manifold.

Our focus in these notes is on property A, which was introduced by G. Yu in \cite{y1} 
as a weak form of amenability. He proved 
that a property A metric space $X$ \emph{coarsely embeds} into (looks
like a subspace of, on the large scale) Hilbert space, that 
such coarse embeddability implies coarse Baum-Connes, and from here the strong
Novikov conjecture if $X$ is the underlying coarse space of a
countable discrete group $G$ (two finiteness restrictions were later removed
in \cite{hig} by Higson, and in more generality in \cite{sty} by G. Skandalis, J-L. Tu and Yu).

The above is background and motivation.  We will focus on rather less complicated
aspects of property A, which is a relatively elementary `large-scale' (or \emph{coarse}) invariant
of metric spaces.  It is part of the appeal of A that such a seemingly
simple property really does have a great many connections with, and
implications for, other areas of mathematics.  These include 
amenability, group representations and various coarse
embeddability problems, all of which this piece starts to explore.  The classical theory of 
\emph{kernels} touches most of these areas, and helps to unify a lot of the material.

The two major omissions we make are the relatively sophisticated K-theoretic consequences
mentioned above, and some finer invariants associated with property A.  The latter unfortunately meant
leaving out many of the results showing that particular classes of groups have property A, as well
as certain \emph{quasi-isometry} invariants, for example those studied in \cite{gk4}, \cite{no4} and \cite{tes},
all of which is interesting and reasonably accessible; references are included in the bibliography.
Our approach does have the benefit that almost everything discussed in these notes is preserved up to the very
flexible standard of \emph{coarse equivalence}.

We assume no prior knowledge of coarse geometry, introducing
ideas as they become necessary.  Coarse geometric language is limited to metric spaces;
introducing abstract coarse theory would unfortunately have taken us a little too far
afield.
The idea (...) is that the piece should be relatively self-contained and
accessible to starting-level graduate students (apart possibly from chapter \ref{analysis}); 
the prerequisites are confined to a
small amount of functional analysis and operator algebra theory.  Many calculations are
gone over in some detail, with an audience new to the common techniques in
this area in mind.  Apologies to those readers who find the pace
slow as a result.  Virtually nothing is new; our aim was only to
survey the area, and present some results in an accessible manner.
Nowhere does a lack of a reference imply any originality whatsoever! 

Despite the attempt to keep the exposition simple, we hope that the collection of 
all these results in one place will be of some use to experts in the area.
Also, we conclude with an
annotated bibliography covering many A-related topics, some, as mentioned above, outside the scope of the main piece.\\

The following is a rough synopsis; see the start of each chapter for a more detailed
outline.\\
Chapter \ref{definition} introduces property A and some basic notions of 
coarse geometry, then gives a selection of equivalent definitions.  It also
reproduces Yu's original proof of coarse embeddability in Hilbert space from property A.\\
Chapter \ref{permprop} provides our first examples of property A spaces
and proves some basic permanence properties.  It then introduces the underlying coarse space of a 
countable discrete group, and proves that the class of property A groups is
closed under certain group-theoretic operations.  Such groups provide our most important
examples of property A spaces.\\
Chapter \ref{kerembed} introduces the theory of kernels, and gives some
characterisations of coarse embeddability into Hilbert space.  It
continues with some more of the classical theory of kernels, and a discussion
of coarse embeddability into $l^p$ spaces.\\
Chapter \ref{analysis} starts with a brief introduction to the uniform
Roe algebra of a (uniformly discrete) metric space, and the reduced $C^*$-algebra of a
group.  It then proves the equivalence of property A with nuclearity of the
former and exactness of the latter.\\
Chapter \ref{beyond} discusses the `correct' definition of property A for non-discrete
spaces.\\
Chapter \ref{amenable} introduces the relationship between
property A and amenability (A is often called a `weak form of amenability').  It 
continues with some concrete constructions that further highlight this
relationship.\\  
Chapter \ref{nona} studies some spaces that do not have property A (despite
the impressive consequences of A, finding such examples is actually quite hard).  This is one
of the most active areas of research directly connected to property A.\\
Chapter \ref{bib} contains (very brief) comments on the items in the 
bibliography; some of these are not used in the main discussion, but are included 
for completeness.

\subsection*{Notation}
If $\xi$ is a map from a set $X$ to some space of functions
we often write $\xi_x$ for the image of $x$ under $\xi$ (the aim is to make
speaking of values $\xi_x(y)$ of $\xi_x$ clearer).\\
$\delta_x$, on the other hand, is reserved for the characteristic function of $x$.\\
We use $e$ for the trivial element of a group, while $1$ is the multiplicative identity
in an algebra, or the trivial group.
$\mathbf{1}$ is reserved for the constant map of value one.\\
$f|_S$ denotes the restriction of the map $f$ to some subset $S$ of its domain.\\
$B_d(x,R)$ and $B_X(x,R)$ both denote the open ball of radius $R$ about $x$ in a metric space $(X,d)$, i.e.$ \{y\in{}X:d(x,y)<R\}$.  We
usually omit the subscripts, and also write $\bar{B}(x,R)$ for  the closed ball $\{y\in{}X:d(x,y)\leq{}R\}$.\\
If $X$ is a normed space, $S(X)$ denotes its unit sphere, i.e. $\{x\in{}X:\|x\|=1\}$.

Throughout, a reference of the form \ref{permprop} is to the second chapter, one of the form \ref{groups} is to the third
section of chapter \ref{permprop}, and one of the form \ref{dlim} to the fifth definition, remark or result in
section \ref{groups}.

\subsection*{Acknowledgments}
I would particularly like to thank Erik Guentner for many illuminating discussions, and his
hospitality, during a visit to the university of Hawaii.  The anonymous referee also made 
several helpful comments, for which I am grateful.  Finally, I would like to thank my advisor John Roe for continuing support,
inspiration and prodding me in this interesting (\emph{I} think; the reader is left to judge for him or herself :) ) direction in the first place.

\pagebreak

\tableofcontents

\pagebreak

\section{Property A and coarse geometry} \label{definition} 

\subsection{Introduction}\label{defintro}

\begin{propa}[Property A, \cite{y1}] \label{propa} 
A discrete metric space $(X, d)$ has \emph{property A} if for all
$R, \epsilon>0$, there exists a family $\{A_x\}_{x\in{}X}$ of finite, non-empty
subsets of $X\times\mathbb{N}$ such that:
\begin{itemize}
\item for all $x,y\in{}X$ with $d(x,y)\leq{}R$ we have 
$\frac{|A_x\triangle{}A_y|}{|A_x\cap{}A_y|}<\epsilon$
\item there exists $S$ such that for each $x\in{}X$, if $(y,n)\in{}A_x$,
 then $d(x,y)\leq{}S$
\end{itemize}
\end{propa}

Note that, following previous authors, we have dropped the condition
`$(x,1)\in{}A_x$ for all $x\in{}X$', which makes the definition a little easier
to work with.
Our first goal is to give a precise meaning to the statement `property A
is a large-scale invariant'.  We define a class of maps
that preserve large-scale, or \emph{coarse}, structure.  See the book 
\cite{roe1} for a self-contained introduction to coarse geometry.

\begin{coarse}
\label{coarse}
If $f:X\to{}Y$ is a map of metric spaces, it is said to be:
\begin{itemize}
\item \emph{bornologous} if for all $R>0$ there exists $S>0$ such that
$d(x,x')<R$ implies $d(f(x), f(x'))<S$ 
(`points that start close cannot end up too far apart').
\item \emph{proper} if for each $x\in{}X$, and each $R>0$ there exists $S>0$ 
such that $f^{-1}(B(f(x), R))\subseteq{}B(x,S)$
(`points that end up close cannot have started out too far apart').
\end{itemize}
If $f$ satisfies both of these conditions, it is called a \emph{coarse map}.\\
Two coarse maps $f,g:X\to{}Y$ are \emph{close} if $\{d(f(x),g(x)):x\in{}X\}$
is a bounded set.\\
$X$ and $Y$ are \emph{coarsely equivalent} if there exist coarse maps 
$f:X\to{}Y$ and $g:Y\to{}X$ such that $f\circ{}g$ and $g\circ{}f$ are close to
the identities on $X$ and $Y$ respectively.  Either of the maps $f$, $g$ is
then called a \emph{coarse equivalence}.
\end{coarse}

Note that if the metric space $X$ is such that closed bounded sets are compact, then
`proper' in the sense above is the same thing as proper in the sense of classical topology - 
compact sets are pulled back to compact sets.

\begin{coarseinv} \label{cinv}
Property A is invariant under coarse equivalence.
\end{coarseinv}

\begin{proof}
Say $X$ and $Y$ are coarsely equivalent metric spaces, so there exist coarse
maps $f:X\to{}Y$ and $g:Y\to{}X$ satisfying the properties above.\\
Assume that $X$ has property A, and let $R,\epsilon>0$.\\
As $g$ is bornologous there exists $R'>0$ such that if $d(y,y')<R$, then
$d(g(y),g(y'))<R'$.  Choose a family of subsets $\{A_x\}_{x\in{}X}$ of
$X\times\mathbb{N}$ that satisfies the conditions in the definition of property
A for $R'$, $\epsilon$ as above, and for some $S$.\\
Let $y_0$ be any element of $Y$. For each $y\in{}Y$ let
$n_y=|(f^{-1}(y)\times\mathbb{N})\cap{}A_{g(y_0)}|$ and define
\begin{displaymath}
B_{y_0}=\bigcup_{y\in{}Y}\{(y,1),(y,2),...,(y,n_y)\}.
\end{displaymath}
Note that the sets $\{(f^{-1}(y)\times\mathbb{N})\cap{}A_{g(y_0)}\}$ partition
$A_{g(y_0)}$, so that $|B_{y_0}|=|A_{g(y_0)}|$ is finite.
Extend this construction to build a family $\{B_y\}_{y\in{}Y}$.\\
Then for any $y,y'\in{}Y$ we have:
\begin{displaymath}
|B_y\cap{}B_{y'}|\geq{}|A_{g(y)}\cap{}A_{g(y')}|~\textrm{and}
\end{displaymath}\begin{displaymath}
|B_y\triangle{}B_{y'}|\leq{}|A_{g(y)}\triangle{}A_{g(y')}|,
\end{displaymath}
so that if $d(y,y')<R$, then $d(g(y),g(y'))<R'$ and
\begin{displaymath}
\frac{|B_y\triangle{}B_{y'}|}{|B_y\cap{}B_{y'}|}\leq
\frac{|A_{g(y)}\triangle{}A_{g(y')}|}{|A_{g(y)}\cap{}A_{g(y')}|}<\epsilon.
\end{displaymath}
Moreover, if $(y',n)\in{}B_y$,
$(f^{-1}(y')\times\mathbb{N})\cup{}A_{g(y)}\not=\emptyset$; 
let $(x,n)$ be in this set.\\
Then by assumption on the family $\{A_x\}_{x\in{}X}$, 
$d(g(y),x)\leq{}S$, whence $d(fg(y),y')\leq{}S'$ for some $S'$ as $f$ is 
bornologous.
Finally, $d(y,y')\leq{}S''$ for some $S''$, as $f\circ{}g$ is close to the 
identity on $Y$.\\
The family $\{B_y\}_{y\in{}Y}$ thus has the required properties.
\end{proof}

\begin{cinvrem} \label{cinvrem}
In fact, all we actually needed in the above proof was the existence of a
bornologous map $f:X\to{}Y$ such that $f\circ{}g$ is close to the identity
on $Y$.  This follows simply from the fact that $g$ is a \emph{coarse embedding} 
as defined below. Hence the above proposition 
actually proves that if $X$ has property A and there exists a coarse embedding 
of $Y$ into $X$, then $Y$ has property A too; in particular, property A is 
inherited by subspaces.
\end{cinvrem}

As already remarked, one of Yu's main motivations for introducing property A
was that it implies coarse embeddability into Hilbert space, which in turn
implies the strong Novikov conjecture.  A coarse embedding is simply a
coarse equivalence with its image; the following gives a formal definition.

\begin{coarseembed} \label{cembed}
A map $f:X\to{}Y$ of metric spaces is \emph{effectively proper} if for all
$R>0$ there exists $S>0$ such that for all $x\in{}X$, 
$f^{-1}(B(f(x), R))\subseteq{}B(x,S)$.\\
A \emph{coarse embedding} is an effectively proper, bornologous map.\\
This is (clearly) equivalent to the existence of non-decreasing maps
$\rho_1,\rho_2:\mathbb{R}^{+}\to\mathbb{R}^+$ such that $\rho_1(t)\to+\infty$ as
$t\to+\infty$ and
\begin{displaymath}
\rho_1(d(x,y))\leq{}d(f(x),f(y))\leq\rho_2(d(x,y))
\end{displaymath}
for all $x,y\in{}X$.
\end{coarseembed}

\begin{cerem} \label{cerem}
Coarse embeddings can be considered as `coarsely injective' maps.
Correspondingly, $f:X\to{}Y$ is `coarsely surjective'
if for some $R>0$ the collection $\{B(y,R):y\in{}f(X)\}$ covers $Y$ (we also
say that $f(X)$ is \emph{coarsely dense} in $Y$).
  
If $f$ is a coarse equivalence as introduced in definition 
\ref{coarse}, then it is coarsely injective and surjective in the
senses introduced above.  The inverse up to closeness required by 
\ref{coarse} is also `coarsely bijective' in this sense (the proof of these
comments is an easy exercise).  Conversely, a coarsely bijective map
gives rise to a coarse equivalence.  These remarks are all easy to prove, and make good exercises 
for those new to this area.

It can in fact be proved (we will not use this formalism, but hope it clarifies some of the above for 
those used to the terminology)
that there is a category whose objects are all (or some appropriate subclass of all) 
metric spaces, and whose morphisms are
equivalence classes of coarse maps under the equivalence relation of closeness.  Isomorphism,
injectivity, and surjectivity in this category correspond to the notions above.

Note that in the literature what we have called a coarse embedding is often
referred to as an \emph{uniform embedding} (after M. Gromov, see ~\cite{gr1}).
We use the term coarse embedding, as uniform embedding means something 
different to Banach space geometers (see the book \cite{bl}).
\end{cerem}

The next result is Yu's original proof that property A implies coarse
embeddability into Hilbert space.  It is redundant for our exposition, as the characterisations of 
coarse embeddability in \ref{eh2} reprove it.  Nonetheless, we include it
as it provides a blueprint for other arguments
(e.g. \ref{geneh2}), and is a nice example of a direct coarse
geometric proof (it is actually a generalisation of an `equivariant' argument, due to
Bekka, Cherix and Valette, that can be found in \cite{frr}).  From now on, we will often say 
`coarsely embeddable' to mean
`coarsely embeddable in Hilbert space'.

\begin{embedhilbert}[\cite{y1}] \label{yembed}
A discrete metric space $X$ with property A coarsely embeds into Hilbert space.
\end{embedhilbert}

\begin{proof}
Start by defining a Hilbert space $\mathcal{H}$ by
\begin{displaymath}
\mathcal{H}=\bigoplus_{k=1}^\infty{}l^2(X\times\mathbb{N})
\end{displaymath}
Using property A for $X$, define a sequence of families of subsets of $X$,
$\{A_x^k\}_{x\in{}X}$, $k=1,2,3,...$ that satisfy
\begin{itemize}
\item There exist $S_k>0$ such that if $(x_1,n),(x_2,m)\in{}A_x^k$ for some
$m,n,k$, then $d(x_1,x_2)\leq{}S_k$.
\item 
$\frac{|A_{x_1}^k\triangle{}A_{x_2}^k|}{|A_{x_1}^k\cap{}A_{x_2}^k|}
<\frac{1}{2^{2k+1}}$
whenever $d(x_1,x_2)\leq{}k$.
\end{itemize}
Note first that
\begin{align*}
& \left\|\frac{\chi_{A_{x_1}^k}}{|A_{x_1}^k|^\frac{1}{2}}
- \frac{\chi_{A_{x_2}^k}}{|A_{x_2}^k|^\frac{1}{2}}
\right\|_{l^2(X\times\mathbb{N})}^2 \\
& =\sum_{x\in{}A_{x_1}^k\backslash{}A_{x_2}^k}\frac{1}{|A_{x_1}^k|}
+\sum_{x\in{}A_{x_2}^k\backslash{}A_{x_1}^k}\frac{1}{|A_{x_2}^k|}
+\sum_{x\in{}A_{x_1}^k\cap{}A_{x_2}^k}
\left(\frac{1}{|A_{x_1}^k|^{\frac{1}{2}}}-\frac{1}{|A_{x_1}^k|^{\frac{1}{2}}}\right) \\
& \leq\frac{|A_{x_1}^k\triangle{}A_{x_2}^k|}{|A_{x_1}^k\cap{}A_{x_2}^k|}
+|A_{x_1}^k\cap{}A_{x_2}^k|\left(\frac{|A_{x_1}^k|+|A_{x_2}^k|-
2|A_{x_1}^k|^{\frac{1}{2}}|A_{x_2}^k|^{\frac{1}{2}}}{|A_{x_1}^k||A_{x_2}^k|}\right) \\
& \leq\frac{|A_{x_1}^k\triangle{}A_{x_2}^k|}{|A_{x_1}^k\cap{}A_{x_2}^k|}
+\frac{|A_{x_1}^k|+|A_{x_2}^k|-
2|A_{x_1}^k\cap{}A_{x_2}^k|}{|A_{x_1}^k\cap{}A_{x_2}^k|} \\
& \leq2\frac{|A_{x_1}^k\triangle{}A_{x_2}^k|}{|A_{x_1}^k\cap{}A_{x_2}^k|}
<\frac{1}{2^{2k}}.
\end{align*}
Hence \begin{displaymath}
\left\|\frac{\chi_{A_{x_1}^k}}{|A_{x_1}^k|^\frac{1}{2}}
-\frac{\chi_{A_{x_2}^k}}{|A_{x_2}^k|^\frac{1}{2}}
\right\|_{l^2(X\times\mathbb{N})}<\frac{1}{2^k}.
\end{displaymath}
Now, choose $x_0\in{}X$, and define $f:X\to{}\mathcal{H}$ by
\begin{displaymath}
f(x)=\bigoplus_{k=1}^\infty
\left(\frac{\chi_{A_x^k}}{|A_x^k|^\frac{1}{2}}
-\frac{\chi_{A_{x_0}^k}}{|A_{x_0}^k|^\frac{1}{2}}\right).
\end{displaymath}
If $m\leq{}d(x_1,x_2)<m+1$ we can make the estimate 
\begin{align*}
\|f(x_1)-f(x_2)\|^2
& =\sum_{k=1}^\infty\left\|\frac{\chi_{A_{x_1}^k}}{|A_{x_1}^k|^\frac{1}{2}}
-\frac{\chi_{A_{x_2}^k}}{|A_{x_2}^k|^\frac{1}{2}}\right\|
_{l^2(X\times\mathbb{N})}^2 \\
& \leq\sum_{k=1}^m\left\|\frac{\chi_{A_{x_1}^k}}{|A_{x_1}^k|^\frac{1}{2}}
-\frac{\chi_{A_{x_2}^k}}{|A_{x_2}^k|^\frac{1}{2}}\right\|
_{l^2(X\times\mathbb{N})}^2+1 \\
& \leq\sum_{k=1}^m\left(\left\|\frac{\chi_{A_{x_1}^k}}{|A_{x_1}^k|^\frac{1}{2}}\right\|
_{l^2(X\times\mathbb{N})}
+\left\|\frac{\chi_{A_{x_2}^k}}{|A_{x_2}^k|^\frac{1}{2}}\right\|
_{l^2(X\times\mathbb{N})}\right)^2+1 \\
& =(2m)^2+1.
\end{align*}
Hence $\|f(x_1)-f(x_2)\|\leq2m+1$.\\
On the other hand, from the choice of the families $\{A_x^k\}_{x\in{}X}$ it
follows that $S_k\geq{}k$, whence
$Q_t:=|k\in\mathbb{N}:S_k<{}y|$ exists for all $t\in\mathbb{R}^+$.\\
Let $x_1,x_2\in{}X$, and say $Q_{d(x_1,x_2)}=m$.  Then
\begin{displaymath}
\|f(x_1)-f(x_2)\|^2
=\sum_{k=1}^\infty\left\|\frac{\chi_{A_{x_1}^k}}{|A_{x_1}^k|^\frac{1}{2}}
-\frac{\chi_{A_{x_2}^k}}{|A_{x_2}^k|^\frac{1}{2}}\right\|
_{l^2(X\times\mathbb{N})}^2
\geq{}4m,
\end{displaymath}
as each term where $d(x_1,x_2)>S_k$ contributes $4$ to the total sum, as
if $d(x_1,x_2)<S_k$, $A_{x_1}^k\cap{}A_{x_2}^k=\emptyset$.\\
It now follows that the maps 
$\rho_1,\rho_2:\mathbb{R}^+\to\mathbb{R}^+$ defined
by $\rho_1(t)=2\sqrt{Q_t}$, and $\rho_2(t)=2t+1$ are non-decreasing, 
$\rho_1(t)\to+\infty$ as $t\to+\infty$, and
\begin{displaymath}
\rho_1(d(x_1,x_2))\leq{}d(f(x_1),f(x_2))\leq\rho_2(d(x_1,x_2)),
\end{displaymath}
i.e. \ref{cembed} is satisfied.
\end{proof}

Note that if $X$ is a separable space, then it is coarsely equivalent to a countable subspace
(cf. \ref{diseq}), which we can pass to.  It follows that if $f:X\to\mathcal{H}$ is a coarse 
embedding, then the closure of span$f(X)$ is a separable Hilbert space into which $X$ coarsely embeds.

\subsection{Equivalent definitions}

The one (long) theorem in this section gives some equivalent formulations of 
property A.  In order to prove some of these equivalences, we need to restrict to the following class 
of spaces. 

\begin{bg} \label{bg}
A discrete metric space $X$ has \emph{bounded geometry} if for all $r>0$
there exists $N_r$ such that $|B(x,r)|<N_r$ for all $x\in{}X$.
\end{bg} 

This class includes many interesting examples, in particular, all countable discrete groups (see section \ref{groups}).

We include all of the
reformulations of property A in one place for ease of reference, but have not introduced the
concepts needed to understand them all, believing that it would be better to introduce
some examples first (in the next chapter).  The definitions of kernels and some of their
properties are given in \ref{kernels}, and that of the uniform Roe algebra in
\ref{roealg}; we will not use either (8) or (9) below until after these two points.

(1) $\Rightarrow$ (2) and (3) $\Rightarrow$ (1) come from \cite{hr}.  The
equivalences of (2), (3) (partially) and (8) are from \cite{tu}.  The full
version of (3) is from \cite{dr}.  (6) is from \cite{dg2}, while the proof of (6) $\Rightarrow$ (2) that
does not require bounded geometry is due to N. Wright. 
The collection of the most of the conditions together, and the remaining equivalences, are from \cite{bnw}.
We make the following two definitions: the first follows \cite{bnw} and makes the statement a little more
concise; the second states exactly what we mean by a partition of unity in (6) (note in particular that 
there are no topological restrictions on our partitions of unity).

\begin{equivdef} \label{equivdef}
Throughout, a function $x\mapsto{}\xi_x$ from $X$ to a Banach space will be said 
to have $(R,\epsilon)$ \emph{variation} if $d(x,y)\leq{}R$ implies 
$\|\xi_x-\xi_y\|<\epsilon$.
\end{equivdef}

\begin{partition}
Let $X$ be a metric space and let $\mathcal{U}=\{U_i\}_{i\in{}I}$ be a cover (i.e. for all $x\in{}X$
there exists $i$ such that $x\in{}U_i$).\\
A \emph{partition of unity on $X$ subordinate to $\mathcal{U}$} is a collection $\{\phi_i\}_{i\in{}I}$
of functions $\phi_i:X\to[0,1]$ such that $\phi_i$ is supported in $U_i$ for all $i$, and $\sum_i\phi(x)=1$ for
all $x\in{}X$.\\
The \emph{diameter} of a subset $U$ of $X$ is diam$(U)=\sup\{d(x,y):x,y\in{}U\}$ (possibly infinite).
\end{partition}

\begin{equivalences}
\label{equiv}
Let $X$ be a bounded geometry discrete metric space. 
The following are equivalent:
\begin{enumerate}
\item $X$ has property A.

\item There exists $p$ with $1\leq{}p<\infty$ such that for all $R,\epsilon>0$ there exists a map $\xi:X\to{}l^1(X)$ such that:
\begin{enumerate} 
\item $\|\xi_x\|_p=1$ for all $x\in{}X$; 
\item $\xi$ has $(R, \epsilon)$ variation; 
\item there exists $S>0$ such that for each $x$, $\xi_x$ is supported in 
$\bar{B}(x;S)$.
\end{enumerate}

\item For all $p$ with $1\leq{}p<\infty$ and all
$R,\epsilon>0$ there exists a map $\eta:X\to{}l^p(X)$ such that:
\begin{enumerate} 
\item $\|\eta_x\|_p=1$ for all $x\in{}X$; 
\item $\eta$ has $(R, \epsilon)$ variation; 
\item there exists $S>0$ such that for each $x$, $\eta_x$ is supported in 
$\bar{B}(x;S)$.
\end{enumerate}

\item For all $p$ with $1\leq{}p<\infty$ and all
$R,\epsilon>0$ there exists a map $\eta:X\to{}l^p(X)$ such that:
\begin{enumerate} 
\item $\|\eta_x\|_p=1$ for all $x\in{}X$; 
\item $\eta$ has $(R, \epsilon)$ variation; 
\item for all $\delta>0$ there exists $S>0$ such that for all $x\in{X}$,\\
$\|\eta_x|_{\bar{B}(x,S)}\|_p>1-\delta$.
\end{enumerate}

\item There exists $\delta<1$ such that for all $p$ with $1\leq{}p<\infty$ and all
$R,\epsilon>0$ there exists a map $\eta:X\to{}l^p(X)$ such that:
\begin{enumerate} 
\item $\|\eta_x\|_p=1$ for all $x\in{}X$; 
\item $\eta$ has $(R, \epsilon)$ variation; 
\item there exists $S>0$ such that for each $x\in{}X$, $\|\zeta_x|_{\bar{B}(x,S)}\|_p>1-\delta$ and
$\|\zeta_x|_{\bar{B}(x,R+S)\backslash\bar{B}(x,S)}\|_p<\epsilon$.
\end{enumerate}

\item For all $R, \epsilon>0$ there exists a cover $\mathcal{U}=\{U_i\}_{i\in{}I}$ of
$X$, a 
partition of unity $\{\phi_i\}_{i\in{}I}$ subordinated to $\mathcal{U}$ and $S>0$ such that:
\begin{enumerate}
\item if $d(x,y)\leq{}R$, then $\sum_{i\in{}I}|\phi_i(x)-\phi_i(y)|<\epsilon$;
\item \emph{diam}$(U_i)\leq{}S$ for all $i$.
\end{enumerate}

\item For all $R,\epsilon>0$ there exists a Hilbert space $\mathcal{H}$ and a map
$f:X\to{}\mathcal{H}$ such that:
\begin{enumerate} 
\item $\|f(x)\|_{\mathcal{H}}=1$ for all $x\in{}X$; 
\item $f$ has $(R, \epsilon)$ variation; 
\item there exists $S>0$ such that $d(x,y)>S$ implies that 
$\langle{}f(x),f(y)\rangle_{\mathcal{H}}=0$.
\end{enumerate} 

\item For all $R, \epsilon>0$ there exists a normalised, finite propagation,
symmetric, positive type kernel $k:X\times{}X\to\mathbb{R}$ with 
$(R,\epsilon)$ variation.

\item For all $R, \epsilon>0$ there exists a finite propagation, positive
type kernel\\ $k:X\times{}X\to\mathbb{C}$ with $(R,\epsilon)$ variation and such
that convolution with $k$ defines a bounded operator in the uniform Roe algebra
$C^{*}_{u}(X)$.

\end{enumerate}
\end{equivalences}

The conditions can be considered as being split into two rough parts:
\begin{itemize}
\item The original definition, (1), and maps with $(R,\epsilon)$ variation 
and `bounded' support, (2), (3), (4),
(5) and (6).  Although (6) looks somewhat different to the others, we hope the proof
makes its similarity apparent (having said this, it sometimes seems easier to use;
see section \ref{genpermprop}).  We prove (1) $\Rightarrow$ (2) $\Rightarrow$ (3) $\Rightarrow$ (4) 
$\Rightarrow$ (5) $\Rightarrow$ (2),  (3) $\Rightarrow$ (1),  (3) $\Rightarrow$ (6) and (6) $\Rightarrow$ (2). 
\item Kernels, (7), (8), (9); the sense in which (7) is a `kernel 
condition' will (we hope) become apparent from theorem \ref{kerthe}.  
We prove (3) $\Rightarrow$ (7)
$\Rightarrow$ (8) $\Rightarrow$ (9) $\Rightarrow$ (2).
\end{itemize}

It is worth pointing out that of all these implications, only
(8) $\Rightarrow$ (9) and (2) $\Rightarrow$ (1) require the bounded geometry
assumption.

We have chosen to reproduce most calculations in a fair amount of detail, 
as a lot of the techniques used are very common in this
area (having said that, there is little of a `conceptual' nature in the
proof). 

\begin{proof}
(1) implies (2):\\
Let $R,\epsilon>0$ and let $\{A_x\}_{x\in{}X}$ be a family of subsets of 
$X\times\mathbb{N}$ as in (1).\\
For each $x\in{}X$, define a map $\xi_x:X\to\mathbb{C}$ by
\begin{displaymath}
x'\mapsto\frac{|(\{x'\}\times\mathbb{N})\cap{}A_x|}{|A_x|}.
\end{displaymath}
Note that
\begin{displaymath}
\|\xi_x\|_1=\sum_{x'\in{}X}\frac{|(\{x'\}\times\mathbb{N})\cap{}A_x|}{|A_x|}=1,
\end{displaymath}
so $\xi:X\to{}l^1(X)$ is well defined and satisfies $\|\xi_x\|_1=1$ 
for all $x\in{}X$.\\
Now, say $d(x_1,x_2)\leq{}R$.  Then firstly
\begin{align*}
\|\xi_{x_1}.|A_{x_1}|-\xi_{x_2}.|A_{x_2}|\|_1 &
\leq\sum_{x\in{}X}|(\{x\}\times\mathbb{N})\cap{}A_{x_1}|-
|(\{x\}\times\mathbb{N})\cap{}A_{x_2}| \\
& \leq|A_{x_1}\triangle{}A_{x_2}|.
\end{align*}
Secondly
\begin{displaymath}
\epsilon>\frac{|A_{x_1}\triangle{}A_{x_2}|}{|A_{x_1}\cap{}A_{x_2}|}=
\frac{|A_{x_1}|+|A_{x_2}|-2|A_{x_1}\cap{}A_{x_2}|}{|A_{x_1}\cap{}A_{x_2}|},
\end{displaymath} whence \begin{displaymath}
2+\epsilon>\frac{(2+\epsilon)|A_{x_1}\cap{}A_{x_2}|}{|A_{x_1}}
>1+\frac{|A_{x_2}|}{|A_{x_1}|},\end{displaymath} so by symmetry
\begin{displaymath}
1+\epsilon>\frac{|A_{x_2}|}{|A_{x_1}|}>\frac{1}
{1+\epsilon}.
\end{displaymath}
Combining these two comments, we conclude that
\begin{align*}
\|\xi_{x_1}-\xi_{x_2}\| & \leq\left\|\xi_{x_1}-\xi_{x_2}.\frac{|A_{x_2}|}{|A_{x_1}|}\right\|_1
+\left\|\xi_{x_2}-\xi_{x_2}.\frac{|A_{x_2}|}{|A_{x_1}|}\right\|_1 \\
& \leq{}\frac{|A_{x_1}\triangle{}A_{x_2}|}{|A_{x_1}|}+
\frac{\epsilon}{1+\epsilon}
\leq{}\frac{|A_{x_1}\triangle{}A_{x_2}|}{|A_{x_1}\cap{}A_{x_2}|}+\epsilon
<2\epsilon.
\end{align*}
Finally, note that if $\xi_x(x')\not=0$, then $(x',n)\in{}A_x$ for some $n$,
whence $d(x,x')\leq{}S$. Hence $\xi_x$ is supported in $\bar{B}(x,S)$.\\

(2) implies (3):\\
Let $R,\epsilon>0$ and let $\xi:X\to{}l^p(X)$ be a map as given in (2).\\
By replacing $\xi_x$ by $x'\mapsto|\xi_x(x')|$, (which still gives a map $\xi$ with the
properties in (2)) assume that each $\xi_x$ positive-valued.\\
Now take any $q$ with $1\leq{}q<\infty$, and define a map $\eta:X\to{}l^q(X)$ by setting 
$\eta_x(x')=\xi_x(x')^{\frac{p}{q}}$.  Note that
\begin{displaymath}
\|\eta_x\|_q=\left(\sum_{x'\in{X}}(\xi_x(x')^{\frac{p}{q}})^q\right)^\frac{1}{q}=1.
\end{displaymath}
Further, there exists $S>0$ such that $\eta_x$ is supported in $\bar{B}(x,S)$,
as the same is true for $\xi$.  Finally, note that we have the inequality
$|a-b|^q\leq|a^q-b^q|$ for any $a,b\geq0$ (use elementary calculus), whence
if $d(x_1,x_2)\leq{}R$,
\begin{align*}
\|\eta_{x_1}-\eta_{x_2}\|_q^q & =\sum_{x\in{}X}|\eta_{x_1}(x)-\eta_{x_2}(x)|^q
\leq\sum_{x\in{}X}|\eta_{x_1}(x)^q-\eta_{x_2}(x)^q| \\
& =\|\xi_{x_1}-\xi_{x_2}\|_p<\epsilon.
\end{align*}

(3) implies (4):\\
A map as in (3) automatically has the property in (4).\\

(4) implies (5):\\
Take $\eta$ as in (4) for fixed constants $R,\epsilon$.  Take $\zeta$ equal
to this $\eta$.  Fix $\delta$ with $0<\delta<1$, and take $\delta'=\min\{\delta,\epsilon^p\}$.\\
By (4), then, there exists $S>0$ such that for all $x\in{}X$,
$\|\zeta|_{\bar{B}(x,S)}\|_p^p>1-\delta'>1-\delta$, whence
$\|\zeta|_{\bar{B}(x,S)}\|_p>1-\delta$.
It follows moreover that
$\|\zeta_x|_{\bar{B}(x,R+S)\backslash\bar{B}(x,S)}\|_p^p<\delta'\leq\epsilon^p$, whence
$\|\zeta_x|_{\bar{B}(x,R+S)\backslash\bar{B}(x,S)}\|_p<\epsilon$.\\

(5) implies (2):\\
Let $R,\epsilon>0$ be fixed, and let $\zeta$ satisfy the conditions in (5) for $p=1$, some fixed $0<\delta<1$ and
parameters $R,\epsilon', S$.
Define $\theta_x=\zeta_x|_{B(x,R+S)}$.  $\theta_{x_1}-\theta_{x_2}$ is then equal to
\begin{align*}
\zeta_{x_1}|_{\bar{B}(x_1,R+S)\backslash\bar{B}(x_2,R+S)}
& -\zeta_{x_2}|_{\bar{B}(x_2,R+S)\backslash\bar{B}(x_1,R+S)} \\
& +(\zeta_{x_1}-\zeta_{x_2})|_{\bar{B}(x_1,R+S)\cap\bar{B}(x_2,R+S)}.
\end{align*}
If $d(x_1,x_2)\leq{}R$, we have 
$\bar{B}(x,R+S)\backslash\bar{B}(y,R+S)\subseteq\bar{B}(x,R+S)
\backslash\bar{B}(y,S)$; the
first term above thus has $l^1$-norm at most $\epsilon'$, and similarly 
for the second.  On the other hand, the fact that $\zeta_x$ has 
$(R,\epsilon)$ variation ensures that the third term is also bounded in $l^1$-norm by $\epsilon$.\\
Hence $\|\theta_{x_1}-\theta_{x_2}\|_1<3\epsilon'$.\\

Define $\xi:X\to{}l^1(X)$ by $\xi_x=\theta_x/\|\theta_x\|_1$.  It clearly satisfies properties (a) and (c) of (2) (with $\xi_x$ supported in $\bar{B}(x,S+R)$).\\
Moreover, for
$d(x_1,x_2)\leq{}R$ we have the estimates
\begin{align*}
\|\xi_{x_1}-\xi_{x_2}\|_1 & =\left\|\frac{\theta_{x_1}}{\|\theta_{x_1}\|_1}+
\frac{\theta_{x_2}}{\|\theta_{x_2}\|_1}-\frac{\theta_{x_1}}{\|\theta_{x_1}\|_1}-
\frac{\theta_{x_2}}{\|\theta_{x_2}\|_1}\right\|_1 \\
& \leq\frac{1}{\|\theta_{x_1}\|_1}\|\theta_{x_1}-\theta_{x_2}\|_1+
\|\theta_{x_2}\|_1\left|\frac{1}{\|\theta_{x_1}\|_1}-\frac{1}{\|\theta_{x_1}\|_1}\right| \\
& =\frac{1}{\|\theta_{x_1}\|_1}\|\theta_{x_1}-\theta_{x_2}\|_1+
\frac{|\|\theta_{x_2}\|_1-\|\theta_{x_1}\|_1|}{\|\theta_{x_1}\|_1}
<\frac{6}{1-\delta}\epsilon'.
\end{align*}
As this expression tends to zero as $\epsilon'$ does, this completes the proof.\\

(3) implies (1):\\
Let $R,\epsilon>0$ and say $\xi:X\to{}l^1(X)$ is as in (3) with respect to $p=1$ and parameters $R,\epsilon'$.\\
In particular, there exists $S>0$ such that
for all $x\in{}X$, Supp$(\xi_x)\subseteq\bar{B}(x,S)$,  whence by bounded 
geometry there exists $N\in\mathbb{N}$ such that 
$|\textrm{Supp}(\xi_x)|<N$ for all $x\in{}X$.\\
By replacing $\xi_x$ by the positive-valued function 
$x'\mapsto|\xi_x(x')|$, assume $\xi_x$
takes only positive values for each $x\in{}X$.\\
Further, let $M\in\mathbb{N}$ satisfy $M>\frac{N}{\epsilon'}$, and define
$\theta_x:x'\mapsto{}j/M$, where $j\in\mathbb{N}$, $0\leq{}j\leq{}M$, and
$j-1<M\xi_x(x')<j$.\\  
Then $\theta_x$ satisfies 
$\|\xi_x-\theta_x\|_1\leq\frac{N}{M}\leq\epsilon'$, whence if $d(x_1,x_2)<R$,
$\|\theta_{x_1}-\theta_{x_2}\|<3\epsilon'$.\\
For each $x$, define $A_x\subseteq{}X\times\mathbb{N}$ by stipulating that
$(x',j)\in{}A_x$ if and only if $j>0$ and $j/M\leq\theta_x(x')$.  Note that
\begin{displaymath}
|A_{x_1}\triangle{}A_{x_2}|=
\sum_{(x,n)\in{}X\times\mathbb{N}}M|\theta_{x_1}(x)-\theta_{x_2}(x)|
=M\|\theta_{x_1}-\theta_{x_2}\|_1=|A_{x_1}|.\|\theta_{x_1}-\theta_{x_2}\|_1,
\end{displaymath}whence
\begin{displaymath}
\frac{|A_{x_1}\triangle{}A_{x_2}|}{|A_{x_1}|}
=\|\theta_{x_1}-\theta_{x_2}\|_1<3\epsilon'.
\end{displaymath}
Now, the left hand side on the line above is equal to:
\begin{displaymath}
\frac{|A_{x_1}|+|A_{x_2}|-2|A_{x_1}\cap{}A_{x_2}|}{|A_{x_1}|}
=\frac{2M-2|A_{x_1}\cap{}A_{x_2}|}{2M}.
\end{displaymath}
Putting this together, we get that
$|A_{x_1}\cap{}A_{x_2}|>M(1-\frac{3}{2}\epsilon')$.  It finally follows that
\begin{displaymath}
\frac{|A_{x_1}\triangle{}A_{x_2}|}{|A_{x_1}\cap{}A_{x_2}|}
=\frac{2M-2|A_{x_1}\cap{}A_{x_2}|}{|A_{x_1}\cap{}A_{x_2}|}
<\frac{2M-2(M(1-\frac{3}{2}\epsilon'))}{M(1-\frac{3}{2}\epsilon')}
=\frac{3\epsilon'}{1-\frac{3}{2}\epsilon'}.
\end{displaymath}
As the expression on the right tends to zero as $\epsilon'$ does, this gives
the result.\\

(3) implies (6):\\
Let $R,\epsilon>0$.\\
Let $\xi:X\to{}l^1(X)$ be a map as in (3) with respect to $p=1$ the parameters $R,\epsilon$.  Assume
without loss of generality that each $\xi_x$ is positive-valued.\\
Define for each $x\in{}X$ a map $\phi_x:X\to[0,1]$ by setting $\phi_x(y)=\xi_y(x)$.
Note that $\sum_y\phi_x(y)=\sum_y\xi_x(y)=\|\xi_x\|_1=1$, whence $\{\phi_x\}_{x\in{}X}$ is
a partition of unity on $X$.\\
Note also that each $\phi_x(y)=\xi_y(x)$ can only be non-zero when $d(x,y)<S$ for some $S$, whence
$\phi_x$ is supported in $\bar{B}(x,S)$ for all $x$, i.e. $\{\phi_x\}$ is subordinated to the 
cover $\{\bar{B}(x,S):x\in{}X\}$, which consists of sets of uniformly bounded diameter.\\
Finally, note that if $d(y,y')\leq{}R$, then we get that
\begin{displaymath}
\sum_{x\in{}X}|\phi_x(y)-\phi_x(y')|=\sum_{x\in{}X}|\xi_y(x)-\xi_{y'}(x)|=\|\xi_y-\xi_{y'}\|_1<\epsilon.
\end{displaymath}

(6) implies (2):\\
Let $R,\epsilon>0$, and let $\mathcal{U}=\{U_i\}_{i\in{}I}$ and $\{\phi_i\}_{i\in{}I}$ have the properties
in (6) with respect to $R,\epsilon$; say in particular that the diameters of the $U_i$ are bounded by $S$.
Endow $I$ with the discrete metric, i.e. $d(i,j)=0$ if $i=j$ and is $1$ otherwise, and consider
$X\times{}I$ with the $l^1$ product metric, i.e. $d((x,i),(y,j))=d_X(x,y)+d_I(i,j)$.\\
For each $i\in{}I$, choose $y_i\in{}U_i$, and
define $\xi:X\times{}I\to{}l^1(X\times{}I)$ by setting 
\begin{displaymath}
\xi_{(x,i)}(y,j)=\left\{\begin{array}{ll} \phi_{j}(x), & y=y_j \\
0, & \textrm{otherwise.}\end{array}\right.
\end{displaymath}
Note then that 
\begin{displaymath}
\|\xi_{(x,i)}\|_1=\sum_{(y,j)\in{}X\times{}I}\xi_{(x,i)}(y,j)=\sum_{j\in{}I}\phi_{j}(x)=1,
\end{displaymath}
and that each $\xi_{(x,i)}(y,j)$ can be non-zero only if $x,y\in{}U_j$, whence
$\xi_{(x,i)}$ is supported in $B((x,i),S+2)$.\\
Finally, note that if $d((x,i),(y,j))\leq{}R$, then $d(x,y)\leq{}R$, and
\begin{displaymath}
\|\xi_{(x,i)}-\xi_{(y,j)}\|_1=\sum_{(z,k)\in{}X\times{}I}|\xi_{(x,i)}(z,k)-\xi_{(y,j)}(z,k)|
=\sum_{k\in{}I}|\phi_k(x)-\phi_k(y)|<\epsilon.
\end{displaymath} 
It follows that $X\times{}I$ satisfies (2), whence $X$ does too, as it is coarsely equivalent to
$X\times{}I$, and (2) is a coarse invariant (this is easy to check).\\

(3) implies (7):\\
Set p=2, and $\mathcal{H}=l^2(X)$.\\

(7) implies (8):\\
Let $f$ be as in (7), with constants $R,\epsilon$.\\
Define $k:X\times{}X\to\mathbb{R}$ by $k(x,y)=$Re$\langle{}f(x),f(y)\rangle$.
It is of finite propagation by the condition on $f$ in (7); it is
clearly symmetric; it is of
positive type, as for any $x_1,...,x_n\in{}X$ and any
$\lambda_1,...,\lambda_n\in\mathbb{R}$, we have that
\begin{displaymath}
\sum_{i,j=1}^n\lambda_i\lambda_jk(x_i,x_j)
=\textrm{Re}\left\langle\sum_{i=1}^n\lambda_if(x_i),
\sum_{i=1}^n\lambda_if(x_i)\right\rangle\geq0;
\end{displaymath}
it is normalised as $k(x,x)=$Re$\langle{}f(x),f(x)\rangle=\|f(x)\|^2=1$.\\
Finally, note that if $d(x,y)\leq{}R$, then $\|f(x)-f(y)\|^2<\epsilon^2$.
It follows by expanding 
$\|f(x)-f(y)\|^2=\|f(x)\|^2+\|f(y)\|^2-2$Re$\langle{}f(x),f(y)\rangle$,
that $0<1-k(x,y)<\epsilon^2/2$.\\

(8) implies (9):\\
Convolution with $k$ as given (provisionally) defines a linear
operator $T_k:l^2(X)\to{}l^2(X)$ by:
\begin{displaymath}
(T_k\xi)_x=\sum_{y\in{}X}k(x,y)\xi_y
\end{displaymath}
(cf. \ref{c*facts1}).  
As $k$ is finite propagation, say of propagation $S$, and $X$ is bounded 
geometry there exists $N=\sup_{x\in{}X}|B(x,S)|$ such
that for each $x\in{}X$, $k(x,y)\not=0$ for at most $N$ $y\in{}X$.\\
Note also that $|k(x,y)|\leq\sqrt{|k(x,x)||k(y,y)|}=1$, 
as the Cauchy-Schwarz inequality can be proved for $k$ (proof: as $k$
is of positive type and symmetric, for any $t\in\mathbb{R}$, we get that
$t^2k(x,x)+2tk(x,y)+k(y,y)\geq0$; the discriminant of this quadratic is thus
non-positive which gives Cauchy-Schwarz in the usual way).\\
We conclude that
\begin{align*}
\|T_k(\xi)\|_2^2 &
=\sum_{x\in{}X}\left|\sum_{y_\in\bar{B}(x,S)}k(x,y)\xi_y\right|^2 \leq\sum_{x\in{X}}\left(\sum_{y\in\bar{B}(x,S)}|\xi_y|\right)^2 \\
& \leq{}\sum_{y\in{}X}N\sum_{x\in\bar{B}(y,S)}|\xi_y|^2  \leq{}\sum_{y\in{}X}N^2|\xi_y|^2\leq{}N^2\|\xi_y\|_2^2.
\end{align*}
Hence $T_k$ is well-defined and bounded by $N$.\\
It is also of finite propagation in the sense of
definition \ref{finprop}; in fact, in the notation of that definition, 
$\sup\{d(x,y):P_{\{x\}}TP_{\{y\}}\not=0\}\leq{}S$, which completes the 
proof.\\

(9) implies (2):\\
Let $k$ be as given in (9) with respect to some $R,\epsilon$, and denote the
corresponding convolution operator in $C_u^*$ by $T_k$ (cf. \ref{c*facts1}).\\
Note that $T_k$ is a positive operator, as shown in \ref{c*facts1}.
It thus has a (unique) positive square root $T_l$, corresponding to some positive type kernel $l$.
Note that as elements of $C_u^*$ consist of limits of operators of
bounded support, we can find an operator $T_m$ (associated to a kernel $m$) of bounded support that 
satisfies
\begin{displaymath}
\|T_l-T_m\|<\min\left(\epsilon,\frac{\epsilon}{2(\|T_l\|+\epsilon)}\right).
\end{displaymath}
Py passing to $\frac{1}{2}(T_m+T_m^*)$, we may assume that $T_m$ is self-adjoint. 
As $T_m$ is of finite propagation, the kernel $m$ is also of finite propagation (from \ref{c*facts1} again).  Note also that $\|T_m\|\leq\|T_l\|+\epsilon$, so
\begin{displaymath}
\|T_k-T_m^2\|=\|T_l^2-T_m^2\|\leq\|T_l\|\|T_l-T_m\|+\|T_l-T_m\|\|T_m\|<\epsilon.
\end{displaymath}
Now, for each $x\in{}X$, let $\theta_x\in{}l^2(X)$ be given by $\theta_x(y)=m(x,y)$.  Note then that
\begin{displaymath}
\langle\theta_x,\theta_y\rangle=
\sum_{z\in{}X}m(x,z)\overline{m(y,z)}=\sum_{z\in{}X}m(x,z)m(z,y),
\end{displaymath}
i.e. the kernel $(x,y)\mapsto\langle\theta_x,\theta_y\rangle$ maps $(x,y)$ to the $(x,y)^{\textrm{th}}$ 
matrix entry of $T_m^2$.  Moreover, $T_m^2$ differs from $T_k$ by at most $\epsilon$ (in operator norm), whence for any $(x,y)\in{}X\times{}X$, the $(x,y)^{\textrm{th}}$
matrix coefficient of $T_m^2$ differs from the $(x,y)^{\textrm{th}}$ matrix coefficient of $T_k$ by at most $\epsilon$, i.e. that $\langle\theta_x,\theta_y\rangle$ differs from $k(x,y)$
by at most $\epsilon$.
Now, letting $d(x,y)\leq{}R$, we get that:
\begin{displaymath}
|1-\langle\theta_x,\theta_y\rangle|
\leq|\langle\theta_x,\theta_y\rangle-k(x,y)|+|1-k(x,y)|<2\epsilon
\end{displaymath}
i.e. $(x,y)\mapsto\langle\theta_x,\theta_y\rangle$ has $(R,2\epsilon)$ 
variation, whence $\theta$ itself has $(R,\sqrt{6\epsilon})$ variation, by the 
relationship 
\begin{displaymath}
\|\theta_x-\theta_y\|_2^2=
\|\theta_x\|_2^2+\|\theta_y\|_2^2-2\textrm{Re}\langle\theta_x,\theta_y\rangle.
\end{displaymath}
Continuing, note that $\|\theta_x\|_2^2\geq1-2\epsilon$, whence we can set
$\eta_x=\frac{\theta_x}{\|\theta_x\|_2}$, and use the estimates in the proof that (5) implies (2)
to conclude that $\eta$ has $(R,2\sqrt{6\epsilon/(1-2\epsilon)})$ variation.  
This implies that $\eta$ satisfies parts (a) and (b) from part (2) (for a suitable initial choice of $\epsilon$) with respect to $p=2$.
Further, as $m$ is of finite propagation, there exists $S>0$ such that if $d(x,y)>S$, then $m(x,y)=0$, whence also $\eta_x(y)=0$.  This says that $\eta$ satisfies (c) from part (2) as well, 
which completes the proof.
\end{proof}

\pagebreak

\section{Permanence properties} \label{permprop}

\subsection{Introduction and basic examples} \label{permpropintro}

There are two main parts to this chapter: \ref{genpermprop} deals with some closure properties
of the class of (bounded geometry, discrete) metric spaces with property A; section \ref{groups}
specialises to the case of countable discrete groups, which furnish our most important class of examples.

Before moving onto these, however, we give a few very basic examples of property A spaces; together with the closure properties of section
\ref{groups}, they allow us to form a surprisingly large class of groups with property A (we will explain in that section what it means for a discrete group, not usually considered a metric object, to have property A).

\begin{finite} \label{finite}
Any finite metric space has property A.
\end{finite}

\begin{trees} \label{trees}
The vertex set of any tree $\mathcal{T}$ (i.e. graph with no loops) has property A.  To be more precise, define a metric on the vertex set $V$ of $\mathcal{T}$ by setting $d(v,w)$ to be number of edges in the (unique) path of shortest length between $v$ and $w$.  The claim is that $(V,d)$ has property A.
\end{trees}

\begin{proof} (\cite{y1}, 2.4)
Assume $\mathcal{T}$ is infinite.  Fix $R,\epsilon>0$ and assume $\epsilon<1$.\\
Let $V$ denote the vertex set of $\mathcal{T}$, and pick a geodesic ray (i.e. image of
an isometry $\mathbb{N}\to{}V$, with the obvious metric on $\mathbb{N}$), say $\gamma_0$.\\
For each $v\in{}V$, let $\gamma_v$ be the unique (by the uniqueness of 
shortest paths in a tree) geodesic ray starting from $v$ that follows the same path as $\gamma_0$ for
an infinite distance. Define for each $v\in{}V$ 
\begin{displaymath}
A_v=\{(w,1)\in{}V\times\mathbb{N}:w\in\gamma_v, d(w,v)\leq3R/\epsilon+1\}.
\end{displaymath}
Note that $3R/\epsilon\leq|A_v|\leq{}3R/\epsilon+1$ for all $v\in{}V$.\\
Note also that if $d(v,w)\leq{}R$, then there are at most $2R$ elements in $A_v\triangle{}A_w$ 
(uniqueness of shortest paths in a tree again; this bound is not optimum), whence there are at
least $3R/\epsilon-R$ elements in $A_v\cap{}A_w$.  It follows that
\begin{displaymath}
\frac{|A_v\triangle{}A_w|}{|A_v\cap{}A_w|}\leq\frac{2R}{3R/\epsilon-R}=\frac{2\epsilon}{3-\epsilon}<\epsilon.
\end{displaymath}
Thus property A as in the original definition is satisfied.
\end{proof}

It follows, for example, that $\mathbb{Z}$ (with the usual metric $d(n,m)=|n-m|$) has property A.
Note that we do not assume that each vertex has finite (or even countable) valency; we will need this in the proof 
that property A is preserved under amalgamated free products in \ref{fprod}. 

\subsection{General permanence properties}\label{genpermprop}

This section covers a variety of closure properties of the class of metric spaces with property A.
Perhaps the most intuitive of these are \ref{products} and \ref{unions}, which say
that property A is preserved under finite unions and direct products.  We also use some
of the results in the next section to get permanence properties of the class of countable discrete groups
with property A.  It is an easy corollary of our techniques that \emph{finite asymptotic dimension} (defined in \ref{fad} below) implies property A;
we finish the section with a proof of this (the result first appeared in \cite{hr}).

All spaces in this section are assumed to be discrete and of bounded geometry.  This
is not strictly necessary, but fits in better with the
emphasis of the piece as a whole, and avoids concerns about exactly when the various 
formulations of property A are equivalent in more general situations.  Having said this,
we will need that bounded geometry is not necessary for \ref{fibreing} in the proof
that property A is preserved under amalgamated free products, \ref{fprod}.

The definition below is meant to capture what it means for a collection of spaces
to have property A `uniformly'.  We state it in terms of \ref{equiv}, (6), but 
it is clear how to generalise it to any of the equivalent definitions in \ref{equiv}. 

\begin{equi} \label{equi}
A collection of metric spaces $\{X_i\}_{i\in{}I}$ is said to have 
\emph{equi-property A}, or \emph{property A uniformly}, if 
for all $R,\epsilon>0$ there exists $S>0$ and covers $\mathcal{U}_i$ for each $X_i$
satisfying the conditions in \ref{equiv} (6) above, and such that the elements of each 
$\mathcal{U}_i$ all have diameter at most $S$.
\end{equi}

The following basic lemma is important: it allows us to conclude that if $\{X_i\}$ is an equi-property A
family, and $\{X_i'\}$ is such that each $X_i'$ is coarsely equivalent to $X_i$ `uniformly in $i$',
then $\{X_i'\}$ is an equi-property A family too.

\begin{subequi} \label{subequi}
Say $\{X_i\}_{i\in{}I}$ and $\{X_i'\}_{i\in{}I}$ are families of metric spaces
and that $\{X_i\}$ has equi-property A.\\
Say moreover that the $X_i'$ are \emph{uniformly coarsely embeddable} in the $X_i$, i.e. that there exist
coarse embeddings $f_i:X_i'\to{}X_i$ and non-decreasing maps $\rho_1,\rho_2:\mathbb{R}^+\to\mathbb{R}^+$ 
such that $\rho_1(t)\to+\infty$ as $t\to+\infty$ and
\begin{displaymath}
\rho_1(d(x,y))\leq{}d(f_i(x),f_i(y))\leq\rho_2(d(x,y))
\end{displaymath}
for all $x,y\in{}X_i'$, all $i$.\\
Then $\{X_i'\}_{i\in{}I}$ is an equi-property A family.
\end{subequi}

\begin{proof}
Let $R,\epsilon>0$ and let $R'=\rho_2(R)$.\\
By equi-property A for $\{X_i\}$ there exists covers $\mathcal{U}_i=\{U_i^j\}_{j\in{}J_i}$ of $X_i$ and $S>0$ such that
the the diam($U_i^j$) is bounded by $S$ for all $i,j$, and subordinate partitions of unity
$\{\phi_i^j\}$ with the properties in \ref{equiv}, (6) with respect to the parameters $R',\epsilon$.\\
For each $i$, define a cover $\mathcal{V}_i=\{V_i^j\}$ of $X_i'$ by setting $V_i^j=f_i^{-1}(U^i_j)$ and
a subordinate partition of unity by setting $\psi_i^j=\phi_i^j\circ{}f_i$.\\
Note then that:
\begin{itemize}
\item the diameter of each $V_i^j$ is bounded by any $S'$ such that $\rho_1(S')\geq{}S$ (and such an $S'$
exists as $\rho_1$ is non-decreasing and unbounded);
\item if $x,y\in{}X_i$ and $d(x,y)\leq{}R$, then $d(f_i(x),f_i(y))\leq{}R'$, whence 
\begin{displaymath}
\sum_{j\in{}J_i}|\psi_i^j(x)-\psi_i^j(y)|
=\sum_{j\in{}J_i}|\phi_i^j(f_i(x))-\phi_i^j(f_i(y))|<\epsilon.
\end{displaymath}
\end{itemize}
\end{proof}

Considering the special case where $X'$ is actually a subspace of $X$, 
inspection of the proof above reveals that a subspace of a property A space has property A 
`at least as well' as the original space.

The next result is the fundamental `gluing lemma' from which all the other results
in this section follow.

\begin{gluing}[\cite{dg2}, theorem 3.1] \label{gluing}
Say $X$ is a metric space such that for all $R,\epsilon>0$ there exists
a partition of unity $\{\phi_i\}_{i\in{}I}$ on $X$ such that:
\begin{itemize}
\item for all $x,y\in{}X$, if $d(x,y)\leq{}R$, then 
$\sum_{i\in{}I}|\phi_i(x)-\phi_i(y)|<\epsilon$;
\item $\{\phi_i\}$ is subordinated to an equi-property A cover $\mathcal{U}=\{U_i\}_{i\in{}I}$.
\end{itemize}
Then $X$ has property A.
\end{gluing}

\begin{proof}
Let $R,\epsilon>0$, and let $\mathcal{U}=\{U_i\}_{i\in{}I}$ be an equi-property A cover of $X$ and 
$\{\phi_i\}$ a subordinate partition of unity satisfying the properties above with respect to
$R, \epsilon/2$.\\
Note that as $\mathcal{U}$ is equi-property A, so too is $\mathcal{U}_R=\{U_i(R)\}_{i\in{}I}$, where
$U_i(R)=\cup_{x\in{}U_i}\bar{B}(x,R)$, by the previous lemma ($\rho_1(t)=t-2R, \rho_2(t)=t$ in the
notation of that lemma; the maps $f_i$ send a point $x\in{}U_i(R)$ to any point of $U_i$ within
$R$ of it).\\
Hence for each $i$, equi-property A for $\mathcal{U}_R$ implies that there exists
a cover $\mathcal{V}_i=\{V_i^j\}_{j\in{}J_i}$ of $U_i(R)$ and subordinate partition of unity
$\{\psi^j_i\}$ satisfying the properties in \ref{equiv}, (6) with respect
to $R,\epsilon/2$ and
such that there is a bound on the diameters of the $V^j_i$ that is uniform in both $i$ and $j$.\\
Extend each $\psi^j_i$ by setting it equal to zero outside of $U_i$, and define
maps $\theta_{ij}:X\to[0,1]$ by setting $\theta_{ij}(x)=\phi_i(x)\psi_i^j(x)$.\\
Note first that $\{\theta_{ij}\}$ is subordinated to the uniformly bounded cover
$\{V_i^j\}$ of $X$, and that
\begin{displaymath}
\sum_{i\in{}I, j\in{}J_i}\theta_{ij}(x)=\sum_{i\in{}I}\phi_i(x)\sum_{j\in{}J_i}\psi_i^j(x)=1
\end{displaymath}
by Fubini's theorem.\\
Hence $\{\theta_{ij}\}$ is a partition of unity on $X$ subordinated to a uniformly bounded cover.\\
Moreover, note that
\begin{align*}
\sum_{i,j}|\theta_{ij}(x)-& \theta_{ij}(y)| \\
& \leq\sum_i\phi_i(x)\sum_j|\psi^j_i(x)-\psi^j_i(x)|
+\sum_i|\phi_i(x)-\phi_i(y)|\sum_j\psi_i^j(y) ~~(*)
\end{align*}
by the triangle inequality.\\
Now, say $d(x,y)\leq{}R$.  If $x\in{}U_i$, then $y\in{}U_i(R)$, whence 
\begin{displaymath}
\sum_i\phi_i(x)\sum_j|\psi^j_i(x)-\psi^j_i(x)|<\epsilon/2
\end{displaymath}
by assumption on $\{\psi_i^j\}$, and as $\sum_i\phi_i(x)=1$.\\
Moreover, the second term in $(*)$ above is also bounded by $\epsilon/2$, as
$\sum_j\psi_i^j(y)$ is either $1$ (if $y\in{}U_i$) or $0$ (otherwise).\\
It now follows that $X$ has \ref{equiv}, (6).
\end{proof}

The first corollary is one of a family of results sometimes called `fibreing lemmas': the main motivating
examples are fibre bundles $p:X\to{}Y$ with base space $Y$ and total space $X$.
Other results in this family cover 
finite asymptotic dimension \cite{bd} (see \ref{fad} for a definition) and coarse embeddability in 
Hilbert space \cite{dg2}.  We will
use it to show that property A for countable discrete groups is preserved under extensions and free
products with amalgam in the next section.

\begin{fibreing}[\cite{dg2}, corollary 3.3] \label{fibreing}
Say $p:X\to{}Y$ is a bornologous map of
metric spaces, that $Y$ has property A, and that for any
uniformly bounded cover $\{U_i\}$ of $Y$, the pullbacks $\{p^{-1}(U_i)\}$ form an equi-property A
family of subspaces of $X$.  Then $X$ has property A.
\end{fibreing}

\begin{proof}
Let $R,\epsilon>0$.  As $p$ is bornologous, there exists $R>0$ such that
$d(x,x')\leq{}R$ gives $d(p(x),p(x'))\leq{}R'$.  By property A for $Y$ there exists
a uniformly bounded cover $\{U_i\}$ of $Y$ and subordinate partition of unity
$\{\phi_i\}$ as in \ref{equiv} (6) with respect to the parameters $R',\epsilon$.\\
Then $\{\phi_i\circ{}p\}$ is a partition of unity subordinate to the equi-property A
cover $p^{-1}(U_i)$, whence the result by \ref{gluing}.
\end{proof}

\begin{products} \label{products}
Let $X,Y$ be metric spaces with property A and say that 
$X\times{}Y$ is endowed with a metric that restricts to the metric
on $X$ (respectively $Y$) on each fibre $X\times\{y\}$ ($\{x\}\times{}Y$).  Say
also that $d((x,y),(x',y'))\geq\max\{d(x,x'),d(y,y')\}$ for all pairs $(x,y),(x',y')$ of elements
of $X\times{}Y$
(e.g. any of the $l^p$ metrics, $1\leq{}p\leq\infty$, satisfy these conditions).\\
Then $X\times{}Y$ has property A.
\end{products}

\begin{proof}
Let $p:X\times{}Y\to{}Y$ be the projection map.  It is bornologous and for any uniformly
bounded cover $\{U_i\}$ of $Y$ the pullbacks $p^{-1}(U_i)$ are uniformly coarsely equivalent to
$X$, so have property A uniformly, whence the result by the above. 
\end{proof} 

The next result is part of the folklore of the subject.  We will use it in the three subsequent corollaries to prove that property A is closed under 
certain types of union. 

\begin{covers}[\cite{bell}, proposition 1] \label{covers}
Let $\mathcal{U}=\{U_i\}_{i\in{}I}$ be a cover of a metric space $X$ with multiplicity $k$ (i.e. each point of $X$ is
contained in at most $k$ elements of $\mathcal{U}$) and Lebesgue number $L$ (i.e. any ball of radius at
most $L$ is wholly contained in one element of $\mathcal{U}$).\\
Then there exists a partition of unity $\{\phi_i\}_{i\in{}I}$ such that
\begin{displaymath}
\sum_{i\in{}I}|\phi_i(x)-\phi_i(y)|\leq\frac{(2k+2)(2k+3)}{L}d(x,y)
\end{displaymath}
for all $x,y\in{}X$.
\end{covers}

\begin{proof}
Define $\phi_i:X\to[0,1]$ for each $i\in{}I$ by setting
\begin{displaymath}
\phi_i(x)=\frac{d(x,X\backslash{}U_i)}{\sum_{j\in{}I}d(x,X\backslash{}U_j)};
\end{displaymath}
it is well-defined, as if not $x\not{\in}{}V$ for all $V\in\mathcal{U}$.\\
Note that $d(x,y)\geq|d(x,X\backslash{}U_i)-d(y,X\backslash{}U_i)|$ by the triangle inequality and
that for any $x$, $d(x,X\backslash{}U_i)\geq{}L$ for at least one $i\in{}I$, as at least
one of the $U_i$ contains $B(x,L)$.  It follows that
\begin{align*}
|\phi_i & (x) -\phi_i(y)| 
 \leq\left|\frac{d(x,X\backslash{}U_i)}{\sum_{j\in{}I}d(x,X\backslash{}U_j)}-
\frac{d(y,X\backslash{}U_i)}{\sum_{j\in{}I}d(x,X\backslash{}U_j)}\right| \\
&+\left|\frac{d(y,X\backslash{}U_i)}{\sum_{j\in{}I}d(x,X\backslash{}U_j)}-
\frac{d(y,X\backslash{}U_i)}{\sum_{j\in{}I}d(y,X\backslash{}U_j)}\right| \\
& \leq\frac{d(x,y)}{\sum_{j\in{}J}d(x,X\backslash{}U_j)}
+\frac{d(y,X\backslash{}U_i)}{\sum_{j\in{}I}d(y,X\backslash{}U_j)}.
\frac{\sum_{j\in{}I}|d(y,X\backslash{}U_j)-d(x,X\backslash{}U_j)|}{\sum_{j\in{}I}d(x,X\backslash{}U_j)} \\
& \leq\frac{2k+3}{L}d(x,y),
\end{align*}
as at most $2k+2$ in the rightmost sum are non-zero, by assumption on the multiplicity of
$\mathcal{U}$.\\
The result now follows, as at most $2k+2$ terms in the sum
$\sum_{i\in{}I}|\phi_i(x)-\phi_i(y)|$ are non-zero.
\end{proof}

The next three corollaries on unions of spaces with property A are easy once one has made the 
observation \ref{covers}, 
and more results along the same line are derivable (see for example \cite{bnw}, where
it is proved that a locally finite space has property A if and only if the family consisting of all its finite subspaces has equi-property A).  
Note, however, that one
\emph{cannot} simply say that the union of a countable equi-property A family has property A;
this follows from the existence of the countable non-property A spaces discussed in chapter
\ref{nona}. 

\begin{unions}[\cite{dg2}, corollary 4.5] \label{unions}
Let $X=\cup_{i=1}^nX_i$ be a bounded geometry discrete metric space such that each $X_i$ (with 
the induced metric) has property A.  Then $X$ has property A.
\end{unions}

\begin{proof}
Expand each $X_i$ by $L>0$ to get $U_i=\cup_{x\in{}X_i}\bar{B}(x,L)$; the resulting cover of $X$ has 
multiplicity at most $n$ and Lebesgue number at least $L$.  It has property A uniformly as each
$U_i$ is coarsely equivalent to each $X_i$ and there are only finitely many of them.\\
The result is now immediate from \ref{gluing} and \ref{covers} by making $L$ suitably
large.
\end{proof}

\begin{unions2}[\cite{dg2}, corollary 4.6] \label{unions2}
Say $X$ is a metric space covered by an equi-property A family $\mathcal{U}$, and that
for every $L>0$ there exists a property A subspace $Y_L$ of $X$ such that the family
$\mathcal{F}=\{U\backslash{}Y_L:U\in\mathcal{}U\}$ is $L-separated$, i.e. 
$inf\{d(x,y):x\in{}V,y\in{}V'\}>L$ for all elements $V\not=V'$ of $\mathcal{F}$.\\
Then $X$ has property A.
\end{unions2}

\begin{proof}
The cover $\mathcal{V}=\mathcal{F}\cup\{Y_L\}$ of $X$ is an equi-property A family.
Expand each element to get $\mathcal{V'}=\{V(L/2):V\in\mathcal{V}\}$.\\
This cover then has Lebesgue number at least $L/2$, and multiplicity at most $2$.  The result
now follows from \ref{covers}.
\end{proof}

\begin{unions3} \label{unions3}
Let $X_1\subseteq{}X_2\subseteq{}X_3\subseteq{}...$ be an increasing sequence of
bounded metric spaces, and let $X=\cup_{n=1}^{\infty}X_n$.  Assume also that
any bounded subset of $X$ is contained in some $X_n$ (for example, these conditions
are satisfied by the sequence $(B(x,n))_{n=1}^{\infty}$ for any fixed $x\in{}X$, or for an increasing sequence
of finite subsets of a locally finite space).\\
If the $X_n$ have property A uniformly, then $X$ has property A.
\end{unions3}

\begin{proof}
Let $L>0$.  Then there exists a subsequence $(X_{n_k})$ of $(X_n)$ such
that
\begin{displaymath}
\left(\bigcup_{x\in{}X_{n_k}}\bar{B}(x,L)\right)\cap
\left(\bigcup_{x\in{}X_{n_{k+2}}\backslash{}X_{n_{k+1}}}\bar{B}(x,L)\right)
=\emptyset;
\end{displaymath}
for given $X_{n_k}$, choose $X_{n_{k+1}}$ to contain $\cup_{x\in{}X_{n_k}}\bar{B}(x,3L)$, which
is possible by the assumptions on boundedness.\\
It follows that by setting
\begin{displaymath}
U_k=\bigcup_{x\in{}X_{n_{k+1}}\backslash{}X_{n_k}}\bar{B}(x,L)
\end{displaymath}
for each $k\geq{}1$ we get a cover $\mathcal{U}=\{U_k\}_{k=1}^{\infty}$ of multiplicity at most $2$ and
Lebesgue number at least $L$; it is equi-property A as the $U_k$ are uniformly coarsely equivalent
to the sets $X_{n_{k+1}}\backslash{}X_{n_k}$, which in turn are subspaces of the equi-property A sequence
$(X_{n_k})_{k=1}^{\infty}$.\\
The result now follows from \ref{gluing} and \ref{covers} by choosing $L$ suitably large.
\end{proof}

Finally in this section, we define what it means for a metric space to have \emph{finite asymptotic dimension}; 
it is then easy to use \ref{covers} to show that this new notion implies property A.
Asymptotic dimension was introduced by 
M. Gromov in \cite{gr1}, section 1.E; the idea is that it is a large-scale
analogue of the classical notion of \emph{topological covering dimension}.  It is a coarse 
invariant that has been extensively investigated (see e.g.
\cite{bd} or chapter 9 of \cite{roe1} for some results and further references). 

\begin{fad} \label{fad}
A metric space $X$ is said to have \emph{finite asymptotic dimension} if there exists
$k\geq{}0$ such that for all $L>0$ there exists a uniformly bounded cover of $X$ of Lebesgue
number at least $L$ and multiplicity $k+1$.  The least possible such $k$ is the \emph{asymptotic dimension} of $X$.
\end{fad}

\begin{fad2}\label{fad2}
Any metric space with finite asymptotic dimension also has property A.
\end{fad2}

\begin{proof}
Say $X$ has asymptotic dimension $k$.  Let $R,\epsilon>0$, and let 
\begin{displaymath}
L>\frac{R}{(2k+3)(2k+5)\epsilon}.  
\end{displaymath}
By definition of finite asymptotic dimension, then, there exists a uniformly bounded cover $\mathcal{U}=\{U_i\}_{i\in{}I}$
of $X$ with multiplicity at most $k+1$ and Lebesgue number $L$.  By the \ref{covers} there exists a partition of unity
$\{\phi_i\}_{i\in{}I}$ subordinate to $\mathcal{U}$ such that for all $x,y\in{}X$
\begin{align*}
\sum_{i\in{}I}|\phi(x)-\phi(y)|\leq\frac{(2(k+1)+1)(2(k+1)+3)}{L}d(x,y) <{}\frac{\epsilon}{R}d(x,y).
\end{align*}
In particular, then, if $d(x,y)\leq{}R$, we have that
\begin{displaymath}
\sum_{i\in{}I}|\phi(x)-\phi(y)|<\epsilon
\end{displaymath}
This shows that $X$ satisfies part (6) of \ref{equiv}, and thus has property A.
\end{proof}

\subsection{Property A for groups} \label{groups}

Countable discrete groups provide our most important examples of spaces with property A.  
We start with an explanation of how it makes sense to talk about a discrete group having
property A at all - it is perhaps not obvious how to make such a group into a
metric space.   

Nontheless, there really is
a canonical geometry on any countable discrete group.  Our large-scale
perspective is crucial: the geometry we define is only determined up to the loose 
standard of coarse equivalence.  As it happens, coarse
equivalence preserves quite a large amount of group-theoretic structure; a lot of the relevant ideas
originated in an inspiring paper by Gromov \cite{gr1}.  We, however,
will restrict ourselves to applications relevant to property A. 

After making the necessary geometric definitions, the bulk of this section deals with
closure properties of the class of discrete groups with property A: we 
prove closure under the taking of subgroups (this is trivial), direct limits, extensions,
free products with amalgam and HNN extensions.  This allows quite a large range 
of examples to be built up from scratch
(including the `classical' groups that one often meets first in an algebraic 
treatment of group theory: finite, abelian, nilpotent, solvable and free groups).  Our exposition
follows that in \cite{dg} and \cite{dg2}; another approach, based on the characterisation
of property A in theorem \ref{amenacta} and using more functional-analytic machinery, is followed in chapter 5 of \cite{bo}.

We start by defining a metric on any finitely generated group.  Although we will subsequently generalise the definition to any countable group,
the finitely generated case provides much of the intuition behind the subject, so we start here.

\begin{wmetric} \label{wmetric}
Let $G$ be a finitely generated discrete group, and $S=\{g_1,...,g_n\}$ a 
symmetric (i.e $g\in{}S$ implies $g^{-1}\in{}S$) generating set.  
The \emph{length} of $g\in{}G$ is $|g|_S=\inf\{k:g=h_1h_2...h_k 
~\textrm{for some} ~h_1,...,h_k\in{}S\}$.  The associated \emph{left-invariant 
word metric} is $d_S(g,h)=|g^{-1}h|_S$.
\end{wmetric}

This defines a metric on $G$ satisfying the usual axioms.  It is
left-invariant in the sense that for any $h,g,g'\in{}G$, 
$d_S(hg,hg')=d_S(g,g')$.  It is also of bounded 
geometry,
as for any $g$ $|B_{d_S}(g,r)|=|B_{d_S}(e,r)|$ by left invariance, and if $|S|=n$
this is bounded by $n^r$, the number of words of length at most $r$ in $n$
letters.

This metric is closely 
related to the metric on the
Cayley graph of $G$ associated to $S$, say $\mathcal{G}_S$. $d_S$ is the metric on this
giving length one to each side, restricted to the vertex set (equivalently,
$d_S(g,h)$ is the least number of edges needed to travel between the vertices $g$
and $h$).

For example, the metric on the free group $F_k$ with
respect to some generating set $\{x_1,...,x_k\}$ is the restriction of the metric on its
Cayley graph - the infinite degree $k$ tree - to the vertex set.  

\begin{integers} \label{integers}
It thus follows from \ref{trees} that $F_k$ (with the word metric associated to its usual
generating set) has property A for all $k$ (including 
$F_1=\mathbb{Z}$). 
\end{integers}

Throughout this section, definition \ref{wmetric} provides our motivating 
example.
The metric it defines has the very important property of being \emph{coarsely 
invariant} in the sense that if $S$, $S'$ are two finite generating sets for
$G$, then the identity map $(G,d_S)\to{}(G,d_{S'})$ is a coarse equivalence.
We are thus entitled to say that `$F_k$ has property A' without mentioning a
specific choice of metric.

As we have already mentioned, however, one can go further: there exists a unique (up to coarse equivalence)
left-invariant bounded geometry metric on any countable discrete group.  In
general, the word metric associated to a finite generating set has much better
properties than an arbitrary left-invariant bounded geometry metric (see
e.g. chapter 1 of \cite{roe1}); nontheless,
none of the results of this piece make use
of these (with the exception of some of those in \ref{propT}; however, that section is based on property T groups
which happen to automatically be finitely generated), so we prove them for any countable discrete group.
It is also in some ways quite natural that we generalise this far,
as \ref{subgp} and \ref{dlim} can both take one out of the class of
finitely generated groups; moreover our main motivating result \cite{y1} was extended
to all countable discrete groups in \cite{sty}.  Note, however, that once we have proved \ref{subgp} and \ref{dlim},
we can conclude that a group $G$ has property A if and only if all of its
finitely generated subgroups do.

\begin{wmce} \label{wmce}
Let $G$ be a countable discrete group.  Then there exists a left-invariant
metric $d$ on $G$ such that $(G,d)$ is a bounded geometry space.\\
Moreover, if $d'$ is another metric on $G$ with these properties, then the
identity map $(G,d)\to{}(G,d')$ is a coarse equivalence.
\end{wmce}

\begin{proof} 
Existence:\\
Order inverse pairs of elements of $G$ as $e=g_0,g_1^{\pm1},g_2^{\pm1},...$,
and define $f:G\to\mathbb{N}$ by $f(g_n^{\pm1})=n$.\\
For each $g\in{}G$, define
\begin{displaymath}
|g|=\inf\left\{f(g_{i_1})+...+f(g_{i_m}): g=\prod_{j=1}^mg_{i_j}^{\epsilon_j},
m\in\mathbb{N}, \epsilon_j=\pm1\right\}.
\end{displaymath}
$|.|$ satisfies $|g|=0$ if and only if $g=e$, $|g|=|g^{-1}|$ and
$|gh|\leq|g|+|h|$; the latter follows as the right hand side is the infimum 
over a smaller set than the left hand side.\\
It follows that $d(g,h):=|g^{-1}h|$ defines a left-invariant metric on $G$.\\
Moreover, $|B(e,r)|\leq(2r+1)^r$; this is the number of words in the at most
$2r+1$ elements of $G$ that satisfy $f(g)\leq{}r$ (and as such is a very bad
estimate; it is, however, finite).  By left invariance, then, $(G,d)$ has
bounded geometry.\\
Uniqueness:\\
Let $d'$ be another bounded geometry, left-invariant metric on $G$.
For any $R>0$, the map $Id:(G,d)\to(G,d')$ takes the (finite) ball $B_d(e,R)$ 
to a finite set, so one contained in some ball $B_{d'}(e,S)$.  Similarly, it 
pulls back the finite ball $B_{d'}(e,R)$ to a finite, so bounded, set.  By left
invariance it is thus bornologous and proper, so coarse.  Similarly, the
identity map in the other direction is a coarse map, whence coarse equivalence.
\end{proof}

We assume from now on (as the above entitles us to) that any countable discrete
group is equipped with an integer valued, bounded geometry left-invariant metric.  Such
a metric has associated to it a \emph{proper length function} given by 
$|g|=d(g,e)$; we will also use this without comment.  Here `proper' means
that the pullback of a compact (i.e. bounded) set is compact (i.e. finite).  Such a
length function is symmetric in the sense that $|g|=|g^{-1}|$.

Throughout the rest of this section all groups are assumed to be countable and discrete.
\ref{wmce} allows us to use any left-invariant bounded geometry metric we can find,
which we often exploit (e.g. it often makes sense to give a subgroup the restricted
metric from the original group).

\begin{subgp} \label{subgp}
If $G$ has property A, and $H$ is a subgroup of $G$, then $H$ has property A.
\end{subgp}

\begin{proof}
\ref{wmce} shows that the inclusion map
$H\hookrightarrow{}G$ is a coarse embedding.  The result is now immediate
from remark \ref{cinvrem}.
\end{proof}

\begin{dlim} \label{dlim}
Let $G$ be the direct limit of the sequence $H_1\hookrightarrow{}H_2\hookrightarrow{}H_3\hookrightarrow{}...$
of property A groups with injective connecting maps.
Then $G$ has property A.
\end{dlim}

\begin{proof} (\cite{dg}, theorem 3.1)
We use \ref{equiv}, condition (7).\\
Let $d$ be any left-invariant, bounded geometry metric on $G$, and let $R,\epsilon>0$.\\
Note first that as $(G,d)$ is of bounded geometry, $B_d(e,R)$ is finite, hence contained
in $H_N$ for some $N$.\\
Equip $H_N$ with the restriction of d, which defines a left-invariant, bounded geometry
metric on $H_N$;
as $H_N$ has property A, we can find a Hilbert space $\mathcal{H}$
and a map $\xi:H_N\to\mathcal{H}$ that satisfies the conditions in \ref{equiv}, (7)
relative to this metric.\\
Now, let $X\subseteq{}G$ be a set of coset representatives for $(G:H_N)$ in $G$, so each
$g\in{}G$ can be written uniquely $x_gh_g$ for some $x_g\in{}X$, $h_g\in{}H$.\\
If $g,g'$ are in the same coset of $H_N$, we conclude that $x_g=x_{g'}$,
whence
\begin{displaymath}
d(g,g')=d(x_gh_g,x_{g'}h_{g'})=d(h_g,h_{g'})
\end{displaymath}
by left invariance of d.\\
Construct now a map $\hat{\xi}:G\to\mathcal{H}\otimes{}l^2(X)$ 
(the `$\otimes$' is assumed a Hilbert space tensor product, so is a
completion of the algebraic vector space inner product) by
setting $\hat{\xi}_g=\xi_{h_g}\otimes\delta_{x_g}$, and show that it has the properties
in \ref{equiv} (7):\\
Checking 7(a):\\
$\|\hat{\xi}_g\|=\|\xi_{h_g}\|.\|\delta_{x_g}\|=1$.\\
Checking 7(b):\\
Say $d(g,g')<R$, whence $d(e,g^{-1}g')<R$ by left invariance, and so
$g^{-1}g'\in{}H_N$ by assumption on $N$.  It follows that $g$, $g'$ are in the
same coset of $H_N$.\\
Hence by the comment above, $d(h_g,h_{g'})=d(g,g')<R$, so we get that
\begin{align*}
\|\hat{\xi}_g-\hat{\xi}_{g'}\| &
=\|\xi_{h_g}\otimes\delta_{x_g}-\xi_{h_{g'}}\otimes\delta_{x_{g'}}\|
=\|(\xi_{h_g}-\xi_{h_{g'}})\otimes\delta_{x_g}\| ~\textrm{as} ~x_g=x_{g'} \\
& =\|\xi_{h_g}-\xi_{h_{g'}}\|<\epsilon
\end{align*}
by assumption on $\xi$.\\
Checking 7(c):\\
By assumption on $\xi$ there exists $S>0$ such that if $h,h'\in{}H_N$ and 
$d(h,h')\geq{}S$, we can conclude that $\langle\xi_h,\xi_{h'}\rangle=0$.\\
Say now that $g,g'\in{}G$, and that $d(g,g')\geq{}S$.  We have
\begin{displaymath}
\langle\hat{\xi}_g,\hat{\xi}_{g'}\rangle
=\langle\xi_{h_g},\xi_{h_{g'}}\rangle.\langle\delta_{x_g},\delta_{x_{g'}}\rangle.
\end{displaymath}
This is zero unless $x_g=x_{g'}$, i.e. $g,g'$ are in the same coset
of $H_N$.  In this case, $d(h_g,h_{g'})=d(g,g')\geq{}S$, and then by assumption
on $\xi$, $\langle\xi_{h_g},\xi_{h_{g'}}\rangle=0$.
\end{proof}

It follows immediately that a countable discrete group has property A if and only if
all of its finitely generated subgroups do; simply enumerate the elements of $G$ as
$g_1,g_2,g_3,...$, and $G$ is then the direct limit of the increasing union of
subgroups generated by the first $n$ elements.  

Note that this result is rather different from \ref{unions3}; the 
group structure on the direct
limit is essential to enable one to drop the boundedness assumptions necessary for that result.  On other hand,
our next result depends on \ref{fibreing}.

\begin{ext} \label{ext}
Say $1\to{}K\to{}G\to{}H\to{}1$ is a short exact sequence of groups, and that
$K$ and $H$ have property A.  Then $G$ has property A.
\end{ext}

\begin{proof} (\cite{dg2}, theorem 5.1)
Let $d_G$ be a left invariant bounded geometry metric on $G$
and let the metric on $K$ be the restriction of $d_G$; this is a left-invariant 
bounded geometry metric on $K$, so $(K,d_G|_{K\times{}K})$ has property A.
Let $p:G\to{}H$ be the quotient map, and define a length
function on $H$ by setting
\begin{displaymath}
|h|_H=\min\{|g|_G:p(g)=h\},
\end{displaymath}
where $|.|_G$ is the length function on $G$ associated to $d_G$.\\
Note that this length function is symmetric, satisfies $|hh'|_H\leq|h|_H+|h'|_H$ and
$|h|_H=0$ if and only if $h=e$, as $|.|_G$ has these properties.  Moreover 
\begin{displaymath}
|\{h\in{}H:|h|_H\leq{}r\}|\leq|\{g\in{}G:|g|_G\leq{}r\}|,
\end{displaymath}
whence $|.|_H$ is proper and thus defines a left-invariant bounded geometry
metric $d_H(h,h')=|h^{-1}h')|_H$ on $H$; in particular, $(H,d_H)$ has property A.  Note further that with this
metric the map $p$ is contractive (hence bornologous), and that the minimum in the expression
defining $|.|_H$ is actually obtained as $|.|_G$ is integer-valued.\\
Now, by \ref{fibreing} it will be sufficient to prove that
for any uniformly bounded cover $\{U_i\}$ of $H$, the pullbacks $p^{-1}(U_i)$ form
an equi-property A family.\\
Note, however, that for such a cover, there exist $S>0$ and elements $h_i$ of $H$ such
that $p^{-1}(U_i)\subseteq{}p^{-1}(B_H(h_i,S))$ for all $i$, whence it suffices to prove that
the sets $p^{-1}(B_H(h_i,S))$ form an equi-property A family (by \ref{subequi}).\\
Finally, note that if $g_i$ is such that $p(g_i)=h_i$, then
$g_ip^{-1}(B_H(e,S))=p^{-1}(B_H(h_i,S))$, whence $p^{-1}(B_H(h_i,S))$ is isometric 
to $p^{-1}(B_H(e,S))$ for all $i$. Hence it suffices to prove that this has property A for
any $S$.\\
In fact, this is true as it is coarsely equivalent to $K$.  For let $g\in{}p^{-1}(B_H(e,S))$
and let $g'\in{}G$ be such that $p(g')=p(g)$ and $|g'|_G=|p(g)|_H<S$; such a $g'$ exists by choice
of the length function on $H$.\\
Note then, however, that $gg'^{-1}\in{}K$, and $d_G(g,gg'^{-1})=|g'^{-1}|_G<S$, i.e. for any element
of $p^{-1}(B_H(e,S))$, there exists an element of $K$ within $S$ of it, whence $K$ is coarsely
dense in $p^{-1}(B_H(e,S))$, whence coarsely equivalent to it.
\end{proof}

It follows that all solvable (hence also all nilpotent) groups have property A, for:
\begin{itemize}
\item all finitely generated abelian groups have property A by their structure theorem, \ref{finite}, 
\ref{integers} and the last result;
\item all countable abelian groups have property A by \ref{dlim};
\item all solvable groups have property A by the last result, as any solvable
group is arrived at by taking a finite number of extensions of abelian groups.
\end{itemize}
 
\begin{produ} \label{produ}
The class of groups with property A is closed under semidirect products, and in particular, direct products.
\end{produ}

\begin{proof}  These are special cases of extensions. \end{proof}

\begin{wreath}
If $G$ and $H$ are property A groups, then so too is the reduced wreath
product $G\wr{}H$.
\end{wreath}

\begin{proof}
$G\wr{}H$ is formed by taking $\oplus_{h\in{}H}G=\lim_n\oplus_{H_n}G$, where $H_n$ is the 
first $n$ elements of $H$ under any ordering, then taking a semidirect product with $H$ with respect 
to the action of $H$ that permutes the different copies of $G$.
\end{proof}

The next result of this section, that property A is preserved under free products with amalgamation,
is the most involved.  It is not, however, too hard once we introduce some very useful technology, mainly
from \cite{ser}.  The first step is to construct the so-called \emph{Bass-Serre tree} of a free 
product.

\begin{bassserre}
Let $F=G\ast_KH$ be the free product of $G$, $H$ over a common subgroup $K$.\\
Define a vertex set by setting $V=(F:G)\sqcup(F:H)$, the disjoint union of the cosets of $G$ and $H$ in $F$;
similarly an edge set $E=(F:K)$, the cosets of $K$.\\
Say the element $fK$ of $E$ connects the elements $fG$, $fH$ of $V$.\\
Then this construction gives a tree $\mathcal{T}_F$ (the \emph{Bass-Serre tree} associated 
to the decomposition $F=G\ast_KH$) on which $F$ acts by left-multiplication.
\end{bassserre}

\begin{proof}
$\mathcal{T}_F$ really is a tree:
\begin{itemize}
\item It is well-defined as $fK=f'K$ implies that $fG=f'G$ and $fH=f'H$.
\item It is connected. For example, let $g_1h_1...g_nh_nG$ be a vertex.  The
sequence of edges $g_1h_1...g_nh_nK,g_1h_1...g_ng_nK,...,g_1h_1K,g_1K$ takes it to 
$g_1h_1...g_nh_nH=g_1h_1...g_nH$,
to $g_1h_1...h_{n-1}g_nG=g_1h_1...h_{n-1}$ to ... to $G$.
\item There are no non-trivial loops.  For consider a sequence of edges
\begin{displaymath}
f_1G\to{}f_1H=f_2H\to{}f_2G=f_3{}G\to...\to{}f_nH=f_1H\to{}f_1G.
\end{displaymath} 
This implies that $f_i^{-1}f_{i+1}$ is an element of $H$ if $i$ is odd and $G$ if $i$
is even.  We also have that $\Pi{}f_i^{-1}f_{i+1}=e$; by definition of the 
amalgamated free product, this can only happen if at least one of the $f_i^{-1}f_{i+1}$
is actually in $K$, and the whole product $\Pi{}f_i^{-1}f_{i+1}$ is of the form $f^{-1}f$
for some $f\in{}F$.
Hence the loop is trivial.
\end{itemize}
\end{proof}

The next stage constructs a metric on the total space (disjoint union of all the
cosets that give the vertices), which uses the Bass-Serre tree to 
`organise' the cosets (themselves metric spaces) making up its vertices.

\begin{metrictree} \label{metrictree}
Let $F, \mathcal{T}_F$ be as above, and let $d_F$ be a left-invariant bounded geometry metric
on $F$.\\
For notational convenience, write an element $x$ of the vertex $v\in{}V$ (i.e. $v\subseteq{}F$ is a coset of either $G$
or $H$ and $x$ is an element of $v$) as $x_v$. Metrise the disjoint union of all the cosets in $V$ by giving it the largest metric satisfying
\begin{displaymath}
d(x_v,y_w)\leq\left\{\begin{array}{ll}
d_F(x_v,y_w), & v=w, ~\textrm{so $x_v,y_w$ are elements of the same coset in $V$} \\
1, & v\not=w, ~\textrm{but $x_v=y_w$ in $F$.} \end{array}\right.
\end{displaymath}
Denote the resulting space by $(X_F,d)$.\\
Then:
\begin{itemize}
\item $X_F$ is locally finite;
\item $F$ acts on $X_F$ by isometries by setting $f.x_v=(fx)_{fv}$; 
\item the map $p:X\to\mathcal{T}_F$
defined by $p(x_v)=v$ is equivariant and contractive (whence bornologous); 
\item the inclusions
$fG\hookrightarrow{}X_F$, $fH\hookrightarrow{}X_F$ are isometries.
\end{itemize}
\end{metrictree}

\begin{proof}
Note first that the metric on $X_F$ makes sense.  For define $\hat{d}$ on 
$S=\{(x_v,y_w)\in{}X\times{}X:v=w ~\textrm{or}~ x_v=y_w\}$ by setting
\begin{displaymath}
\hat{d}(x_v,y_w)=\left\{\begin{array}{ll}
d_F(x_v,y_w), & v=w, ~\textrm{so $x_v,y_w$ are elements of the same coset in $V$} \\
1, & v\not=w, ~\textrm{but $x_v=y_w$ in $F$.} \end{array}\right.
\end{displaymath}
Note that any metric $d$ on $X_F$ satisfying the above must be no larger than $\hat{d}$ on $S$;
on the other hand, any metric satisfying the above must also satisfy
\begin{displaymath}
d(x,y)\leq\inf\left\{\sum_{i=1}^n\hat{d}(x_i,x_{i+1}):
n\in\mathbb{N}, (x_i,x_{i+1})\in{}S, x_1=x,x_{n+1}=y\right\}
\end{displaymath}
by the triangle inequality.  By connectedness of $\mathcal{T}_F$, however, the expression on the
right hand side of the above always exists, so \emph{is} a metric (cf. \ref{wmce}) and 
must be the metric $d$ on $X_F$.

In order to get the comments above, we claim that $d$ is give concretely by the formula
\begin{displaymath}
d(x_v,y_w)=\inf\left\{\sum_{i=1}^p\hat{d}(x_i,x_{i+1})\right\}=
d_{\mathcal{T}_F}(v,w)+\inf\left\{\sum_{i=1, ~\textrm{$i$ odd}}^p\hat{d}(x_i,x_{i+1})\right\},
\end{displaymath}
where the infimum is over all sequences $x_1,...,x_{p+1}$ such that
\begin{enumerate}
\item $p=2d_{\mathcal{T}_F}(v,w)+1$;
\item $x_1=x_v,x_{p+1}=y_w$;
\item pairs $(x_i,x_{i+1})$ are both in the same coset (vertex) if $i$ is odd,
or such that $x_i=x_{i+1}$ if $n$ is even;
\item $x_i\in{}v_i$, and the vertices $v=v_1,v_3,...,v_p=w$ form a geodesic (shortest) path
between $v$ and $w$ in $\mathcal{T}_F$.
\end{enumerate}
Given this, the comments above are all (nearly) immediate.\\
Proof of claim:\\
Denote by $\tilde{d}$ the metric (it is an easy check to see that it is one) given by the right hand
side of the expression above; note that it agrees with $\hat{d}$ on $S$, whence $d\geq\tilde{d}$.\\
On the other hand, note that if $p$ in the above were not fixed, we would have $d\leq\tilde{d}$ by the
triangle inequality; it thus suffices to show that if $x_1,...,x_{2n}$ is any sequence of
elements in $X$ satisfying (2) and (3) above with $x_i\in{}v_i$, but $v_1,v_3,...,v_{2n-1}$ 
\emph{not} a geodesic path in $\mathcal{T}_F$, then elements can be removed to get a new
sequence $y_1,...,y_{2n-4}$ with $\sum_{i=1}^{2n-1}d(x_i,x_{i+1})\geq\sum_{i=1}^{2n-5}d(y_i,y_{i+1})$.\\
Now, given such a sequence $x_1,...,x_{2n}$, there exists $i$ with $v_{2i+1}=v_{2i+5}$; this says simply that
the path has gone back on itself, which is the only way tree paths fail to be
geodesic.  Hence we have a situation like
\begin{displaymath}
\begin{array}{ccccccccc}
... & \to & v_{2i}=v_{2i+1} & \to & v_{2i+2}=v_{2i+3} & \to & v_{2i+4}=v_{2i+5} & \to & ... \\
 & & x_{2i}~~~x_{2i+1} & & x_{2i+2}~~~x_{2i+3} & & x_{2i+4}~~~x_{2i+5} & & 
\end{array}
\end{displaymath}
where each $x_j$ is placed beneath the vertex it belongs to, and 
$v_{2i}=v_{2i+1}=v_{2i+4}=v_{2i+5}$.\\
Remove $x_{2i+1},x_{2i+2},x_{2i+3},x_{2i+4}$ from the sequence, noting that
$v_{2i}=v_{2i+5}$, so the resulting sequence has the same form, and label the new
sequence $y_1,...,y_{2n-4}$.\\
It now follows that $\sum_{j=1}^{2n-5}d(y_j,y_{j+1})$ is the same as 
$\sum_{j=1}^{2n-1}d(x_j,x_j+1)$, except that $d(x_{2i},x_{2i+1})+d(x_{2i+1},x_{2i+2})+
d(x_{2i+2},x_{2i+3})+d(x_{2i+3},x_{2i+4})+d(x_{2i+4},x_{2i+5})$ is replaced by
$d(x_{2i},x_{2i+5})$; the latter is smaller by the triangle inequality, which completes
the proof of the claim.
\end{proof}

We now have the necessary machinery to prove the main result.

\begin{fprod} \label{fprod}
Let $G$, $H$ have property A, and let $K$ be a common subgroup.  Then the
amalgamated free product $F=G\ast_K{}H$ has property A.
\end{fprod}

\begin{proof}
Note first that $F$ acts freely by isometries on $X_F$ as above.  Hence the orbit map
$F\to{}X$, $f\mapsto{}fx_0$ for some fixed $x_0\in{}X_F$ is an injection.  In fact, one can define a length 
function $|.|'$ on $F$ by setting $|f|'=d(x_0,fx_0)$.  This is symmetric, as $F$ acts by isometries
and $d$ is symmetric, and satisfies 
$|gh|'=d(x_0,ghx_0)\leq{}d(x_0,gx_0)+d(gx_0,ghx_0)=d(x_0,gx_0)+d(x_0,hx_0)=|g|'+|h|'$ and
$|g|'=0$ if and only if $g=e$ as $F$ acts freely.  Moreover,
it is proper as $F$ acts freely and $X_F$ is locally finite, whence there are
at most finitely many translates $fx_0$ within any ball about $x_0$ in $X_F$.  Hence
$d'(g,h)=|g^{-1}h|'$ defines a left-invariant, bounded geometry metric on $F$, which is coarsely
equivalent to the original metric by \ref{wmce}.\\
It will thus suffice to prove that $X_F$ has the property \ref{equiv}, (6); as $F$ is coarsely
equivalent to a bounded geometry (as $F$ is bounded geometry) subspace of this, it will follow 
that $F$ has property A.

Now, as the map $p:X_F\to\mathcal{T}_F$ defined above is contractive, and $\mathcal{T}_F$ has 
property A by \ref{trees}, it will be sufficient (by \ref{fibreing}) to prove that for any cover
$\{\bar{B}(v,n):v\in{}V\}$ of $\mathcal{T}_F$, the pullbacks $p^{-1}(\bar{B}(v,n))$ form an
equi-property A family; as for any fixed $v_0$, all $v\in{}V$, $p^{-1}(\bar{B}(v,n))$ is isometric to
$p^{-1}(\bar{B}(v_0,n))$, it suffices to prove that this has the property \ref{equiv}, (6).\\
Further, $p^{-1}(\bar{B}(v_0,n))$ is the union of all the cosets $v$ (as subspaces of $X_F$) such that 
$d_{\mathcal{T}_F}(v_0,v)\leq{}n$; note that
each such coset is isometric to $G$ or $H$, so has property A by assumption.

To prove that $p^{-1}(\bar{B}(v_0,n))$ has \ref{equiv} (6), then, proceed by induction on $n$.  For the base case,
$p^{-1}(\bar{B}(v_0,0))=v_0$, which is isometric to either $G$ or $H$, so has
property A.\\
For the inductive step, assume that $p^{-1}(\bar{B}(v_0,n-1))$ has property A, and note that
\begin{displaymath}
p^{-1}(\bar{B}(v_0,n))=p^{-1}(\bar{B}(v_0,n-1))\cup\bigcup_{d_{\mathcal{T}_F}(v_0,v)=n}v;
\end{displaymath}
we will use \ref{unions2} to prove that this has the property \ref{equiv} (6).

Let $L>0$.\\
For each $v\in{}V$ such that $d_{\mathcal{T}_F}(v_0,v)=n$, let $\tilde{v}$ consist of those
elements of $v$ at distance no more than $L/2$ from $p^{-1}(\bar{B}(v_0,n-1))$, and let
\begin{displaymath}
Y=p^{-1}(\bar{B}(v_0,n-1))\cup\bigcup_{d_{\mathcal{T}_F}(v_0,v)=n}\tilde{v}.
\end{displaymath}
It is coarsely equivalent to $p^{-1}(\bar{B}(v_0,n-1))$, so has \ref{equiv} (6) by induction.\\
Moreover, the family $\{v\backslash\tilde{v}:d_{\mathcal{T}_F}(v_0,v)\}$ is $L$-separated (in
$p^{-1}(\bar{B}(v_0,n))$; this is false in $X_F$). This follows from the explicit form for
$d$ given in \ref{metrictree}.  For if $x_v\in{}v\backslash\tilde{v}$, $y_w\in{}w\backslash\tilde{w}$ for
$v,w$ distinct elements of $\{v\in\mathcal{T}_F:d(v,v_0)=n\}$, then any sequence of
elements as in \ref{metrictree} must pass through $p^{-1}(\bar{B}(v_0,n-1))$ (as 
$\mathcal{T}_F$ is a tree), whence $d(x_v,y_w)\geq{}L$.\\
The result now follows from \ref{unions2}.
\end{proof}

\begin{hnn}
Let $G$ be a property A group, and $\alpha:H\to{}K$ an isomorphism between two
subgroups of $G$.  Then the HNN extension $G*_\alpha$ has property A.
\end{hnn}

\begin{proof}
HNN extensions are built up from direct limits, free products with amalgamation and a
semi-direct product with $\mathbb{Z}$ (which has A), so this is immediate from the above.
\end{proof}

\begin{cecprop}
Almost all of the closure properties above also hold for (countable, discrete) groups that are
coarsely embeddable in Hilbert space; see \cite{dg} and \cite{dg2}.  The exception is the
theorem on extensions, where the best known result is that if
$1\to{}K\to{}G\to{}H\to{}1$ is short exact, $K$ is coarsely embeddable
and $H$ has A, then $G$ is also coarsely embeddable.

It is also worth pointing out that property A is (probably) not closed under quotients.
This would follow from a claim of Gromov \cite{gr2} (briefly discussed in section \ref{nonagp}) 
that there exist finitely 
generated groups that do not have property A; any such group is of course 
a quotient of some $F_n$.  As quotients are a special case of general inductive limits, property A is (probably)
not closed under these either.
\end{cecprop}

We conclude this section with a list of some classes of groups known to have property A
(a class somewhat larger than those that can be built out of free groups using 
the properties in this section!).  There are also criteria for property A 
that use some of the techniques of geometric group theory; see \cite{hr} (we look at this in section 
\ref{aamen}), \cite{bell2}, \cite{bell}, \cite{cn}, \cite{gu} and \cite{tu}, and 
the notes on each in chapter \ref{bib}.

\begin{itemize}
\item Amenable groups; we prove this in section \ref{amenableintro}.
\item Groups of finite asymptotic dimension; we proved this in 
\ref{fad}.
\item Groups with a one-relation presentation, \cite{gu}.
\item Groups that are \emph{hyperbolic in the sense of Gromov}. This follows from the above and
the result of \cite{roe4} (due originally to Gromov) that such groups have finite asymptotic dimension.
\item Groups that are hyperbolic relative to a collection of property A subgroups; this is proved in
\cite{oz2}.
\item Finitely generated right-angled Artin groups and
finitely generated Coxeter groups. This follows from \cite{cn}, where it is proved that
groups that act properly cocompactly by isometries on a (finite dimensional) \emph{CAT(0) cube complex} have property A.
\item  Any discrete subgroup of a connected Lie group.  This follows as such a group
admits an amenable action on a compact homogeneous space arising as a quotient of the
the original group; see \cite{adr}, 2.2.14 and 5.1.1.
\item Any countable subgroup of $GL(n,F)$ for any $n$ and any field $F$.  This is a result
of E. Guentner, N. Higson and S. Weinberger, \cite{ghw}. 
\end{itemize}

\pagebreak

\section{Kernels and coarse embeddability} \label{kerembed}

\subsection{Introduction}\label{kerembedintro}

One motivation for this chapter is that coarse embeddability into certain Banach spaces has various implications for a metric space. Foremost for us is the result of G. Yu \cite{y1} that coarse embeddability into
Hilbert space implies coarse Baum-Connes, but see also \cite{ky}, which does something similar for more general uniformly convex Banach spaces. A further motivation is simply that coarse embeddability into Hilbert space has a fairly close relationship with property A (and therefore with everything else in this piece).  Finally, another active area of research studies the distortion of embeddings into Banach spaces from a quantitative point of view; we will not look at this here (such questions require a little more precision that our purely coarse methods allow), but see \cite{tes} for some applications of property A in this context.

The chapter has two main sections:

Section \ref{embedhil} introduces the theory of \emph{kernels}.  For us, a kernel is a map $X\times{}X\to\mathbb{R}$ 
(or $\mathbb{C}$), where $X$ is some metric space.  Certain types of these are linked to 
the $C^*$-algebras $C^*_u(X)$ and $C^*_r(G)$ (chapter \ref{analysis}), amenability 
(chapter \ref{amenable}), and also property A, so tie together a lot of ideas in this piece.  There is also a
classical theory linking kernels to embeddings in Hilbert space; we formulate and prove certain results
connecting kernels and coarse embeddability, which in particular clarify the relationship of the latter
to property A.  A good general reference for material on positive and negative type kernels in chapter 8 of \cite{bl}.

Section \ref{lp} starts with some more classical results on kernels, which
are used to relate coarse
embeddability in Hilbert space to coarse embeddability in $l^p$ spaces.  This is due to P. Nowak, \cite{no1}, \cite{no2}. 
The main results of this section
are that $l^p$ coarsely embeds in $l^2$ for $0<p\leq2$, and that $l^2$ coarsely embeds
in $l^p$ for $1\leq{}p\leq\infty$.  The former is a direct application of the theory
of kernels we initially develop (and was proved in a different context by
\cite{bdk}), the second a generalisation of a result
of M. Dadarlat and E. Guentner \cite{dg} characterising coarse embeddability 
(which we have reproduced in \ref{eh2}, (2)).

One expects coarse embeddability into $l^p$ to get progressively easier as $p$ increases;
indeed any space $X$ coarsely embeds in $l^{\infty}(X)$ (see e.g. \cite{dr}).  
These results show there is
a little `leeway' in this.  They are also of interest as it was once conjectured
(\cite{dr}, 4.4) that property A (for a separable space) is equivalent to coarse 
embeddability in $l^1$; the results of section \ref{nonace} 
(also due to Nowak, \cite{no3}) show this to be false.  
Further, there is a version of the coarse Baum-Connes 
conjecture related to $l^1$ (this is related to the work of V. Lafforgue on so-called
\emph{Banach KK-theory}); one might expect that
$l^1$ coarse embeddability has a similar relationship to this as $l^2$ coarse
embeddability does to the usual `$l^2$ coarse Baum-Connes conjecture'.  Note that $l^1$ is \emph{not}
uniformly convex, so the results of \cite{ky} do not apply to spaces that coarsely embed in it.

\subsection{Coarse embeddability into Hilbert space} \label{embedhil}

We start with the basic definitions.

\begin{kernels} \label{kernels}
Let $X$ be a metric space.

A kernel is a map $k:X\times{}X\to\mathbb{C}$ or $\mathbb{R}$.\\
 
$k$ is \emph{self-adjoint} (or \emph{symmetric} if real-valued) if 
$k(x,y)=\overline{k(y,x)}$.\\

$k$ has $(R,\epsilon)$ \emph{variation} if $d(x,y)\leq{}R$ implies 
$|k(x,y)-1|<\epsilon$.\\

$k$ has \emph{finite propagation} if there exists $S>0$
such that $d(x,y)>S$ implies $k(x,y)=0$ .  The \emph{propagation} of $k$ is
the smallest such $S$.\\ 

$k:X\times{}X\to\mathbb{C}$ (or $\mathbb{R}$) is of \emph{positive type} if 
for all finite sequences
$x_1,...,x_n$ of elements of $X$ and 
$\lambda_1,...,\lambda_n$ of complex (or real) numbers,
\begin{displaymath} 
\sum_{i,j=1}^n\lambda_i\overline{\lambda_j}k(x_i,x_j)\geq{}0.
\end{displaymath}

$k:X\times{}X\to\mathbb{R}$ is of \emph{negative type} if for all 
finite sequences
$x_1,...,x_n$ of elements of $X$ and $\lambda_1,...,\lambda_n$ of
real numbers such that $\sum_{i=1}^n\lambda_1=0$,
\begin{displaymath}
\sum_{i,j=1}^n\lambda_i\lambda_jk(x_i,x_j)\leq{}0.
\end{displaymath}

A positive (negative) type kernel is \emph{normalised} if for all $x\in{}X$, 
$k(x,x)=1$ ($k(x,x)=0$).
\end{kernels}

The following records some basic facts about kernels.

\begin{kernelrem} \label{kernelrem}
Let $X$ be a set.
\begin{enumerate}
\item If $k$ is of positive type, then $-k$ is of negative type.
\item Kernels of either type are preserved by addition, multiplication by
non-negative scalars and pointwise limits.
\item A constant map to a non-negative scalar is a kernel of positive type.
\item A real-valued constant map is a kernel of negative type.
\item A map of the form $k:(x,y)\mapsto{}f(x)\overline{f(y)}$  where
$f:X\to\mathbb{C}$ is a kernel of positive type.
\item $k$ is a kernel of positive type on $X$ if and only if for all finite
sequences $x_1,...x_n$, the matrix $[k(x_i,x_j)]_{i,j=1}^n$ is positive. 
\end{enumerate}
\end{kernelrem}

The proofs are simple calculations; for example for (5)
\begin{displaymath}
\sum_{i,j}\lambda_i\overline{\lambda_j}k(x_i,x_j)
=\sum_{i,j}\lambda_i\overline{\lambda_j}f(x_i)\overline{f(x_j)}
=\left(\sum_i\lambda_if(x_i)\right)\overline{\left(\sum_i\lambda_if(x_i)\right)}
\geq0.
\end{displaymath}
Note that the relationship between kernels of
negative and positive type is more subtle than a simple change of sign; indeed there
are non-zero maps that are both.

The following examples are our most important (indeed the
theorem following shows that they are in a sense the \emph{only}
examples).

\begin{kerex} \label{kerex}
Let $\mathcal{H}$ be any (real or complex) Hilbert space, with inner product
$\langle,\rangle$ and norm $\|.\|$.  Let $X$ be a metric space and 
$f:X\to\mathcal{H}$ be any map.
\begin{itemize}
\item The map $(x,y)\mapsto\|f(x)-f(y)\|^2$ defines a normalised negative type 
kernel on $X$.
\item The map $(x,y)\mapsto\langle{}f(x),f(y)\rangle$ defines a positive type kernel on
$X$.  It is normalised if each $f(x)$ is a unit vector.
\end{itemize}
\end{kerex}

\begin{proof}
We only prove the first; the second is similar (and essentially the same as the proof of
(5) above, which is the one-dimensional case).\\
Let $x_1,...,x_n\in{}X$, and $\lambda_1,...,\lambda_n\in\mathbb{R}$ be 
such that $\sum\lambda_i=0$. Then
\begin{align*}
& \sum_{i,j}\lambda_i\lambda_j\|f(x_i)-f(x_j)\|^2
=\sum_{i,j}\lambda_i\lambda_j(\|f(x_i)\|^2+\|f(x_j)\|^2-
2\textrm{Re}\langle{}f(x_i),f(x_j)\rangle) \\
& ~=\sum_i\lambda_i\|f(x_i)\|^2\left(\sum_j\lambda_j\right)
+\sum_j\lambda_j\|f(x_j)\|^2\left(\sum_i\lambda_i\right)
-2\textrm{Re}\sum_{i,j}\langle{}f(x_1),f(x_j)\rangle \\
& ~=-2\textrm{Re}\left\langle\sum_i\lambda_if(x_i),
\sum_i\lambda_if(x_i)\right\rangle\leq0,
\end{align*}
where the last equality is by assumption on the $\lambda_i$.\\
Normalisation is immediate.
\end{proof}

The next theorem shows that these two examples are `universal' in the sense that for
any kernel of positive or negative type, there exists a Hilbert space such that it
arises in one of these ways.  It thus gives a clear link between kernels and Hilbert space.
All of the results on kernels in this section (that are not specifically related
to coarse embeddings) are from the classic papers \cite{sch1} and \cite{sch2}.

\begin{kerthe} \label{kerthe}
\begin{itemize}
\item If $k:X\times{}X\to\mathbb{R}$ is a symmetric normalised kernel of negative type, 
then there exists a real Hilbert space $\mathcal{H}$ and a map 
$f:X\to\mathcal{H}$ such that $k(x,y)=\|f(x)-f(y)\|^2$.
\item If $k:X\times{}X\to\mathbb{C}$ ($\mathbb{R}$) is a self-adjoint (symmetric)
kernel of positive type, 
then there exists a complex (real) Hilbert space $\mathcal{H}$ and a map 
$f:X\to\mathcal{H}$ such that $k(x,y)=\langle{}f(x),f(y)\rangle$.
\end{itemize}
\end{kerthe}

\begin{proof}
We start with the case where $k$ is negative type.  Let $C_c^{(0)}(X)$ be the vector space of finitely supported functions
$f:X\to\mathbb{R}$ that satisfy $\sum_{x\in{}X}f(x)=0$, and define
\begin{displaymath}
\langle{}f,g\rangle=-\frac{1}{2}\sum_{x,y\in{}X}k(x,y)f(x)g(y).
\end{displaymath}
This is clearly bilinear, and is symmetric as $k$ is.  Moreover,
\begin{displaymath}
\langle{}f,f\rangle=-\frac{1}{2}\sum_{x,y\in{}X}k(x,y)f(x)f(y)\geq0,
\end{displaymath}
as $k$ is of negative type.\\
By these properties, it satisfies the Cauchy-Schwarz inequality (proof:
for any real $t$, and $f,g\in{}C_c^{(0)}(X)$, 
$0\leq\langle{}tf+g,tf+g\rangle=t^2\langle{}f,f\rangle+2t\langle{}f,g\rangle
+\langle{}g,g,\rangle$; the discriminant associated to this (real) quadratic is thus 
non-positive, which is exactly the Cauchy-Schwarz inequality).\\
It follows that $E=\{f\in{}C_c^{(0)}(X):\langle{}f,f\rangle=0\}$ is a subspace, so
we can take the quotient $C_c^{(0)}(X)/E$; the map induced on this by
$\langle,\rangle$ is a genuine inner product.\\
Complete in the induced norm to get a Hilbert space $\mathcal{H}$.\\
Finally, pick any point $x_0\in{}X$ and define $f:X\to{}\mathcal{H}$ by
setting $f(x)=\delta_x-\delta_{x_0}$.  Then
\begin{displaymath}
\|f(x)-f(y)\|^2=\langle\delta_x-\delta_y,\delta_x-\delta_y\rangle
=k(x,y)-\frac{1}{2}k(x,x)-\frac{1}{2}k(y,y)=k(x,y),
\end{displaymath}
as $k$ is normalised.

For the case of $k$ positive type, the proof is essentially the same, except we start by setting
$C_c(X)$ to be the space of all finitely supported functions on $X$, and
defining a pseudo-inner product by
\begin{displaymath}
\langle{}f,g\rangle=\sum_{x,y\in{}X}k(x,y)f(x)g(y).
\end{displaymath}
\end{proof}

These ideas are used (more or less implicitly) in many connected areas, such as the
GNS construction giving the existence of faithful representations of $C^*$-algebras, and
the equivalence of group a-T-menability and the existence of a proper isometric action on 
Hilbert space.  See \ref{stinespring} and \ref{replem} for some more applications.

Note that a normalised, positive type
complex-valued kernel is self-adjoint and a normalised negative type kernel is 
symmetric.  For example in the positive type case, substituting in 
$(\lambda_x,\lambda_y)=(1,1),(1,i)$ respectively give that
$2+k(x,y)+k(y,x)\geq{}0$ and $2-ik(x,y)+ik(y,x)\geq0$. Hence the imaginary
parts of $k(x,y)$, $k(y,x)$ cancel and the real parts are the same.
There is thus some redundancy in the above.

\begin{mokerex} \label{mokerex}
Consider the following positive type kernel on a set $X$: 
\begin{displaymath}
k(x,y)=\left\{\begin{array}{ll}
1, & x=y \\
0, & \textrm{otherwise.}\end{array}\right.
\end{displaymath}
The construction above yields $l^2(X)$ with the usual inner product.

On the other
hand, consider the constant function
$\mathbf{1}$ on $X\times{}X$ (which is also a positive type kernel).  The corresponding
pseudo-inner product is (of course!) highly singular, and the resulting Hilbert space (after 
taking the quotient by functions of length zero) is just $\mathbb{C}$.
\end{mokerex}

Before proving the main theorem of this section, we mention some links
between real and complex valued kernels of positive type, and between real and complex 
Hilbert spaces.  This shows that our occasional expedient jumps
between real and complex are not cheating!

\begin{krecom} \label{krecom}
A symmetric, positive type, (normalised) real-valued kernel $k$ on $X$ is also a
self-adjoint, positive type, (normalised) complex-valued kernel, 
i.e. if $\sum_{i,j}\lambda_i\lambda_jk(x_i,x_j)\geq{}0$
for all sequences $(\lambda_i)_{i=1}^n$ of real numbers, then also
$\sum_{1,j}\lambda_i\overline{\lambda_j}k(x_1,x_j)\geq{}0$
for all $n$-tuples $(\lambda_i)_{i=1}^n$ of complex numbers.
\end{krecom}

To prove it, just note that the imaginary parts of
$\sum_{1,j}\lambda_i\overline{\lambda_j}k(x_1,x_j)$ cancel. 

\begin{recom} \label{recom}
$X$ coarsely embeds into a real Hilbert space if and only if it coarsely embeds 
into a complex Hilbert space.
\end{recom}

\begin{proof}
The identities
\begin{align*}
\langle{}x,y\rangle_{\mathbb{R}} & =\textrm{Re}\langle{}x,y\rangle \\
\langle{}x,y\rangle_{\mathbb{C}} & =\frac{1}{2}\langle{}x,y\rangle
+\frac{1}{2}\langle{}ix,iy\rangle+\frac{1}{2}i\langle{}x,iy\rangle
-\frac{1}{2}i\langle{}ix,y\rangle
\end{align*}
convert a complex inner product into a real one, and a real inner product into a complex
one respectively.  Note that the former always works, as a complex vector
space is automatically also real; in order to make the second applicable, we need to
expand the space to include all imaginary multiples of a given othonormal basis.
\end{proof}

The following is our main result of this section.

\begin{eh2} \label{eh2}
Let $X$ be any metric space.  The following are equivalent:
\begin{enumerate}
\item $X$ is coarsely embeddable in Hilbert space.
\item For every $R,\epsilon>0$, there exists a (real) Hilbert space $\mathcal{H}$ and a
map $\eta:X\to\mathcal{H}$ that satisfies:
\begin{enumerate}
\item $\|\eta_x\|=1$ for all $x\in{}X$;
\item $\eta_x$ has $(R,\epsilon)$ variation (see definition \ref{equivdef});
\item
$\lim_{S\to\infty}\sup\{|\langle\eta_x,\eta_y\rangle|:d(x,y)\geq{}S\}=0$
\end{enumerate}
\item For all $R,\epsilon>0$, there exists a normalised symmetric kernel 
$k:X\times{}X\to\mathbb{R}$ of positive type and $(R,\epsilon)$ variation, and 
such that $\lim_{S\to\infty}\sup\{|k(x,y)|:d(x,y)\geq{}S\}=0$.
\item There exists a symmetric, normalised negative type kernel $k$ on $X$ 
and maps 
$\rho_1,\rho_2:\mathbb{R}^+\to\mathbb{R}^+$ such that $\rho_1(t)\to+\infty$
as $t\to+\infty$ and
\begin{displaymath}
\rho_1(d(x,y))\leq{}k(x,y)\leq\rho_2(d(x,y))
\end{displaymath}
(such a kernel is sometimes called \emph{effective}).
\end{enumerate}
\end{eh2}

\begin{proof}
(1) implies (2) (\cite{dg}, proposition 2.1):\\
Let $f:X\to\mathcal{H}$ be a coarse embedding of $X$ into a Hilbert $\mathcal{H}$, so there
exist functions $\rho_1,\rho_2:\mathbb{R}^+\to\mathbb{R}^+$ satisfying the
properties in \ref{cembed}.  By the lemma, we can assume that $\mathcal{H}$ is
real.  Define:
\begin{displaymath}
\textrm{exp}(H)=\mathbb{R}\oplus{}\mathcal{H}\oplus(\mathcal{H}\otimes{}\mathcal{H})\oplus
(\mathcal{H}\otimes{}\mathcal{H}\otimes{}\mathcal{H})\oplus...
\end{displaymath}
and define $\exp:X\to\exp(\mathcal{H})$ by:
\begin{displaymath}
\textrm{exp}(h)=1\oplus{}h\oplus\left(\frac{h}{\sqrt{2!}}\otimes{}h\right)
\oplus\left(\frac{h}{\sqrt{3!}}\otimes{}h\otimes{}h\right)\oplus...
\end{displaymath}
Now, $\langle\exp(h),\exp(h')\rangle=e^{\langle{}h,h'\rangle}$ for all
$h,h'\in{}H$.  For any $t>0$ define:
\begin{displaymath}
\eta_x^t=e^{-t\|f(x)\|^2}\exp(\sqrt{2t}f(x)).
\end{displaymath}
Note then that
\begin{displaymath}
\langle\eta_x^t,\eta_y^t\rangle
=e^{-t(\|f(x)\|^2+\|f(y)\|^2)}e^{2t\langle{}f(x),f(y)\rangle}
=e^{-t\|f(x)-f(y)\|^2},
\end{displaymath}
whence $\|\eta_x^t\|=1$ for all $x,t$, and by assumption on $f$,
\begin{displaymath}
e^{-t\rho_2(d(x,y))^2}\leq\langle\eta_x^t,\eta_y^t\rangle
\leq{}e^{-t\rho_1(d(x,y))^2}.
\end{displaymath}
Let $t_0=\frac{\epsilon}{1+\rho_2(R)^2)}$, and set $\eta_x=\eta_x^{t_0}$.\\
We complete the argument with the calculations:
\begin{displaymath}
\sup_{d(x,y)\geq{}S}|\langle\eta_x,\eta_y\rangle|
=\sup_{d(x,y\geq{}S}e^{-t_0\|f(x)-f(y)\|^2}\to0
\end{displaymath}
as $S\to\infty$, as $f$ is effectively proper.  Also, if $d(x,y)<R$
\begin{displaymath}
\|\eta_x-\eta_y\|^2=2-2\langle\eta_x,\eta_y\rangle
\leq2\left(1-\exp\left(\frac{-\epsilon\rho_2(d(x,y))}{1+\rho_2(R)^2}\right)\right)
\leq2(1-e^{-\epsilon})<2\epsilon.
\end{displaymath}

(2) implies (3):\\
Let $\eta$ be as in (2) and define $k(x,y)=\langle\eta_x,\eta_y\rangle$.\\
(3) is now immediate from the relationship: 
$\|\eta_x-\eta_y\|^2=\|\eta_x\|^2+\|\eta_y\|^2-2\langle\eta_x\eta_y\rangle$.\\

(3) implies (4) (\cite{roe1}, 11.16):\\
By (3), for each $n\in\mathbb{N}^+$ we can find $R_n>\max\{n,R_{n-1}\}$ and 
normalised self-adjoint kernels
$k_n$ of positive type satisfying
\begin{displaymath}
|k_n(x,y)-1|<2^{-n} ~\textrm{for} ~d(x,y)<n; ~|k_n(x,y)|<\frac{1}{2}
~\textrm{for} ~d(x,y)>R_n.
\end{displaymath}
Define
\begin{displaymath}
k(x,y)=\sum_{n=1}^{\infty}(1-k_n(x,y)),
\end{displaymath}
which exists, as the sum has only finitely many terms larger than $2^{-n}$ for
each pair $(x,y)$.\\
We note the following two properties of $k$:
\begin{itemize}
\item If $d(x,y)<n$, we get that at most $n$ terms in the sum above are greater
than $2^{-n}$; the remaining terms are bounded by $1$ as each $k_n$ is
normalised.  Hence $|k(x,y)|\leq{}2d(x,y)+1$.
\item If $Q_t=\min\{n:t\leq{}R_n\}$ for $t\in\mathbb{R}^+$, and $d(x,y)<Q_t$,
then at least $Q_t$ terms in the sum above are at least $\frac{1}{2}$.
Hence $|k(x,y)|\geq{}\frac{1}{2}R(d(x,y))$.
\end{itemize} 
Note also that $k$ is a pointwise limit of linear combinations of negative type
kernels, so of negative type; it is normalised as each $k_n$ is.\\
If we set $\rho_1(t)=\frac{1}{2}Q_t$, $\rho_2(t)=2t+1$, we get (4).\\

(4) implies (1):\\
The construction of \ref{kerthe} builds a Hilbert space $\mathcal{H}$ and a map
$f:X\to\mathcal{H}$ such that $k(x,y)=\|f(x)-f(y)\|^2_{\mathcal{H}}$.\\
It follows immediately that
\begin{displaymath}
\sqrt{\rho_1(d(x,y))}\leq\|f(x)-f(y)\|^2\leq\sqrt{\rho_2(d(x,y))},
\end{displaymath}
whence $f$ is a coarse embedding.
\end{proof}

\begin{eh2rem} \label{eh2rem}
Note the striking similarity between (2) and (3) in the proposition above and 
(7) and (8) respectively from theorem \ref{equiv} (whether the Hilbert space is
real or complex in either case is irrelevant).
This gives an alternative proof that property A implies coarse
embeddability in a Hilbert space.

It is also partially responsible for the oft-repeated slogan 
`property A is to coarse embeddability as amenability is to a-T-menability'; 
we will return to this in chapter \ref{amenable}.  Alternatively, one could
say that property A is a `bounded support' property, versus coarse
embeddability, which is a `$C_0$ property'.

Note also that one can replace (2) (c) above by
\begin{displaymath} 
\lim_{S\to\infty}\inf\{\|\eta_x-\eta_y\|:d(x,y)\geq{}S\}=\sqrt2,
\end{displaymath}
but the form we used is more obviously related to \ref{equiv}, (7).  
We will generalise
this alternative form in the next section; see \ref{geneh2}. 
\end{eh2rem}

We conclude with an application of property (2) above.

\begin{ceunions} \label{ceunions}
Say $X=Y\cup{}Z$ is a metric space, and that both $Y$, $Z$ coarsely embed.  Then 
$X$ coarsely embeds.
\end{ceunions}

It is possible to give a proof of this fact similar in spirit to that of \ref{unions}; see
\cite{dg2}.  For variety we sketch a different (admittedly longer) approach.

\begin{proof} (sketch).
By replacing $Y,Z$ with $Y, Y\backslash{}Z$, and noting that coarse embeddability
is inherited by subspaces, assume $Y$, $Z$ are disjoint.\\
Assume also that they are both non-empty.\\
Let $R,\epsilon>0$, and assume without loss of generality that $\epsilon<1$.\\
Let $c=8R/\epsilon$.\\
Now, by coarse embeddability for $Y,Z$ there exist Hilbert spaces $\mathcal{H}_Y$, 
$\mathcal{H}_Z$ and maps
\begin{align*}
\eta^Y & :Y\to\mathcal{H}_Y, \\
\eta^Z & :Z\to\mathcal{H}_Z
\end{align*}
satisfying the conditions above with $(R+4c,\epsilon/8)$ variation.\\
For $t\in\mathbb{R}$, let $\{t\}=\max\{0,\min\{t,1\}\}$.\\
Let $p:X\to{}Z$, $q:X\to{}Y$ be projections such that
$d(x,px)\leq2d(x,Z)$, $d(x,qx)\leq2d(x,Y)$ for all $x\in{}X$.\\
Define $\tilde{\eta}:X\to\mathcal{H}_Y\oplus\mathcal{H}_Z$ by setting
\begin{displaymath}
\tilde{\eta}_x=\left\{\frac{d(x,Y)}{c}\right\}\eta^Z_{px}+
\left\{1-\frac{d(x,Y)}{c}\right\}\eta^Y_{qx}.
\end{displaymath}
Finally, set $\eta_x=\tilde{\eta}_x/\|\tilde{\eta}_x\|_2$, which we claim has the
properties above.  While elementary, the estimates required are a little lengthy; we hope the above gives an idea of how a proof might proceed.
\end{proof}

\subsection{Coarse embeddability into $l^p$ spaces} \label{lp}
As mentioned in the introduction to this section, we prove 
the following two results:
\begin{itemize}
\item $l^p$ coarsely embeds in $l^2$ for $0<p<2$.
\item $l^2$ coarsely embeds in $l^p$ for $1\leq{}p<\infty$.
\end{itemize}

We start with the first, which makes use of some of the machinery of positive and
negative type kernels introduced in \ref{kernels}, introducing a few more of
Schoenberg's results (\cite{sch1} and \cite{sch2}) along the way.

\begin{schur} \label{schur}
The pointwise (Schur) product of kernels $k$, $l$ of positive type is also of positive 
type.
\end{schur}

\begin{proof}
Let $m(x,y)=k(x,y)l(x,y)$ be the pointwise product.\\
Let ${x_1,...,x_n}\subseteq{}X$ and let $M$ be the $n\times{}n$ matrix
$[m(x_p,x_q)]$.\\
Then the matrices $K=[k(x_p,x_q)]$, 
$L=[l(x_p,x_q)]$ are positive, whence there exists $n\times{}n$ $A$ with $L=A^*A$, 
i.e. $l(x_p,x_q)=\sum_ra_{pr}\overline{a_{rq}}$.\\
Then for any $\mathbf{\lambda}=(\lambda_1,...,\lambda_n)^T\in\mathbb{C}^n$, we get that
\begin{displaymath}
\mathbf{\lambda}^*M\mathbf{\lambda}=
\sum_r\sum_{p,q}k(x_p,x_q)(\overline{a_{qr}\lambda_q})(a_{pr}\lambda_p)\geq0;
\end{displaymath}
by the comment \ref{kernelrem}, (6), this finishes the proof.
\end{proof}

Note that it follows that if $f(x)=\sum{}a_nx^n$ is a power series converging pointwise
on the range of positive type $k$, and $a_n\geq0$ for all $n$, then $f(k)$ is of positive
type too. 

The next result is often referred to as `Schoenberg's lemma'.

\begin{expneg} \label{expneg} 
Let $X$ be a set.\\
Then $k:X\times{}X\to\mathbb{R}$ is a symmetric negative type kernel on $X$ if and 
only if $(x,y)\mapsto{}e^{-tk(x,y)}$ is symmetric of positive type for all $t\geq0$.
\end{expneg}

\begin{proof}
Say first that $e^{-tk}$ is of positive type for all $t\geq0$.\\
Then by the comments following \ref{kernels}, $1-e^{-tk}$ is of
negative type, whence so too is
\begin{displaymath}
\lim_{t\to{}0}\frac{1-e^{-tk}}{t}=k.
\end{displaymath}
Conversely, say $k$ is of negative type and fix $p\in{}X$.\\
Set $l(x,y)=-k(x,y)+k(x,p)+k(p,y)-k(p,p)$, which an easy calculation shows to be
of positive type.\\
Hence for $t\geq0$, we have  that 
$e^{-tk(x,y)}=e^{tl(x,y)}(e^{-tk(x,p)}e^{-tk(p,y)})e^{tk(p,p)}$.\\
The first term is of positive type by the comment following \ref{schur}; the
middle two terms together are of positive type by symmetricity of $k$ and 
\ref{kernelrem}, (5); the fourth term is a positive constant.\\
Hence $e^{-tk}$ is of positive type by \ref{schur}.
\end{proof}

\begin{schoen}
Let $k$ be a negative type kernel on $X$ such that $k(x,y)\geq0$ for all $x,y\in{}X$.\\
Let $0<\alpha<1$.  Then $k^{\alpha}$ is also of negative type.
\end{schoen}

\begin{proof}
By assumption on $k$, $1-e^{-tk}\geq0$ for all $t\geq0$; this is also a kernel of
negative type by Schoenberg's lemma and the comments in \ref{kernelrem}.
We also have that for any $x\geq0$ and any $0<\alpha<1$,  
$c_{\alpha}=(\int_0^{\infty}(1-e^{-t})t^{\alpha-1}dt)^{-1}$ exists and is
positive, and
\begin{displaymath}
x^{\alpha}=c_{\alpha}\int_0^{\infty}(1-e^{-tx})t^{-\alpha-1}dt.
\end{displaymath}
The result follows.
\end{proof}

Finally, this allows us to prove

\begin{pto2}[\cite{no1}, 4.1] \label{pto2}
$l^p$ embeds in $l^2$ for $0<p\leq{}2$. 
\end{pto2}
\begin{proof}
Note first that the map $k:\mathbb{R}\times\mathbb{R}\to\mathbb{R}^+$ given by
$k(x,y)=|x-y|^2$ is a negative type kernel.  This follows as a special
(one-dimensional) case of \ref{kerex}.\\
Hence by the previous lemma, $k^{\alpha}:(x,y)\to{}|x-y|^{2\alpha}$ is a negative type
kernel for all $0<\alpha<1$.\\
Further, if $x=(x_i),y=(y_i)\in{}l^p$, we can sum to get that the map
$k_p:l^p\times{}l^p\to\mathbb{R}$ given by $(x,y)\mapsto\|x-y\|_p^p$ is a 
symmetric, normalised negative-type kernel for all $0<p\leq2$.  Moreover, it satisfies
\begin{displaymath}
\sqrt[p]{d(x,y)}\leq{}k_p(x,y)\leq\sqrt[p]{d(x,y)},
\end{displaymath}
i.e. $k_p$ satisfies condition (3) in theorem \ref{eh2}; as $l^p$ is separable for
$p<\infty$, we may as well assume the Hilbert space it coarsely embeds into is $l^2$.
\end{proof}

We now move on to the second of Nowak's results mentioned in the introduction to this
section.  The first proposition we prove is a generalisation of \ref{eh2}, part (2) 
(more
precisely, with (c) replaced by the variant given in remarks \ref{eh2rem}). 
This is:

\begin{geneh2}[\cite{no2}, theorem 2] \label{geneh2}
Let $X$ be a metric space, and let $1\leq{}p<\infty$.\\
Say there exists $\delta>0$ such that for every $R,\epsilon>0$, there
exists a map $\xi:X\to{}l^p$ satisfying:
\begin{enumerate}
\item $\|\xi_x\|_p=1$ for all $x\in{}X$;
\item $\xi_x$ has $(R,\epsilon)$ variation (see definition \ref{equivdef});
\item $\lim_{S\to\infty}\sup\{\|\xi_x-\xi_y\|_p:d(x,y)\geq{}S\}\geq\delta.$
\end{enumerate}
Then $X$ coarsely embeds in $l^p$.
\end{geneh2}

\begin{proof}
We begin the proof by saying that if $(X_n,\|.\|_n)_{n=1}^{\infty}$ is a 
sequence of Banach spaces, then we define
\begin{displaymath}
(\sum{}X_n)_p=
\left\{(x_n)_{n=1}^{\infty}:x_n\in{}X_n, \sum_{n=1}^{\infty}\|x_n\|^p_n<\infty\right\},
\end{displaymath}
equipped with the norm
\begin{displaymath}
\|(x_n)\|_p=\left(\sum_{n=1}^{\infty}\|x_n\|_n^p\right)^\frac{1}{p}.
\end{displaymath}
We will be interested in the case $(\sum{}l^p)_p$, which is (clearly) isomorphic
to $l^p$.\\
Now, by the assumption, we can construct for each $n$ a map $\xi^n:X\to{}l^p$ that
satisfies:
\begin{enumerate}
\item $\|\xi^n_x\|_p=1$;
\item $\|\xi^n_x-\xi^n_y\|_p<\frac{1}{2^n}$ whenever $d(x_1,x_2)<n$;
\item there exists $S_n$ such that $\|\xi^n_{x_1}-\xi^n_{x_2}\|_p\geq{}\delta$ 
whenever $d(x_1,x_2)\geq{}S_n$.
\end{enumerate}
We may as well assume that the sequence $(S_n)$ is non-decreasing, and tends to infinity
(certainly condition (2) forces $S_n\to\infty$, and we can switch to a 
subsequence if $(S_n)$ is not non-decreasing).\\
Now, choose $x_0\in{}X$ and define a map $f:X\to(\sum{}l^p)_p(\cong{}l^p)$ by
\begin{displaymath}
f(x)=\bigoplus_{n=1}^{\infty}(\xi_x^n-\xi_{x_0}^n)
\end{displaymath}
Note first that $f$ is well-defined, as all but finitely many of the terms
defining $f(x)$, some $x\in{}X$, have $l^p$ norm less than $2^{-n}$.\\
Now, say that $\sqrt[p]m\leq{}d(x,y)\leq\sqrt[p]{m+1}$.  Then
\begin{align*}
\|f(x_1)-f(x_2)\|_p^p & = \sum_{n=1}^{m}\|\xi^n_{x_1}-\xi^n_{x_2}\|_p^p
+\sum_{n=m+1}^{\infty}\|\xi^n_{x_1}-\xi^n_{x_2}\|_p^p \\
& \leq\sum_{n=1}^{m}(\|\xi^n_{x_1}\|_p+\|\xi^n_{x_2}\|_p)^p+1\leq2^pm+1.
\end{align*}
On the other hand, $Q_t:=|\{k\in\mathbb{N}:S_k<t\}|$ is well-defined and
$Q_t\to+\infty$ as $t\to+\infty$.
Moreover if $Q_{d(x_1,x_2)}=q$, then
\begin{displaymath}
\|f(x_1)-f(x_2)\|_p^p\geq\sum_{n=1}^m\|\xi^n_{x_1}-\xi^n_{x_2}\|_p^p\geq{}q\delta^p.
\end{displaymath}
We now have that
\begin{displaymath}
\delta\sqrt[p]{Q_{d(x_1,x_2)}}\leq\|f(x_1)-f(x_2)\|_p
\leq2\sqrt[p]{d(x_1,x_2)+1},
\end{displaymath}
whence $f$ is a coarse embedding.
\end{proof}

Note the similarity between this proof
and that of \ref{yembed}. The next proposition shows that the condition in the 
previous proposition is `robust' with respect to the index $p$.  In order to be
able to prove this we need the classical \emph{Mazur map}.
Recall that $S(l^p)=\{x\in{}l^p:\|x\|_p=1\}$ denotes the unit sphere of $l^p$.

\begin{mazur}
The \emph{Mazur map} $M_{p,q}:S(l^p)\to{}S(l^q)$ is defined by
\begin{displaymath}
M_{p,q}((x_n)_{n=1}^{\infty})=
\left(|x_n|^{\frac{p}{q}}\textrm{sign}(x_n)\right)_{n=1}^{\infty},
\end{displaymath}
where sign$(x_n)=x_n/|x_n|$ if $x_n\not=0$ and zero otherwise.
\end{mazur}

The important property of this map for us is that if $p>q$ then 
there exists some constant 
$C$ depending only on $p/q$ such that $M_{p/q}$ satisfies:
\begin{displaymath}
\frac{p}{q}\|x-y\|_p\leq\|M_{p,q}(x)-M_{p,q}(y)\|_q\leq{}C\|x-y\|_p^{p/q}
\end{displaymath}
(in the language of Banach space theory, $M_{p,q}$ is a \emph{uniform homeomorphism}).  If
$p<q$, the opposite inequalities hold.\\
For a proof of this, and more details about Mazur maps, see \cite{bl}, section 9.1.

\begin{robgeneh2}[\cite{no2}, theorem 3]
If $X$ is a metric space having the property in \ref{geneh2} with respect to some 
$1\leq{}p<\infty$, then it has it with respect to all such $p$.
\end{robgeneh2}

\begin{proof}
Let $R,\epsilon>0$, and say $X$ satisfies the conditions in \ref{geneh2} for some
$\delta$; let $\xi:X\to{}l^p$ satisfy \ref{geneh2} for $R,\epsilon,S$.\\
Define $\xi':X\to{}l^q$ by $\xi'=M_{p,q}\circ\xi$.\\
Say first that $p>q$. The inequality stated above translates to
\begin{displaymath}
\frac{p}{q}\|\xi_x-\xi_y\|_p\leq\|\xi'_x-\xi'_y\|_q\leq{}C\|\xi_x-\xi_y\|_p^{p/q}.
\end{displaymath}
It follows immediately that:
\begin{itemize}
\item if $d(x,y)<R$, then $\|\xi'(x)-\xi'(y)\|_q<C\epsilon^{p/q}$;
\item if $d(x,y)\geq{}S$, $\|\xi'(x)-\xi'(y)\|_q\geq\frac{p}{q}\delta$.
\end{itemize}
The case $p>q$ is essentially the same, so this gives the result.
\end{proof}

\begin{embcor} 
$l^2$ coarsely embeds in $l^p$ for all $1\leq{}p<\infty$, whence any $X$ that coarsely
embeds in $l^2$ coarsely embeds in $l^p$ for all $1\leq{}p<\infty$.\\
Hence property A implies
coarse embeddability of separable spaces into $l^p$ for all such $p$.
\end{embcor}

Note that it is not too difficult to generalise theorem \ref{yembed} using
\ref{equiv} (3) to get that property A implies coarse embeddability in
$l^p(X)$, $1\leq{}p<\infty$ directly; see \cite{dr}. 
Following the route above gives us additional
results on coarse embeddability, however.  We sum up with:
 
\begin{embthe}[\cite{no2}, theorem 5]
Let $X$ be a separable metric space. Then the following are equivalent:
\begin{itemize}
\item $X$ coarsely embeds in a Hilbert space.
\item $X$ coarsely embeds in $l^p$ for some (hence all) $1\leq{}p\leq{}2$.
\end{itemize}
\end{embthe}

This result is in a sense sharp, as $l^p$ is known not to admit a coarse
embedding
into $l^2$ for $p>2$.  This is due to W. Johnson and N. Randrianarivony, \cite{jr}.  On the
other hand, we are not aware of any bounded geometry
spaces that coarsely embed in $l^p$, $p>2$, but not in $l^2$; indeed, it may be the case that no such
exist.

\pagebreak

\section {Connections with $C^*$-algebra theory} 
\label{analysis}

\subsection{Introduction}

This chapter discusses the relationship of property A to certain $C^*$-algebras that can be 
built out of a group, or a uniformly discrete metric space.  Specifically:
\begin{itemize}
\item A uniformly discrete, bounded geometry space $X$ has property A if and only if the uniform Roe
algebra, $C_u^*(X)$, is nuclear.
\item A countable discrete group $G$ has A if and only if the reduced group 
$C^*$-algebra, $C_r^*(G)$, is exact.
\end{itemize}

We will prove in theorem \ref{nuclexact} that for a countable discrete group $G$,
property A, nuclearity of $C^*_u(G)$ (we abuse notation by using $G$ to refer both to the
group and the underlying coarse space), and exactness of $C^*_r(G)$ are equivalent.
This is a theorem of E. Guentner, J. Kaminker and N. Ozawa, \cite{gk1}, \cite{gk2} and \cite{oz} (the 
results of Higson and Roe, \cite{hr}, and Anantharaman-Delaroche and Renault \cite{adr} provide another approach).
Following \cite{bo}, we then expand
this to get the full result.

The chapter is split into two main parts:

Section \ref{cucr} moves slowly, introducing $C^*_u(X)$ and $C_r^*(G)$, and discussing some of their
properties.  In particular, there is a sense in which the former is generated by positive type kernels, and the
latter by positive type functions (defined in \ref{posfun} below), which we explore.  

Section \ref{nuclexactsect} is somewhat more technical.  We use characterisations of exactness and
nuclearity in terms of finite-dimensional approximation properties to prove the main
results.  The arguments are taken from \cite{bo} (thank you to Nate Brown for providing me with
a preprint of this very readable book).

Throughout, we introduce everything as we need it, and most proofs are included.  
The only $C^*$-algebraic facts we use
without a statement are the spectral calculus, the GNS construction, the Gelfand-Naimark theorem and certain basic facts about 
positive elements.  Nonetheless,
this is probably the most technical chapter of these notes, and some of the arguments
in the section \ref{nuclexactsect}
may be difficult to follow with no background in this area.

Throughout this chapter, $X$ is a uniformly discrete bounded geometry metric space and $G$ is a countable discrete group.
Contrary to our usual convention, we assume that $G$ is equipped with a bounded geometry \emph{right} invariant metric.
Had we used a left-invariant metric on $G$, $C^*_r(G)$ would be isomorphic to a $C^*$-subalgebra of $C^*_u(G)$, but with a right invariant metric it
actually \emph{is} a $C^*$-algebra of $C_u^*(G)$, which simplifies some arguments.

\subsection{$C^*_u(X)$ and $C^*_r(G)$}\label{cucr}

In non-commutative geometry, these two spaces can be regarded as `non-commutative algebras 
of functions' on the coarse space $X$, and on a space of representations of the group $G$ respectively;
each encodes a lot of information about the spaces in question (somewhat analogously to the commutative
$C^*$-algebra $C(Y)$ encoding information about a compact Hausdorff topological space $Y$).
Unfortunately, more detailed comments in this direction are beyond the scope of this piece, but we hope
this provides some motivation.
 
We start with the construction of the uniform Roe algebra, $C_u^*(X)$.

\begin{finprop} \label{finprop}
If $E\subseteq{}X$, we write $P_E\in\mathcal{B}(l^2(X))$ for the orthogonal projection 
onto the (closed) span of $\{\delta_x:x\in{}E\}$.\\
A bounded operator $T$ is of \emph{finite propagation} if there exists $S\geq{}0$ 
such that $P_{\{x\}}TP_{\{y\}}=0$ whenever $d(x,y)>S$.\\
Its \emph{propagation} is the smallest such $S$.
\end{finprop}

Note that $l^2(X)$ is \emph{locally finite dimensional} in the sense that $P_E(l^2(X))$ is
finite dimensional for all bounded $E\subseteq{}X$.

Intuitively, if one thinks of operators on $l^2(X)$ as (infinite) matrices, the finite
propagation operators are those whose entries vanish after a certain distance from the
diagonal.  It is easy to see that these operators are closed under taking linear
combinations, compositions and adjoints, so form a $*$-subalgebra of $\mathcal{B}(l^2(X))$. 

\begin{roealg} \label{roealg}
The \emph{uniform Roe algebra}, or occasionally \emph{uniform translation
algebra}, $C_u^*(X)$, of a metric space $X$ is the closure in 
$\mathcal{B}(l^2(X))$ of the $*$-algebra of finite propagation operators.
\end{roealg}

\begin{roealgrem} \label{roealgrem}
One can also define a (`non-uniform') Roe algebra (or \emph{translation $C^*$-algebra}) of a general metric space $Y$, denoted $C^*(Y)$; it is this object that
appears on the right-hand side of the coarse Baum-Connes conjecture as mentioned in the introduction.
We will not discuss it here, but note that it is \emph{never} exact if $Y$ is unbounded.
See, for example, chapter 6 of \cite{hr2} or \cite{roe1}, section 4.4 for more on this.
\end{roealgrem}

We shall now define the reduced $C^*$-algebra of $G$.

\begin{leftrep}
For any $f:G\to\mathbb{C}$, and $g\in{}G$, we define for all $h\in{}G$,
$(\lambda_gf)(h)=f(g^{-1}h)$.\\
$\lambda_g$ thus defines an element of $\mathcal{B}(l^2(G))$ of norm one.\\
The map $\lambda:G\to\mathcal{B}(l^2(G))$, $g\mapsto\lambda_g$ is called the 
\emph{left regular representation of G}.
\end{leftrep}

\begin{leftreprem}
Note that
$(\lambda_{g_1}\lambda_{g_2}f)(g)=(\lambda_{g_2}f)(g_1^{-1}g)=f(g_2^{-1}g_1^{-1}g)~$
and $(\lambda_g)^*=\lambda_{g^{-1}}=(\lambda_g)^{-1}$,
whence $\lambda$ is a group homomorphism into the unitary operators on 
$l^2(G)$.\\
\end{leftreprem}

\begin{gpalg}
The group algebra of $G$ is the vector space of finitely
supported functions on $G$, with multiplication and adjoint operation defined by
\begin{displaymath}
\left(\sum_{g\in{}G}a_g\delta_g\right)\left(\sum_{g\in{}G}b_g\delta_g\right)
=\sum_{g,h\in{}G}a_gb_h\delta_{gh}
\end{displaymath} and \begin{displaymath}
\left(\sum_{g\in{}G}a_g\delta_g\right)^*=\sum_{g\in{}G}\overline{a_g}\delta_{g^{-1}}.
\end{displaymath}
\end{gpalg}

It follows that $\lambda$ extends by linearity to an injective $*$-homomorphism
$\mathbb{C}[G]\to\mathcal{B}(l^2(G))$ (which we also denote $\lambda$) by
\begin{displaymath}
\lambda:\sum_{g\in{}G}a_g\delta_g\mapsto\sum_{g\in{}G}a_g\lambda_g.
\end{displaymath}

\begin{redgp}
The \emph{reduced $C^*$-algebra of $G$}, $C^*_r(G)$, is the norm closure in 
$\mathcal{B}(l^2(G))$ of $\lambda(\mathbb{C}[G])$.
\end{redgp}

We will also write $\lambda$ for the inclusion map $C^*_r(G)\to\mathcal{B}(l^2(G))$.

Note that one can by analogy define the \emph{right regular representation} of $G$ on
$l^2(G)$ by $(\rho_gf)(h)=f(hg)$.  The norm completion of span$\{\rho_g:g\in{}G\}$ 
in $\mathcal{B}(l^2(G))$ defines a $C^*$-algebra $C^*_{\rho}$, which is isomorphic to 
$C^*_r(G)$ via the map $\rho_g\mapsto\lambda_g$.  It is a $C^*$-subalgebra
of $C^*_u(G)$ in the case $G$ is equipped with a \emph{left} invariant metric.

We now define positive type functions, an analogue of kernels of positive type. 
They are introduced here as 
positive type kernels have a close relationship to $C^*_u(X)$, and positive type functions
have much the same relationship to $C^*_r(G)$.  Positive type functions also have essentially the same relationship
to amenability (see \ref{amendef}, part (3)) as positive type kernels do to property A (see \ref{equiv}, part (8)).

\begin{posfun} \label{posfun}
A \emph{positive type function} on $G$ is a map 
$\phi:G\to\mathbb{C}$ such
that for all $g_1,...,g_n\in{}G$ and all $\lambda_1,...,\lambda_n\in\mathbb{C}$, 
\begin{displaymath}
\sum_{i,j=1}^n\lambda_i\overline{\lambda_j}\phi(g_i^{-1}g_j)\geq0
\end{displaymath}
Such a $\phi$ is \emph{normalised} if $\phi(e)=1$.\\
Negative type functions can be defined analogously.
\end{posfun}

A positive type function is the same thing as an invariant positive type kernel with respect to the left
(respectively, right) action of $G$ on $G\times{}G$ via the identification $k(g,h)=\phi(g^{-1}h)$ 
(resp. $k(g,h)=\phi(gh^{-1})$.

We record the following facts that connect certain structural features of
the $C^*$-algebras defined above and the notions of positive type function and kernel.
They are included both because (we hope!) they aid intuition of how these two $C^*$-algebras and 
property A (and amenability) are related.

\begin{cstarfacts} \label{c*facts1}
A (not necessarily bounded) operator $T$ on $l^2(X)$ defines a kernel $k_T$ on $X$ by
\begin{displaymath}
k_T(x,y)=\langle{}T\delta_y,\delta_x\rangle
\end{displaymath}
A finite propagation kernel $k$ on $X$ defines a 
\emph{convolution operator} $T_k$ on $l^2(X)$ by
\begin{displaymath}
(T_kf)(x)=\sum_{y\in{}X}k(x,y)f(y).
\end{displaymath}
Moreover, $T_k$ is of finite propagation if and only if $k_T$ is, and in this case the two processes are mutually inverse.\\
If $k_T$, $T_k$ are finite propagation, then $T_k$ is bounded if and only if $k_T$ is.\\
If $k_T$, $T_k$ are bounded and finite propagation, then $T_k$ is positive if and only if 
$k_T$ is of positive type.\\
\end{cstarfacts}

\begin{proof}
The proofs are easy calculations.  For example, to prove the positivity result,
note that if $T$ is positive then for any $x_1,...,x_n\in{}X$ and any 
$\lambda_1,...,\lambda_n\in\mathbb{C}$,
the function $f=\sum_{i=1}^n\lambda_i\delta_{x_i}$ is in $l^2(X)$ and 
\begin{displaymath}
\sum_{i,j}\lambda_i\overline{\lambda_j}k_T(x_i,x_j)=\langle{}Tf,f\rangle\geq0.
\end{displaymath}
The converse follows as if $k$ is of positive type, then $\langle{}T_kf,f\rangle\geq{}0$ for
all finitely supported $f\in{}l^2(X)$ by the same identity.  Hence $\langle{}T_kf,f\rangle\geq{}0$ for all $f\in{}l^2(X)$
as the finitely supported functions are dense.
\end{proof}

In order to state some corresponding facts about positive definite functions we need 
the following:

\begin{cdd}
An operator $T$ on $l^2(G)$ is said to be \emph{constant down the diagonals} if
\begin{displaymath}
\langle{}T\delta_{g_1},\delta_{g_2}\rangle=\langle{}T\delta_{h_1},\delta_{h_2}\rangle
~\textrm{whenever}~g_1g_2^{-1}=h_1h_2^{-1}
\end{displaymath}
(draw the matrix of an operator on $l^2(G)$, $G$ cyclic, to see where the name comes 
from).\\
Equivalently, these are the $G$-equivariant operators on $l^2(G)$ with respect to the
right regular representation, i.e.
those for which $T(\rho_gf)=\rho_gT(f)$ for all $f\in{}l^2(G)$ and all $g\in{}G$ (another
alternative description is that they constitute the commutant of $C^*_{\rho}(G)$ in $\mathcal{B}(l^2(G))$).
\end{cdd}

\begin{cstarfacts2} \label{c*facts2}
A (not necessarily bounded) operator $T$ on $l^2(G)$ defines a function $\phi_T$ on $G$ by
\begin{displaymath}
\phi_T(g)=\langle{}T\delta_e,\delta_g\rangle
\end{displaymath}
A finitely supported function $\phi$ on $G$ defines a \emph{convolution operator} 
$T_{\phi}$ on $l^2(G)$ by
\begin{displaymath}
(T_{\phi}f)(g)=\sum_{h\in{}G}\phi(gh^{-1})f(h)
\end{displaymath}
If $T$ is constant down the diagonals, these processes are mutually inverse.\\
Assuming this, we get that:\\
$T_{\phi}$ is of finite propagation if and only if $\phi_T$ is finitely supported;
either implies that $T_{\phi}$ is bounded.\\
If $\phi_T$ is finitely supported and $T_\phi$ is finite propagation, then $\phi_T$ is of positive type if and only if 
$T_{\phi}$ is positive.
\end{cstarfacts2}

The proofs are simple calculations again.  For example, to show that the processes described really are mutually
inverse when $T$ is constant down the diagonals, note that:
\begin{align*}
\phi_{T_{\phi}}(g)=\langle{}T_{\phi}\delta_e,\delta_g\rangle
=(T_{\phi}\delta_e)(g)=\sum_{h\in{}G}\phi(gh^{-1})\delta_e(h)=\phi(g),
\end{align*}
and
\begin{align*}
\langle{}T_{\phi_T}\delta_g,\delta_h\rangle & =(T_{\phi_T}\delta_g)(h)=
\sum_{k\in{}G}\phi_T(hk^{-1})\delta_g(k)=\phi_T(hg^{-1})
=\langle{}T\delta_e,\delta_{hg^{-1}}\rangle \\
& =\langle{}T\delta_g,\delta_h\rangle,
\end{align*}
the last step being by the constant-down-diagonals property.

The final result in this section is the main motivation for the preceding lemmas; we hope it ties some of this
information together.

\begin{redroe} \label{redroe} 
\begin{itemize}
\item $C^*_r(G)$ is the $C^*$-subalgebra of $\mathcal{B}(l^2(G))$ generated by the elements $T_{\phi}$ such that
$\phi$ is a finitely supported positive type function on $G$.
\item $C^*_u(X)$ is the $C^*$-subalgebra of $\mathcal{B}(l^2(X))$ generated by the elements $T_k$ such that
$k$ is a finite propagation, bounded, positive type kernel on $X$.
\item In the case $X$ is the underlying metric space of $G$, $C^*_r(G)$ is a $C^*$-subalgebra of
$C^*_u(X)$.
\end{itemize}
\end{redroe}

\begin{proof}
Write $A$ for the $C^*$-subalgebra of $\mathcal{B}(l^2(G))$ generated by the elements $T_\phi$, where $\phi$ is finitely supported and positive type.\\
To show that $C_r^*(G)$ is contained in $A$, it suffices to show that any positive element $T\in{}C^*_r(G)$ is contained in $A$ 
(as any element in a $C^*$-algebra can be written as a linear combination of four positive elements).  Now, $T$ has a positive square root $\sqrt{T}$, and by definition of $C^*_r(G)$, we may approximate $\sqrt{T}$ by some sequence $(S_n)$ of elements of $\lambda{}(\mathbb{C}[G])$.    Now, $S_n^*S_n$ is positive and in $\lambda(\mathbb{C}[G])$ for each $n$, whence equal to $T_{\phi_n}$ for some finitely supported, positive type $(\phi_n)$.  As $S_n\to\sqrt{T}$ as $n\to\infty$, $T_{\phi_n}\to{}T$, and we are done.\\
To get the converse inclusion, note that if $\phi$ is finitely supported and of positive type
\begin{align*}
\langle{}T_{\phi}\delta_g,\delta_h\rangle=\phi(hg^{-1})
=\left\langle\sum_{k\in\textrm{support}(\phi)}\phi(k)\lambda_k\delta_g,\delta_h\right\rangle,
\end{align*}
whence $T_{\phi}=\sum\phi(k)\lambda_{k}$ is an element of $C^*_r(G)$, and so $A$ is contained in $C^*_r(G)$.

The second identification is proved in a similar manner.

For the final part, note that if $\phi$ is finitely supported and of positive type, then setting $k(g,h)=\phi(gh^{-1})$ defines
a finite propagation, positive type, bounded kernel and $T_k=T_{\phi}$.  The only possibly difficult thing to check is finite propagation of $k$;
note, however, that support$(\phi)\subseteq{}B(e,R)$ for some $R$, whence $k(g,h)\not=0$ forces
$d(e,gh^{-1})<R$, whence $d(g,h)<R$ by right invariance of $d$.

\end{proof}

We could have given a direct proof that $C^*_r(G)\subseteq{}C^*_u(G)$, which would have been quicker, but 
chose instead to discuss kernels, hoping this would make clearer the connection of this chapter to
chapter \ref{kerembed}.

This completes the preliminaries (which are accessible even if one knows very
little about $C^*$-algebras; for an interesting and useful application see $\cite{gk4}$).
The next section introduces some rather more technical notions,
and uses them to prove the main results.

\subsection{Nuclearity and exactness}\label{nuclexactsect}

\begin{c*algdefns} \label{c*algdefns}
A linear map $\phi:A\to{}B$ of unital $C^*$-algebras is 
\emph{completely positive} (or cp) if its extension (by applying $\phi$ entry-wise)
$\phi_n:M_n(A)\to{}M_n(B)$ takes positive elements to positive elements for all
$n\in\mathbb{N}$.\\
A \emph{unital completely positive} (or ucp) map is a cp map that is also unital.\\
A linear map $\theta:A\to{}B$ is \emph{nuclear} if for any finite subset $F\subseteq{}A$ and
any $\epsilon>0$ there exists $n$ and ucp maps $\phi:A\to{}M_n(\mathbb{C})$,
$\psi:M_n(\mathbb{C})\to{}B$ such that $\|\theta(a)-(\psi\circ\phi)(a)\|<\epsilon$ for
all $a\in{}F$.
\end{c*algdefns}

Note that a $*$-homomorphism $\phi$ is automatically cp (ucp if unital), as the extensions $\phi_n$ are also
$*$-homomorphisms, whence preserve positivity.

Note that one can define nuclear maps in the non-unital case, but we will not need this extra complication, so omit it.

To establish the main result, we need a fairly serious collection of preliminaries. 
We start with the following basic example of a cp (or ucp) map.

\begin{ucplem} \label{ucplem}
Let $\pi:A\to\mathcal{B}(\mathcal{H})$ be a $*$-homomorphism of 
$C^*$-algebras, and $T:\mathcal{H}'\to\mathcal{H}$ a bounded operator.\\
Then the map $\phi:a\mapsto{}T^*\pi(a)T$ is cp.  If $\pi$ is
unital and $T$ is an isometry, $\phi$ is ucp.
\end{ucplem}

\begin{proof}
Let $a=[a_{ij}]_{i,j=1}^n$ be a positive element in $M_n(A)$.  As it is a positive element of a 
$C^*$-algebra, we can write $a=b^*b$ for some $b\in{}M_n(A)$.\\
Note that $\phi_n(a)=T_n^*\pi(b)^*\pi(b)T_n$, where
$T_n$ is the diagonal matrix with all (diagonal) entries $T$ in $M_n(\mathcal{B}(\mathcal{H}))$.\\
Hence for any $x\in\mathcal{H}^n$, we get that
\begin{displaymath}
\langle{}\phi_n(a)x,x\rangle=
\langle\pi(b)T_nx,\pi(b)T_nx)\rangle\geq{}0,
\end{displaymath}
where the inner product is that of $\mathcal{H}^n$. Hence $\phi_n$ is cp.\\
The second part follows as if $\pi$ is unital and $T$ an isometry, $\phi(1)=T^*\pi(1)T=T^*T=1$, whence
$\phi_n$ is also unital for all $n$.
\end{proof}

It is a consequence of the next result, called Stinespring's dilation theorem, that every cp (or ucp) map
arises in this way (even for non-unital algebras).  For intuition (and as we need the corollary that follows it), we give
a proof following \cite{bo} 1.5.3.

\begin{stinespring} \label{stinespring}
Let $A$, $B$ be unital $C^*$-algebras, and let $\phi:A\to{}B$ be a cp map.  By faithfully 
representing $B$, assume $B\subseteq\mathcal{B}(\mathcal{H})$ for some Hilbert space $\mathcal{H}$.\\
Then there exists a triple $(\mathcal{H}',\pi,T)$ such that:
\begin{itemize}
\item $\mathcal{H}'$ is a Hilbert space;
\item $T:\mathcal{H}\to\mathcal{H}'$ is a bounded linear operator;
\item $\pi:A\to\mathcal{B}(\mathcal{H}')$ is a representation (i.e. a $*$-homomorphism);
\end{itemize} 
and $\phi(a)=T^*\pi(a)T$ for all $a\in{}A$.
\end{stinespring}

\begin{proof}
Let $A\otimes{}\mathcal{H}$ be the (incomplete) algebraic tensor product, and define a sesquilinear form
on it by setting
\begin{displaymath}
\left\langle\sum_ia_i\otimes{}v_i,\sum_jb_j\otimes{}w_j\right\rangle
=\sum_{i,j}\langle\phi(b_j^*a_i)v_i,w_j\rangle_{\mathcal{H}},
\end{displaymath}
which a similar calculation to that used in the previous lemma shows to be positive semi-definite.  By the argument from
\ref{kerthe}, then, it satisfies a Cauchy-Schwarz type inequality, and the elements $x$ for which
$\langle{}x,x\rangle=0$ form a linear subspace.  Quotient out by it (say $Q$ is the quotient map), and complete in the norm corresponding to $\langle,\rangle$ to get a Hilbert space $\mathcal{H}'$.\\
Now, let $J:\mathcal{H}\to\mathcal{H'}$ be the inclusion of $\mathcal{H}$ in
$A\otimes\mathcal{H}$ by $v\mapsto1_A\otimes{}v$, and define $T=Q\circ{}J:\mathcal{H}\to\mathcal{H}'$, which
is a bounded linear map.\\
Define $\pi:A\to\mathcal{B}(\mathcal{H}')$ by setting 
\begin{displaymath}
\pi(a)Q(\sum_ia_i\otimes{}v_i)=Q(\sum_iaa_i\otimes{}v_i)
\end{displaymath}
on the dense subspace $Q(A\otimes\mathcal{H})$ of $\mathcal{H}'$ and extending to the completion.  It is clearly linear and multiplicative, while $*$-preserving follows as 
\begin{align*}
& \left\langle\pi(a)Q(\sum_ia_i\otimes{}v_i),Q(\sum_jb_j\otimes{}w_j)\right\rangle
=\sum_{i,j}\langle\phi(b_j^*aa_i)v_i,w_j\rangle_{\mathcal{H}} \\
& =\sum_{i,j}\langle\phi((a^*b_j)^*a_i)v_i,w_j\rangle_{\mathcal{H}}
=\left\langle{}Q(\sum_ia_i\otimes{}v_i),\pi(a^*)Q(\sum_jb_j\otimes{}w_j)\right\rangle,
\end{align*}
whence $\pi$ is a $*$-homomorphism.\\
Finally, note that
\begin{displaymath}
\langle{}Q(a\otimes{}w),Tv\rangle=\langle\phi(a)w,v\rangle_{\mathcal{H}'}
=\langle{}T^*(Q(a\otimes{}w)),v\rangle_{\mathcal{H}'},
\end{displaymath}
whence $T^*(Q(a\otimes{}w))=\phi(a)w$. Hence
\begin{align*}
(T^*\pi(a)T)v=T^*Q(a\otimes{}v)=\phi(a)v
\end{align*} 
for all $v\in\mathcal{H}$, which completes the result.
\end{proof}

Note that a state $f$ on $A$ is a ucp map $A\to\mathbb{C}$ (exercise!).  The above proof is a direct generalisation
of the GNS construction, which builds a representation of $A$ out of $f$. 

We will use the following easy:

\begin{stinec}
If $\phi$ is cp, $\phi$ is $*$-preserving.\\
If $\phi$ is ucp, $\phi$ is contractive.
\end{stinec}
\begin{proof}
The first part is immediate, as $(T^*\pi(a^*)T)^*=T^*\pi(a)T$.\\
The second part follows as if $\phi$ is unital then so is $\pi$ by inspection of the above
proof.  It follows that $1=\phi(1)=T^*T$, whence $T$ is an isometry, and so in particular
$\|\phi\|\leq\|T^*\|\|\pi\|\|T\|=1$, so $\phi$ contractive.
\end{proof}

The  next lemma gives a useful characterisation of
completely positive maps from a $C^*$-algebra to a matrix algebra.

\begin{ucplem2} \label{ucplem2}
Let $A$ be a $C^*$-algebra.  There is a one-to-one correspondence
\begin{displaymath}
\{cp ~\textrm{maps} ~\phi:A\to{}M_n(\mathbb{C})\}\leftrightarrow
\{\textrm{positive linear functionals} ~f\in{}M_n(A)^*\}
\end{displaymath}
given by
\begin{displaymath}
f_{\phi}([a_{ij}])=\sum_{i,j}\phi(a_{ij})_{i,j}, ~\phi_f(a)=[f(ae_{ij})]_{i,j=1}^n
\end{displaymath}
where $\phi(a_{ij})_{i,j}$ is the $(i,j)^{th}$ matrix coefficient of 
$\phi(a_{ij})\in{}M_n(\mathbb{C})$ and $e_{ij}$ are the matrix units of $M_n(A)$.
\end{ucplem2}

\begin{proof}
Note first that if $e_1,...,e_n$ are the elements of the standard basis for 
$\mathbb{C}^n$, and $e$ is their
sum, we get that $f_{\phi}([a_{ij}])=\langle\phi([a_{ij}])e,e\rangle$, so $\phi$ cp
implies $f_{\phi}$ positive.\\
Conversely, let $f$ be a positive linear functional on $M_n(A)$ with associated
GNS triple $(\pi,\mathcal{H},v)$.  Assume for the rest of the proof that inner products
are linear in the second variable.\\
Define a linear map $T:\mathbb{C}^n\to\mathcal{H}$ by setting $Te_j=\pi(e_{1j})v$, and note
that $T^*\delta_{ae_{1k}}=\sum_{i=1}^nf(ae_{ik})e_i$ for any $a\in{}A$.\\
Hence for any $a\in{}A$, any $k=1,...,n$, 
$\phi_f(a)(e_k)=[f(ae_{ij}](e_k)=\sum_{i=1}^nf(ae_{ik})e_i$.\\
On the other hand,
\begin{displaymath}
T^*\pi\left(
\begin{bmatrix} a & & \\ & \ddots & \\ & & a \end{bmatrix}
\right)Te_k=T^*(\delta_{ae_{1k}})=\sum_{i=1}^nf(ae_{ik})e_i.
\end{displaymath}
Hence
\begin{displaymath}
\phi_f(a)=T^*\pi\left(
\begin{bmatrix} a & & \\ & \ddots & \\ & & a \end{bmatrix}
\right)T,
\end{displaymath}
which is ucp by the previous lemma.
\end{proof}   

\begin{ucplem3} \label{ucplem3}
Let $A\subseteq{}B$ be $C^*$-algebras and $\phi:A\to{}M_n(\mathbb{C})$ be a ucp map.\\
Then $\phi$ extends to a ucp map defined on all of $B$.
\end{ucplem3}

\begin{proof}
By the previous lemma, $\phi$ defines a positive linear functional $f_{\phi}$ on 
$M_n(A)$; as $\phi$ is unital, $f_{\phi}(1)=1$.\\
By Hahn-Banach, $f_{\phi}$ extends to $f_{\phi'}$ defined on all of $M_n(B)$, which moreover has norm one as 
$f_{\phi}$ does (this follows as $f_{\phi}$ is positive and $f_{\phi}(1)=1$).\\
Note that if $b\in{}M_n(B)$ is positive then
\begin{displaymath}
\|b\|-f_{\phi'}(b)=|\|b\|-f_{\phi'}(b)|=|f_{\phi'}(\|b\|-b)|\leq\|\|b\|-b\|\leq\|b\|
\end{displaymath}
whence $f(b)\geq{}0$ and by the previous lemma, there exists 
ucp $\phi':B\to{}M_n(\mathbb{C})$ extending $\phi$.
\end{proof} 

We next define exactness and nuclearity of $C^*$-algebras;
our definition is non-standard, but equivalent to the `historical' one by a deep
theorem of Choi-Effros and Kirchberg; see for example \cite{bo} or \cite{wass} for a proof 
and discussion.  The original definition is in terms of tensor products:
specifically a $C^*$-algebra $A$ is \emph{nuclear} if for any 
$C^*$-algebra $B$ there is a unique $C^*$-algebra norm on the algebraic tensor
product $B\otimes{}A$ satisfying $\|a\otimes{}B\|=\|a\|\|b\|$.  $A$ is \emph{exact} if the functor $\otimes_{\min}A$ is exact.  We
will not need this, and only mention it for completeness.

\begin{exact}
A unital $C^*$-algebra $A$ is \emph{nuclear} if the identity map $A\to{}A$ is a
nuclear map.\\
A unital $C^*$-algebra $A$ is \emph{exact} if it is nuclearly representable, i.e. if there 
exists
a Hilbert space $\mathcal{H}$ and a faithful representation 
$\pi:A\to\mathcal{B}(\mathcal{H})$ that is also a nuclear map.\\
A (countable discrete) group $G$ is called \emph{exact} if $C^*_r(G)$ is.
\end{exact}

Note that in the case of a group, exactness is equivalent to the
nuclearity of the map $\lambda:C^*_r(G)\to\mathcal{B}(l^2(G))$; one direction is
immediate, while the other follows as a composition of 
a $*$-homomorphism and a nuclear map is nuclear. 
Note also that to show a ucp map nuclear, it will suffice to show that it approximately 
factors through
a nuclear $C^*$-algebra in the sense of the second 
part of \ref{c*algdefns} (i.e. as \ref{c*algdefns}, but
with $M_n(\mathbb{C})$ replaced by some nuclear $C^*$ algebra).   
We will use these facts, and the following 
fundamental class of nuclear $C^*$-algebras.

\begin{nuclabel} \label{nuclabel}
(Unital) abelian $C^*$-algebras are nuclear.
\end{nuclabel}

`Unital' is in parentheses as it is possible to define nuclearity for non-unital $C^*$-algebras, and the
result remains true for these.  We will not use this, however, and thus restricted to a (simpler) definition that
only makes sense in the unital case.

\begin{proof} (\cite{bo} 2.4.2)
Consider any unital abelian $C^*$-algebra, which is of the form $C(X)$, $X$ some compact Hausdorff topological space, by the Gelfand-Naimark theorem.\\
Given a finite set $F$ of functions on $X$ and $\epsilon>0$, we can find
a (finite) covering $\{U_i\}_{i=1}^n$ of $X$, elements $x_i\in{}U_i$ and a subordinate
partition of unity $\{\xi_i\}$ such that each $f\in{}F$ is approximated within
$\epsilon$ by $\sum_{i=1}^nf(x_i)\xi_i$ (in the $C(X)$ norm).\\
Define $\phi:C(X)\to{}M_n(\mathbb{C})$ by sending any $f$ to the diagonal matrix with entries
$f(x_1),...,f(x_n)$.  This is a unital $*$-homomorphism, so certainly ucp.  Define moreover 
$\psi:M_n(\mathbb{C})\to{}C(X)$ as the composition of
\begin{displaymath}
\psi_1:A\mapsto{}\sum_{i=1}^ne_{i,i}Ae_{i,i}~~\textrm{and}~~\psi_2:(a_1,...,a_n)\mapsto\sum_{i=1}^na_i\xi_i.
\end{displaymath}
Here $e_{i,j}$ are the matrix units of $M_n(\mathbb{C})$, so $\psi_1:M_n(\mathbb{C})\to\mathbb{C}^n$ and
$\psi_2:\mathbb{C}^n\to{}C(X)$.  Note that $\psi_1$ is ucp (e.g. by \ref{ucplem}), and that $\psi_2$ is ucp, as 
a unital $*$-homomorphism.  As
\begin{displaymath}
(\psi\circ\phi)f=\sum_{i=1}^nf(x_i)\xi_i
\end{displaymath}
for all $f\in{}C(X)$, choice of the partition of unity completes the proof. 
\end{proof}    

Note that finite sums of, and matrix algebras over, nuclear $C^*$ algebras are also 
nuclear.  We are now ready to
prove the main theorem of this chapter.

\begin{nuclexact} \label{nuclexact}
Let $G$ be a countable discrete group.\\
The following are equivalent:
\begin{enumerate}
\item $C_r^*(G)$ is exact.
\item $G$ has property A.
\item $C^*_u(G)$ is nuclear.
\item $C^*_u(G)$ is exact.
\end{enumerate}
\end{nuclexact}

\begin{proof}

The proof we give follows \cite{bo}, 5.1.5 and 5.5.7.\\

(1) implies (2):\\
Let $R,\epsilon>0$.\\
By exactness there exist $n$ and ucp $\phi:C^*_r(G)\to{}M_n(\mathbb{C})$, 
$\psi:M_n(\mathbb{C})\to\mathcal{B}(l^2(G))$ such that 
\begin{displaymath}
\|(\psi\circ\phi)\lambda_g-\lambda_g\|<\epsilon/2 ~\textrm{for all} ~g\in\bar{B}(e,R)
\end{displaymath}

We start the argument by showing that $\phi$ can be approximated by some $\phi'$ with
finite dimensional support.\\
Let $f_{\phi}\in{}M_n(C^*_r(G))^*$ be as in \ref{ucplem2}, and use \ref{ucplem3}
to extend it to a positive linear functional on $\mathcal{B}(l^2(G))$.\\
Say that an element of $\mathcal{B}(l^2(G))^*$ is of \emph{finite character} if it is
a finite linear combination of elements of the form $T\mapsto\langle{}Tf,f'\rangle$, where
$f,f'\in{}l^2(G)$ are of finite support.\\
Such elements are weak-$*$ dense in $\mathcal{B}(l^2(G))^*$ (obvious!), whence if
$e_{ij}$ are the matrix units in $M_n(C^*_r(G))$, we can find positive $f_{\phi'}$ 
that is arbitrarily close to $f_{\phi}$ on the (finitely many) elements $\lambda_ge_{ij}$
such that $|g|\leq{}R$.\\
By the one-to-one correspondence (and explicit formula 
$\phi'(a)=[f_{\phi'}(ae_{ij})]_{i,j=1}^n$) in \ref{ucplem2}, therefore, it
follows that $\|\phi'(\lambda_g)-\phi(\lambda_g)\|$ is small for all such $g$.\\
Note that $P_{B(e,S)}\phi'P_{B(e,S)}=\phi'$ for some $S>0$, as $\phi'$ corresponds
to a linear functional of finite character. 
Note moreover that as $\phi'(\lambda_e)$ is close to 
$\phi(\lambda_e)=1\in{}M_n(\mathbb{C})$, it is invertible (as everything within
distance one of the identity is); letting $p=\phi'(\lambda_e)$,
use the spectral calculus to define
\begin{displaymath}
\phi'':C^*_r(G)\to{}M_n(\mathbb{C}) ~\textrm{by setting} 
~\phi''(T)=p^{-\frac{1}{2}}\phi'(T)p^{-\frac{1}{2}},
\end{displaymath}
which is also ucp, as $p^{-\frac{1}{2}}$ is positive.\\
Finally, by approximating $\phi$ arbitrarily well in this way, we get a map 
$\theta=\psi\circ\phi''$ such that
\begin{displaymath}
\|\theta(\lambda_g)-\lambda_g\|<\epsilon ~\textrm{as long as} ~|g|\leq{}R
\end{displaymath}
and $P_{B(e,S)}\theta{}P_{B(e,S)}=\theta$ for some $S>0$ (as the same is true for $\phi''$; this also uses that $\psi$ is contractive, as ucp). 

Now, define a kernel $k$ on $G$ by setting
\begin{displaymath}
k(g,h)=\langle\theta(\lambda_{gh^{-1}})\delta_h,\delta_g\rangle
\end{displaymath}
Note that $k(g,g)=1$ as $\theta$ is unital, that $k$ is self-adjoint (as $\theta$ is $*$-preserving), and if
$d(g,h)>2S$ then $\theta(\lambda_{gh^{-1}})=0$ by right invariance of $d$,
so $k$ has propagation at most $2S$.\\
Note moreover that $k$ is of positive type.  For if $g_1,...,g_n$ are elements of
$G$ and $\alpha_1,...,\alpha_n$ are complex numbers, we get that
\begin{align*}
\sum_{i,j}\alpha_i\overline{\alpha_j}k(g_i,g_j)
& =\sum_{i,j}\langle
\alpha_i\theta(\lambda_{g_ig_j^{-1}})\delta_{g_i},\alpha_j\delta_{g_j}\rangle \\
& =\sum_j\left\langle\sum_i\theta(\lambda_{g_ig_j^{-1}})\alpha_i\delta_{g_i},\alpha_j\delta_{g_j}\right\rangle \\
& =\langle\theta_n[(\lambda_{g_ig_j^{-1}})]_{i,j=1}^nv,v\rangle\geq0
\end{align*}
where $v=[\alpha_1\delta_{g_1},\hdots,\alpha_n\delta_{g_n}]^T$.  Positivity follows as $\theta$ is ucp and
\begin{displaymath}
\langle[(\lambda_{g_ig_j^{-1}})]_{i,j=1}^nw,w\rangle=
\left\langle\sum_j\lambda_{g_j^{-1}}w_j,\sum_j\lambda_{g_j^{-1}}w_j\right\rangle\geq0
\end{displaymath}
for any $w=[w_1,...,w_n]^T\in{}l^2(G)^n$.
Hence $k$ is of positive type.\\
Finally note that if $d(g,h)\leq{}R$, we get that
\begin{align*}
|1-\textrm{Re}k(g,h)| & =\frac{1}{2}|2-2\textrm{Re}\langle
\theta(\lambda_{gh^{-1}})\delta_h,\delta_g\rangle| \\
& \leq\frac{1}{2}\|\theta(\lambda_{gh^{-1}})\delta_h-\delta_g\|^2<\epsilon
\end{align*}
as long as $\epsilon<1$, as $\theta(\lambda_g)$ is within $\epsilon$ of $\lambda_g$ for $|g|\leq{}R$.\\
Hence Re$k$ has the properties in \ref{equiv}, (8), whence $G$ has property A.\\

(2) implies (3):\\
Let $F\subseteq{}C^*_u(G)$ be finite and $\epsilon>0$.  We want to show that
there exists a nuclear $C^*$-algebra $B$ and maps $\phi:C^*_u(G)\to{}B$ and 
$\psi:B\to{}C^*_u(G)$ such that $\|(\psi\circ\phi)(T)-T\|<\epsilon$ for all $T\in{}F$.

As the finite propagation operators are dense in $C^*_u(G)$, we may as well assume that
$F$ consists only of these.  Let
$R=\max\{\textrm{propagation}(T):T\in{}F\}$.\\
As $X$ has property A, it follows from \ref{equiv}, (3), that we can find
$\eta:G\to{}l^2(G)$ such that $\|\eta_g\|=1$, $\eta_g$ is supported in $B(g,S)$ 
for some $S>0$, all $g$, and $\eta$ has 
$(R,2\sqrt(\max\{\|T\|:T\in{}F\}.\sup\{|B(g,R)|:g\in{}G\}^{-1}.\epsilon))$ variation.\\
Note that without loss of generality, we can assume each $\eta_x$ positive-valued; in
particular then, $\langle\eta_g,\eta_h\rangle\in\mathbb{R}$ for all $g,h\in{}G$.\\
Define
\begin{align*} 
\phi: & C^*_u(G)\to\prod_{g\in{}G}\mathcal{B}(l^2(B(g,S))) \\
& T\mapsto(P_{B(g,S)}TP_{B(g,S)})_{g\in{}G}.
\end{align*}
$\phi$ is ucp by lemma \ref{ucplem}.\\
Moreover, $G$ is bounded geometry, whence $N=\max\{|B(g,S)|:g\in{}G\}$ exists.\\
Letting $Y_k=\{g\in{}G:|B(g,S)|=k\}$ for $k=1,...,N$, we get that
\begin{displaymath}
\prod_{g\in{}G}\mathcal{B}(l^2(B(g,S)))\cong
\bigoplus_{k=1}^Nl^{\infty}(Y_k)\otimes{}M_k(\mathbb{C})
\cong\bigoplus_{k=1}^N{}M_k(l^{\infty}(Y_k)),
\end{displaymath}
which is a nuclear $C^*$-algebra, as a finite direct sum of matrix algebras
over abelian $C^*$-algebras (there is some redundancy here: $|B(g,S)|$ has the same cardinality for all $g\in{}G$, so most of the $Y_k$ are empty; we write it like this so it generalises directly to the case where the space is not a group).

In order to be able to construct a $\psi$ with the right properties, we define
bounded linear operators $V_g:l^2(G)\to{}l^2(B(g,S))$ by setting
$V_g\delta_h=\eta_h(g)\delta_h$.  Note that 
$V_g^*\delta_h=\overline{\eta_h(g)}\delta_h$.\\
We now set
\begin{align*}
\psi: & \prod_{g\in{}G}\mathcal{B}(l^2(B(g,S)))\to\mathcal{B}(l^2(G)) \\
& (T_g)_{g\in{}G}\mapsto\sum_{g\in{}G}V_g^*T_gV_g,
\end{align*}
which defines a ucp map; it is completely positive by \ref{ucplem}, and unital
as for any $h\in{}G$
\begin{displaymath}
\sum_{g\in{}G}V_g^*V_g\delta_h=
\sum_{g\in{}G}\eta_h(g)\overline{\eta_h(g)}\delta_h=\|\eta_h\|\delta_h=\delta_h,
\end{displaymath}  
i.e. $\sum{}V^*_gV_g$ converges and equals one in the strong operator topology.\\
Note moreover that any operator $T_g\in{}B(g,S)$, $V_g^*T_gV_g$ has propagation at 
most $S$, whence $\psi((T_g)_{g\in{}G})$ has propagation at most $S$ for all
$(T_g)_{g\in{}G}$.  In particular, the image of $\psi$ is contained in
$C_u^*(G)$.

Now, let $k(g,h)=\langle\eta_g,\eta_h\rangle$, which defines a positive type
normalised kernel of propagation at most $S$, and of 
$(R,(\max\{\|T\|:T\in{}F\}.\sup\{|B(g,R)|:g\in{}G\})^{-1}.\epsilon)$ variation, as
$|1-k(g,h)|=\frac{1}{2}\|\eta_g-\eta_h\|^2$, recalling that
$\langle\eta_g,\eta_h\rangle\in\mathbb{R}$ for all $g,h\in{}G$.\\
Define a map (the \emph {Schur multiplier} associated to $k$)
$M_k:C^*_u(G)\to{}C^*_u(G)$ by setting $M_k(T)$ to be the operator with matrix
coefficients
\begin{displaymath}
\langle{}M_k(T)\delta_g,\delta_h\rangle=k(g,h)\langle{}T\delta_g,\delta_h\rangle.
\end{displaymath}
It is well defined and takes $C^*_u(G)$ into itself by finite propagation (and
boundedness) of $k$.  We then have that
\begin{align*}
\langle(\psi\circ\phi)(T)\delta_g,\delta_h\rangle
& =\sum_{g'\in{}G}\langle{}V_{g'}TV_{g'}\delta_g,\delta_h\rangle
=\sum_{g'\in{}G}\langle{}TV_{g'}\delta_g,V_{g'}\delta_h\rangle \\
& =\langle{}T\delta_g,\delta_h\rangle\sum_{g'\in{}G}\eta_g(g')\overline{\eta_h(g')}
=\langle\eta_g,\eta_h\rangle.\langle{}T\delta_g,\delta_h\rangle \\
& =k(g,h)\langle{}T\delta_g,\delta_h\rangle=\langle{}M_k(T)\delta_g,\delta_h\rangle.
\end{align*}
Hence $\psi\circ\phi=M_k$; moreover
\begin{displaymath}
\|M_k(T)-T\|\leq\|T\|\sup_{d(g,h)\leq{}R}|1-k(g,h)|\sup_g|B(g,R)|<\epsilon
\end{displaymath}
by comparing matrix coefficients, and choice of $k$.\\

(3) implies (4):\\
Nuclear $C^*$-algebras are exact.\\

(4) implies (1):\\
Note that if $A\subseteq{}B$ are $C^*$-algebras, and $\pi:B\to\mathcal{B}(\mathcal{H})$ is
a faithful nuclear representation, then $\pi|_A$ is a faithful, nuclear representation too,
whence exactness passes to sub-$C^*$-algebras.\\
The result is now immediate from \ref{redroe}.
\end{proof} 

The following proposition completes this section by generalising part of the above to
any bounded geometry metric space. 

\begin{roenucl} \label{roenucl}
Let $X$ be a bounded geometry metric space.\\
Then $C^*_u(X)$ is nuclear if and only if $X$ has property A.
\end{roenucl}
\begin{proof}
If $X$ has property A, exactly the same proof as used to show (2) $\Rightarrow$ (3)
in the above shows that $C^*_u(X)$ is nuclear.

Conversely, say $C_u^*(X)$ is nuclear and let $R,\epsilon>0$.  We will construct a
map as in \ref{equivdef}, (3).\\
The proof follows \cite{bo} 5.5.7.\\

We claim first that there exists a finite set $F\subseteq\mathcal{B}(l^2(X))$
of partial isometries such that for all $x,y$ with $d(x,y)\leq{}R$
there exists $V\in{}F$ such that $V\delta_x=\delta_y$.\\
As $X$ is bounded geometry, it is countable; order the elements as $(x_i)$.\\
Construct a sequence $(V_j)$ of operators on $l^2(X)$ inductively in $j$, by
inductively (in $i$) defining
\begin{displaymath}
V_j\delta_{x_i}=\left\{
\begin{array}{ll}
\delta_{x}, & x\in\bar{B}(x_i,R), ~V_j\delta_{x_{i'}}\not=\delta_x 
~\textrm{for any $i'<i$ and} \\
 & ~\delta_x\not=V_{j'}\delta_{x_i} ~\textrm{for any} ~j'<j \\
0, & ~\textrm{no such $x$ exists}.\end{array}\right.
\end{displaymath}
It is immediate that each $V_j$ is a partial isometry, and that for any
$x,y$ with $d(x,y)\leq{}R$ there exists $j$ with $V_j\delta_x=\delta_y$.\\
Moreover, if it was the case that $V_j\not=0$ for arbitrarily large $j$, we would
have that there exist $x_i$ with 
$\bar{B}(x_i,R)\cap\bar{B}(x_{i'},R)\not=\emptyset$ for an
arbitrarily large number of $i'<i$.  This in turn implies that there
exist $x_i$ such that $|B(x_i,2R)|$ is arbitrarily large, which contradicts
bounded geometry.\\
We may thus take $F$ to be the collection of all non-zero $V_j$.

Now, by nuclearity there exist $n$ and ucp maps 
$\phi:C^*_u(X)\to{}M_n(\mathbb{C})$, $\psi:M_n(\mathbb{C})\to{}C^*_u(X)$
such that $\|(\psi\circ\phi)(V)-V\|<\epsilon$ for all $V\in{}F$.\\
Set $\mathcal{E}=\mathbb{C}^n\otimes\mathbb{C}^n\otimes{}C^*_u(X)$, and
define a map $\langle,\rangle:\mathcal{E}\times\mathcal{E}\to{}C^*_u(X)$ by
setting
\begin{displaymath}
\langle{}x_1\otimes{}y_1\otimes{}T_1,x_2\otimes{}y_2\otimes{}T_2\rangle
=x_2\cdot{}x_1.y_2\cdot{}y_1.T_1^*T_2,
\end{displaymath}
where $\cdot$ is the usual scalar product on $\mathbb{C}^n$
(this makes $\mathcal{E}$ into what is called a \emph{Hilbert $C^*_u(X)$-module}; for us, however, it is 
just a notational convenience).\\
$\mathcal{E}$ has an associated norm defined by 
$\|e\|^2=\|\langle{}e,e\rangle\|$.\\
Note that $[\lambda_{ij}]\in{}M_n(\mathbb{C})$ acts on $\mathcal{E}$ as
$[\lambda_{ij}]\otimes[\lambda_{ij}]\otimes\psi([\lambda_{ij}])$.\\
Further, if $e_{ij}$ are the matrix units of $M_n(\mathbb{C})$, then
the matrix $[\psi(e_{ij})]_{i,j=1}^n$ is positive in $M_n(C^*_u(X))$, as
it is the image of positive $[e_{ij}]\in{}M_n(M_n(\mathbb{C}))$ under the
extension $\psi_n$ of ucp $\psi$.\\
Let $[b_{ij}]=[\psi(e_{ij})]^{\frac{1}{2}}$ and, letting $\{e_1,...,e_n\}$ be 
the standard basis for $\mathbb{C}^n$, define
\begin{displaymath}
\xi_{\psi}=\sum_{j,k}e_j\otimes{}e_k\otimes{}b_{kj}\in\mathcal{E}.
\end{displaymath}
Note then that for any $[\lambda_{ij}]\in{}M_n(\mathbb{C})$, we get that
\begin{displaymath}
\psi([\lambda_{ij}])=\langle\xi_{\psi},([\lambda_{ij}]\otimes{}1)\xi_{\psi}\rangle
\end{displaymath}
(this is a routine calculation).\\
Now, perturb the $b_{kj}$ so that they become finite propagation
operators, yet still satisfy
\begin{displaymath}
\psi(\phi(T))=\langle\xi_{\psi},(\phi(T)\otimes{}1)\xi_{\psi}\rangle
\end{displaymath}
for all $T$ of propagation at most $\max_{V\in{}F}\{$propagation$(V)\}$.\\
Write $\xi_{\psi}=\sum_j\xi_j\otimes{}b_j$, where 
$\xi_j\in\mathbb{C}^n\otimes\mathbb{C}^n$ and $b_j\in{}C^*_u(X)$ is assumed
to be finite propagation.\\ 
We now define a map $\eta:X\to{}l^2(X)$ by setting
\begin{displaymath}
\eta_x(y)=\|\sum_j\xi_j\langle\delta_y,b_j\delta_x\rangle_{l^2(X)}\|_{\mathcal{E}}.
\end{displaymath}
It remains to show that $\eta$ has the right properties.

Note first that as the operators $b_j$ are of finite propagation, and 
there are only
finitely many of them, there exists $S>0$ such that $\eta_x$ is supported in 
$B(x,S)$ for all $x$.\\
Moreover,
\begin{align*}
\|\eta_x\|^2 & 
=\sum_{y\in{}X}\|\sum_j\xi_j\langle\delta_y,b_j\delta_x\rangle\|^2 \\
& =\sum_{y\in{}X}\|\sum_{j,l}\xi_l\cdot\xi_j
\overline{\langle\delta_y,b_j\delta_x\rangle}\langle\delta_y,b_l\delta_x\rangle \\
& =\sum_{j,l}\xi_l\cdot\xi_j
\langle{}b_j\delta_x,\delta_y\rangle.\langle\delta_y,b_j\delta_x\rangle \\
& =\langle\delta_x,\langle\xi_{\psi},\xi_{\psi}\rangle_{\mathcal{E}}\delta_x\rangle
=\langle\delta_x,\delta_x\rangle=1.
\end{align*} 
Finally, let $d(x,y)\leq{}R$, and chose $V\in{}F$ such that $V\delta_x=\delta_y$.
As ucp maps are contractive, $\phi(V)$ is too, so we get that
\begin{align*}
\langle\eta_x,\eta_y\rangle_{l^2(X)}
& = \sum_{z\in{}X}\|\sum_j\xi_j\langle\delta_z,b_j\delta_x\rangle\|
\|\sum_l\xi_l\langle\delta_z,b_l\delta_y\rangle\| \\
& \geq\sum_{z\in{}X}\|\sum_j\xi_j\langle\delta_z,b_j\delta_x\rangle\|
\|(\phi(V)\otimes1)\sum_l\xi_l\langle\delta_z,b_l\delta_y\rangle\| \\
& =\langle\delta_x,\langle\xi_{\psi},
(\phi(V)\otimes1)\xi_{\psi}\rangle_{\mathcal{E}}\delta_y\rangle \\
& =\langle\delta_x,(\psi\circ\phi)(V)\delta_y\rangle
\geq1-\frac{\epsilon^2}{2}.
\end{align*}
It follows that $\|\eta_x-\eta_y\|_2$ is small.  We have shown that $\eta$ has the properties in \ref{equiv} (2) (for $p=2$), which completes the proof. 
\end{proof}

The proof above may seem a little intimidating on a first reading (it was to the author, anyway; it can be slightly simplified if one knows a little about
Hilbert $C^*$-modules, but we chose to skip this relatively technical prerequisite.).
The basic idea is 
fairly intuitive, however: use bounded geometry to find partial isometries that
allow one to move between elements at most $R$ apart, then use the ucp maps
given by nuclearity to find a map $\eta:X\to{}l^2(X)$ such that if 
$V\delta_x=\delta_y$, then $\eta_y\circ{}V$ is close in $l^2(X)$ to
$\eta_x$.  

It is also worth noting that the technique of the proof can be used to prove amenability
of $G$ from nuclearity of $C^*_r(G)$, except that here the partial isometries
$V$ are replaced by the unitaries $\lambda_g$ arising from the left regular representation (we are still
assuming a \emph{right}-invariant metric on $G$).  
This proof is actually somewhat
`cleaner' due to the greater homogeneity of the
action of the elements $\lambda_g$; this is one way of expressing the
difference between property A and amenability, which is the subject of
chapter \ref{amenable}. 

If one wants to be \emph{really} general, 
the natural forum is groupoids.  The
result that a groupoid is amenable if and only if its reduced $C^*$-algebra is
nuclear subsumes both of the preceding results; the proof is a natural (if one is comfortable with
groupoids)
generalisation.  See the comments at the end of \ref{aamen} for a very brief
overview, and the book \cite{adr}.

An interesting development in this analytic aspect of our discussion is the definition of 
\emph{ghost operators} (which form an ideal in $C_u^*(X)$) by Yu, who used them
give a new proof that certain spaces (expanders, section \ref{expander}) are counterexamples to coarse Baum-Connes; 
it turns out, however, that property A for $X$ is an obstruction to
this method of providing counterexamples.  See section 11.5.2 of \cite{roe1} and the papers \cite{cw2}, \cite{cw} for more discussion.

\pagebreak

\section {Property A for non-discrete spaces}
\label{beyond}

\subsection{Introduction}

Yu's original definition of property A is for discrete spaces, and a lot of
subsequent work has focused on bounded geometry spaces, which cover a large range
of examples, are relevant to coarse Baum-connes and for which the (generally easier
to manipulate)
reformulations of \ref{equiv} hold true.

We shall restrict ourselves to
Yu's original definition in the non-bounded geometry cases (this
will be relevant in \ref{nonace}, but nothing else from
now on), and discuss here
some possible extensions to non-discrete spaces.  Our motivations are:
\begin{enumerate}
\item Warped cones \cite{roe2}, which we shall discuss in section 
\ref{bswc}.
\item Index theory, a general discussion of which is beyond the scope of this 
piece.  The main point is, however, that often both coarse (large-scale) and
topological or differential (small-scale) structures are relevant to problems in this area. 
See e.g. \cite{roe3} for some ideas.
\end{enumerate}

We start with the following well-known principle of coarse geometry.  Recall
that a metric space is \emph{uniformly discrete} if there exists $\delta>0$
such that if $x,y$ are distinct points of $X$ then $d(x,y)\geq\delta$.
If we are interested in the constant $\delta$, we say the space is 
\emph{$\delta$-separated}.

\begin{diseq} \label{diseq}
A metric space $X$ is coarsely equivalent to a $\delta$-separated subspace for
any $\delta>0$.
\end{diseq}

\begin{proof}
Consider any chain of $\delta$-separated subsets, ordered by inclusion.  Its
union is also a $\delta$-separated subset, so we can apply Zorn's lemma to get a 
maximal such, say $Y$.  Note by maximality that every point in $X$ is within
$\delta$ of a point of $Y$.\\
Let $i:Y\to{}X$ be the inclusion map.  It is (trivially) bornologous,
and pulls back $B(i(y),R)$ to a subset of $B(y,R+\delta)$, so a coarse embedding.
Moreover, by the comment above, it is coarsely
surjective in the sense of \ref{cerem}, i.e. $\{B(y,\delta):y\in{}Y\}$
covers $X$.\\
Hence $i$ induces a coarse equivalence (note that an inverse up to closeness is
given by any map $r:X\to{}Y$ that takes a point $x$ to any point of $Y$ within
$\delta$ of it).
\end{proof}

The second part of the proof allows us to conclude
that if $X$ is a metric space,
then any coarsely dense (recall that $Y\subseteq{}X$ is \emph{coarsely
dense} if there exists $R>0$ such that $X\subseteq\cup_{y\in{}Y}B(y,R)$) subset
is coarsely equivalent to $X$.  This suggests the following definition:

\emph{A metric space $X$ has property A if and only if some (hence every) 
coarsely dense, discrete subset of $X$ has property A.}

Note that this still implies coarse embeddability in Hilbert space.  This
definition has two disadvantages, however:
\begin{itemize}
\item It can be difficult to work with: the argument in \ref{diseq} made heavy use of
the axiom of choice, and although coarsely dense subsets can sometimes be
constructed explicitly, to find one that is suitably `well-described' is not always easy.
\item It forces us to restrict our attention to discrete subspaces, which 
destroys metric (or topological, or differential...) information we might be interested in.
\end{itemize}

Bearing this in mind, in the next section, we look at a somewhat more tractable definition suggested by J. Roe \cite{roe2}.

\subsection{A concrete definition}

We start with the definition of bounded geometry for a 
non-discrete space.

\begin{nondbg} \label{nondbg}
A metric space $X$ has \emph{bounded geometry} if there exists 
$\epsilon\geq{}0$ such that for each $r>0$ there
exists $N_r\in\mathbb{N}$ such that for all $x\in{}X$, $B(x,r)$ can
be covered by at most $N_r$ balls of radius $\epsilon$.
\end{nondbg}

It is perhaps more usual to define bounded geometry in terms of
\emph{$\epsilon$-capacity}, which is the maximum cardinality of an
$\epsilon$-separated subset of a space.  $X$ then has bounded
geometry if there exists $\epsilon\geq0$ such that there is a uniform
bound on the $\epsilon$-capacity of $B(x,R)$ for all $x,R$.  The definitions
are essentially the same, however, as the minimal number of $\epsilon$-balls 
needed
to cover a space is bounded below by its $2\epsilon$-capacity and above
by its $\epsilon$-capacity.

Note that the definition for the discrete case (\ref{bg}) simply forces
$\epsilon=0$ in the above.  \ref{nondbg} is \emph{not} equivalent
to \ref{bg} for general discrete spaces.  It is, however,
true that a space has bounded geometry in the sense of this definition
if and only if it is coarsely equivalent to a space with bounded geometry in
our earlier sense.  Note also that this definition gives a property that is coarsely invariant
in full generality;
the earlier one is only coarsely invariant for maps between uniformly discrete spaces. 

Recall now that a \emph{proper} metric space $X$ is one in which
closed balls are compact.  It follows that $X$ is locally compact and 
$\sigma$-compact,
whence we can talk about the space of (complex) Radon measures on 
$X$.  It is the dual space of the $C^*$-algebra $C_0(X)$, and thus comes 
equipped with norm and weak-$*$ topologies.

\begin{nonda}[\cite{roe2}, 2.1] \label{nonda}
Let $X$ be a proper bounded geometry metric space.  $X$ has \emph{property A} 
if  for
any $R,\epsilon>0$ there exists a weak-$*$ continuous map
$\mu:X\to{}(C_0(X))^*$ such that:
\begin{enumerate}
\item $\|\mu_x\|=1$ for all $x\in{}X$;
\item $\mu$ has $(R,\epsilon)$ variation in the sense of \ref{equivdef};
\item there exists $S>0$ such that for each $x$, $\mu_x$ is supported in $\bar{B}(x,S)$
(i.e. if $f\in{}C_0(X)$, $f|_{\bar{B}(x,S)}=0$, then $\int_RFD\mu_x=0$).
\end{enumerate}
\end{nonda}

\begin{nondarem}
If $X$ is discrete, we have that bounded geometry 
(in the sense of \ref{bg})
implies proper, $C_0(X)=c_0(X)$, $(c_0(X))^*\cong{}l^1(X)$ and
weak-$*$ continuity of $\mu$ is automatic, so this reduces to the condition from
\ref{equiv}, (2) (for $p=1$).

We restrict ourselves to the bounded geometry case partly
because we are most interested in bounded geometry spaces, but also to remain
consistent with \ref{equiv}, which we only know gives equivalent conditions
to A for bounded geometry spaces.
\end{nondarem}

The next lemma provides further justification for our calling this new notion
`property A'.

\begin{nondacda}[\cite{roe2}, 2.2]
A proper metric space $X$ of bounded geometry has property A if and only if
some (hence every) coarsely dense discrete subspace has.  
\end{nondacda}

\begin{proof}
Suppose $X$ has property A, let $Y\subseteq{}X$ be uniformly discrete, and 
say $\mathcal{U}=\{B(y:r):y\in{}Y\}$ forms a cover for $X$.\\
Let $\{\phi_y\}_{y\in{}Y}$ be a partition of unity on $X$ subordinate to this
cover.\\
Let $R,\epsilon>0$, and let $\mu:X\to(C_0(X))^*$ be as in the definition
above.  By replacing $\mu_x$ with $|\mu_x|$, assume the image of $\mu$
consists of positive measures. Define $\xi:Y\to{}l^1(Y)$ by setting
\begin{displaymath}
\xi_{y_0}=\sum_{y\in{}Y}\left(\int_X\phi_y(x)d\mu_{y_0}(x)\right)\delta_y,
\end{displaymath}
whence
\begin{displaymath}
\|\xi_{y_0}\|_1=\int_X\sum_{y\in{}Y}\phi_y(x)d\mu_{y_0}(x)=1.
\end{displaymath}
Moreover, if $\mu_y$ is supported in $B(y,S)$, $\xi_y$ is supported in
$B(y,S+2r)$, and if $d(y_1,y_2)<R$,
\begin{align*}
\|\xi_{y_1}-\xi_{y_2}\|_1
& =\sum_{y\in{}Y}\left|\int_X\phi_y(x)d(\mu_{y_1}-\mu_{y_2})(x)\right|
\leq{}\int_X\sum_{y\in{}Y}\phi_y(x)d|\mu_{y_1}-\mu_{y_2}|(x) \\
& =\|\mu_{y_1}-\mu_{y_2}\|<\epsilon.
\end{align*}
Conversely, assume $Y$ is a uniformly discrete, coarsely dense
subspace of $X$ having
property A in the sense of \ref{equiv} (3) for $p=1$.\\
Say again that $\mathcal{U}=\{B(y:r):y\in{}Y\}$ forms a cover for $X$ and let
$\{\phi_y\}_{y\in{}Y}$ be a partition of unity subordinate to $\mathcal{U}$.\\
Let $R,\epsilon>0$ and let $\xi:Y\to{}l^1(Y)$ be as in \ref{equiv} (3)
with respect to the parameters $R+2r,\epsilon, S$.\\
Define $\mu:X\to(C_0(X))^*$ by setting
\begin{displaymath}
\mu_x(f)=\sum_{y\in{}Y}\phi_y(x)\xi_y(f) ~\textrm{where}
~\xi_y(f)=\sum_{y'\in{}Y}f(y')\xi_y(y').
\end{displaymath}
Then $\|\mu_x\|=1$ for all $x\in{}X$, $\mu$ depends (weak-$*$) continuously 
on $x$ and $\mu_x$ is supported in $B(x,S+2r)$.\\
Finally, for $d(x,x')<R$ consider
\begin{align*}
\|\mu_x-\mu_{x'}\|
& =\left\|\sum_{y\in{}Y}\phi_y(x)\xi_y-\sum_{y\in{}Y}\phi_y(x')\xi_y\right\| \\
& =\left\|\sum_{y\in{}B(x,r)}\phi_y(x)\xi_y-
\sum_{y\in{}B(x',r)}\phi_y(x')\xi_y\right\| \\
& \leq\sup_{y\in{}B(x,r), y'\in{}B(x',r)}\|\xi_y-\xi_{y'}\|.
\end{align*}
The last inequality follows as both sums in the previous line are convex combinations
of $\xi_y$ (respectively, $\xi_{y'}$) such that $y\in{}B(x,r)$ (respectively
$y'\in{}B(x',r)$).\\
This in turn is bounded by
\begin{displaymath}
\sup_{d(y,y')<2r+R}\|\xi_y-\xi_{y'}\|\leq\epsilon,
\end{displaymath}
by assumption on $\xi$, as if $d(x,x')<R$, $y\in{}B(x,r)$, $y'\in{}B(x',r)$, then
$d(y,y')<R+2r$.
\end{proof}

\begin{cinvcor}
Say $X$ is proper, and has bounded geometry and property A in the sense of 
\ref{nonda}.\\ 
Then any discrete subspace of $X$ has A in the sense of \ref{propa} and
$X$ coarsely embeds into Hilbert space.\\
Moreover, property A as defined in \ref{nonda} is coarsely invariant.
\end{cinvcor}
\begin{proof}
Some uniformly discrete coarsely dense subspace $Y\subseteq{}Z$ has
bounded geometry in the sense of \ref{bg}, whence it has property A in the
sense of the original definition.\\
Coarse embeddability follows from coarse equivalence of $X$ and $Y$; any discrete
subspace of $X$ is coarsely equivalent to a subspace of $Y$, so has A by 
\ref{cinvrem}; coarse invariance follows from \ref{cinv} as any coarsely 
equivalent spaces have coarsely equivalent discrete subspaces. 
\end{proof}

As an example, we can now say that Euclidean space has property A, as $\mathbb{R}^n$ 
is coarsely equivalent to $\mathbb{Z}^n$, which we know has property A from
\ref{integers} and \ref{produ}. 

We finish this section with the following reformulation of \ref{nonda};
cf. \ref{equiv}, (8).

\begin{nonda2} \label{nonda2}
Let $X$ be a bounded geometry proper metric space.\\
$X$ has property A in the sense of \ref{nonda} if and only if for each
$R,\epsilon>0$ there exists a finite propagation symmetric continuous positive 
type kernel on $X$ of $(R,\epsilon)$ variation.
\end{nonda2}

\begin{proof}
It will be sufficient to prove that such a kernel exists on $X$ if and only if
it exists on some coarsely dense discrete subspace.
If so, then its restriction to any discrete subspace has the
right properties.\\
Conversely, say $Y$ is a coarsely dense discrete subspace; say
$\mathcal{U}=\{B(y,r):y\in{}Y\}$ covers $X$, and let $\{\phi\}_{y\in{}Y}$ be a
partition of unity on $X$ subordinate to $\mathcal{U}$.\\
Let $R,\epsilon>0$, and let $k$ be a kernel on $Y$ with the right properties
with respect to the parameters $R+2r,\epsilon$.\\
Then $\hat{k}:X\times{}X\to\mathbb{R}$ defined by
\begin{displaymath}
k(x,x')=\sum_{(y,y')\in{}Y\times{}Y}k(y,y')\phi_y(x)\phi_{y'}(x')
\end{displaymath}
has the right properties.
\end{proof} 

We will use these versions of property A in our discussion of warped cones 
(see \cite{roe2}) in section \ref{bswc}).

\pagebreak

\section{Connections with amenability} \label{amenable}

\subsection{Introduction} \label{amenableintro}

Amenability is a property of groups, originally introduced by
Von Neumann in his work on so-called paradoxical (=non-amenable) groups.  Like property A,
it controls growth and connectivity to some extent, although is much stronger (A is 
often called a `weak form of amenability').

We start this chapter by defining amenability, and also a-T-menability, for
countable discrete groups.  Although amenability can be defined for any metric
space, we restrict ourselves to this case, as it is here that the relationship
with property A is most apparent (in particular, amenability does \emph{not} imply A for
general  (bounded geometry) metric spaces; see for example \cite{roe1}, chapter 3 for discussion and further references).  
Both amenability and a-T-menability can usefully be defined for a larger class of groups (e.g. locally compact groups; see
for example appendix G of \cite{bhv} and \cite{ccjjv}).

Section \ref{aamen} then starts to 
explore the relationship between A and amenability.
Section \ref{bswc} discusses two
specific constructions: box spaces and warped cones.
A group is an ingredient in both constructions,
and we get (modulo some technicalities) that the 
constructed space has A if and only if the group is amenable.
We use this in section \ref{propT} to construct examples of spaces that do
not have property A. 

The next definition gives some equivalent formulations of amenability.

\begin{amenable} \label{amendef}
Let $G$ be a countable discrete group.\\
$G$ is \emph{amenable} if any of the following equivalent conditions hold:
\begin{enumerate} 

\item
For all 
$R,\epsilon>0$ there exists a finitely supported map $\xi\in{}l^1(G)$ such that 
$\|\xi\|_1=1$ and
\begin{displaymath}
\|g\xi-\xi\|_1<\epsilon ~\textrm{for all} ~g\in\bar{B}(e,R),
\end{displaymath}
where we write $(g\xi)(h)=\xi(g^{-1}h)$.

\item 
For any $1\leq{}p<\infty$ and all 
$R,\epsilon>0$ there exists a finitely supported map $\eta\in{}l^p(G)$ such that 
$\|\eta\|_p=1$ and
\begin{displaymath}
\|g\eta-\eta\|_p<\epsilon ~\textrm{for all} ~g\in\bar{B}(e,R).
\end{displaymath}

\item
For any $R,\epsilon>0$ there exists $S>0$ and
a normalised positive type function (see \ref{posfun}) 
$\phi:G\to\mathbb{C}$ such that:
\begin{enumerate}
\item $|1-\phi(g)|<\epsilon$ for $g\in\bar{B}(e,R)$;
\item $\phi$ is supported in $\bar{B}(e,S)$.
\end{enumerate}
We call (a) \emph{$(R,\epsilon)$-variation} for functions of positive type.

\item
There exists $\phi$ as in (3) such that convolution with $\phi$ defines a positive
operator $T_{\phi}\in{}C^*_r(G)$.

\item 
There exists a 
positive linear functional on $l^{\infty}(G)$, which we denote by 
\begin{displaymath}
M:f\mapsto\int_Gf(g)dg,
\end{displaymath}
which sends the constant function $\mathbf{1}$ to $1$, and is 
invariant in the sense that 
\begin{displaymath}
\int_Gf(gh)dg=\int_Gf(hg)dg=\int_Gf(g)dg
\end{displaymath}
for all $h\in{}G$, all $f\in{}l^{\infty}(G)$.\\
Such an $M$ is called an \emph{invariant mean}.
\end{enumerate}
\end{amenable}

\begin{proof}
We prove (1) $\Leftrightarrow$ (2), (2) $\Rightarrow$ (3) $\Rightarrow$ 
(4) $\Rightarrow$ (2) and
(1) $\Leftrightarrow$ (5).  We only sketch the proofs, as in most cases
they are similar to \ref{equiv}, and amenability is not our main focus.\\

(1) if and only if (2):\\
Set $\eta(g)=\sqrt[p]{|\xi(g)|}$, or $\xi(g)=|\eta(g)|^p$.\\

(2) implies (3):\\
Set $\phi(g)=\langle\eta,g\eta\rangle$.\\

(3) implies (4):\\
This is a simple check.\\

(4) implies (2):\\
Use the square root of $T_{\phi}$ to define a map into $l^2(G)$ with the right
properties; cf. \ref{equiv}, (9) implies (2).\\

(5) if and only if (1):\\
Use (1) to find a sequence $\xi_n$ such that $|g|\leq{}n$ forces $\|g\xi_n-\xi_n\|_1<1/n$.
As $l^1(G)\subseteq{}l^{\infty}(G)^*$, this sequence lies in weak-$*$ compact
$\bar{B}_{l^{\infty}(X)}(0,1)$.  Any weak-$*$ limit point is a left-invariant
mean, say $M_l$ (left invariance means $M_l(gf)=f$ for all $f$, $g$).\\
For any $f\in{}l^{\infty}(G)$, define an element $\mu_f\in{}l^{\infty}(G)$ by
setting $\mu_f(g)=M_l(fg)$, where $fg$ denotes the image of $f$ under the right action of $g$,
i.e. $fg(h)=f(hg^{-1})$.  Define $M:l^{\infty}(G)\to\mathbb{C}$ by
setting $M(f)=M_l(\mu_f)$; as $\mu_{fg}=g^{-1}\mu_f$, this is an
invariant mean.

To get the converse, note that the space of finitely supported $l^1$ functions on $G$
is weak-$*$ dense in $l^{\infty}(G)^*$, and thus an invariant mean can be 
approximated by such a function; a sufficiently good weak-$*$ approximation 
(using finiteness of closed balls in $G$) satisfies the conditions in (1).
\end{proof}

(5) is perhaps the most common definition.  It was popularised
M. Day, who started the modern study of amenability in the 1950s.  (1) is sometimes called the
\emph{Reiter condition}; by analogy \ref{equiv}, (2) with respect to $p=1$ is sometimes called the
\emph{weak Reiter condition}. 

\begin{amenimpa} \label{amenimpa}
If $G$ is a (countable discrete) amenable group, then $G$ has property A.
\end{amenimpa}
\begin{proof}
Let $R,\epsilon>0$, and say $\phi$ satisfies the conditions in the definition
of amenability, part (3).\\
Define a normalised, finite propagation, positive type kernel 
$k:G\times{}G\to\mathbb{C}$ of $(R,\epsilon)$ variation by
setting $k(g,h)=\phi(g^{-1}h)$.\\
The relationships
\begin{align*}
0\leq & \phi(g^{-1}h)+\phi(h^{-1}g)+2\phi(e) \\
0\leq & i\phi(g^{-1}h)-i\phi(h^{-1}g)+2\phi(e)
\end{align*} 
show that $k$ is self-adjoint; they imply that 
$\phi(g^{-1}h)+\phi(h^{-1}g)$, $i\phi(g^{-1}h)-i\phi(h^{-1}g)$ are real,
whence $\phi(g^{-1}h)=\overline{\phi(h^{-1}g)}$.\\
Finally, set $\hat{k}=$Re$k$; this satisfies the conditions in
\ref{equiv}, (8).
\end{proof}

The lemma includes a proof that a finitely supported (normalised) 
function of positive type
on $G$ defines a (normalised) finite propagation kernel of positive 
type.  The converse is far from true, as the many examples of
non-amenable property A groups testifies.  
For example, non-abelian free groups are not amenable, but do have property A by
\ref{integers}.
In general, only a \emph{left-invariant} 
positive type kernel (i.e. $k(g_1,g_2)=k(gg_1,gg_2)$ for any $g,g_1,g_2\in{}G$) 
gives rise to a positive type function; we do, however, get the following partial
converse.

\begin{kerfun} \label{kerfun}
If $G$ is amenable, any bounded positive type kernel gives rise to a bounded,
positive type function.\\
Moreover, the construction preserves the properties of being normalised and of
$(R,\epsilon)$ variation, and creates finitely supported functions out of finite
propagation kernels.
\end{kerfun}

\begin{proof}
The basic idea is to use the invariant mean characterisation of \ref{amendef}, (5) to
construct a left-invariant kernel.\\
Let $k$ be of positive type, then, and for any $g\in{}G$ consider 
the function $k(g\cdot,\cdot):h\mapsto{}k(gh,h)$.\\
As $k$ is bounded, it is an element of $l^{\infty}(G)$, so it makes sense to define
\begin{displaymath}
\phi(h)=\int_Gk(h^{-1}g,g)dg,
\end{displaymath}
where the integral is the invariant mean on $G$.
For elements $\lambda_1,...,\lambda_n$ of $\mathbb{C}$, and $g_1,...,g_n$ of $G$, we get
\begin{align*}
\sum_{i,j}\lambda_i\overline{\lambda_j}\phi(g_i^{-1}g_j)
& = \sum_{i,j}\lambda_i\overline{\lambda_j}\int_Gk(g_jg_i^{-1}g,g)dg
=\sum_{i,j}\lambda_i\overline{\lambda_j}\int_Gk(g_i^{-1}g,g_j^{-1}g)dg \\
& = \int_G\sum_{i,j}\lambda_i\overline{\lambda_j}k(g_i^{-1}g,g_j^{-1}g)dg\geq0
\end{align*}
The penultimate equality follows by invariance of the mean, and the inequality follows by
the mean's positivity, and that of the function 
$h\mapsto{}\sum_{i,j}\lambda_i\overline{\lambda_j}k(g_i^{-1}\cdot,g_j^{-1}\cdot)$.\\
The remaining properties are routine checks.
\end{proof}

This illuminates an important general principle when dealing with amenable groups: 
one can use
the invariant mean to `average out' a function over the group. We do this in 
\ref{aeqamen}, while in
\ref{boxamen} and \ref{wcaimpamen} we use respectively the Haar measure on a 
compact group and
counting measure on a finite group to fulfill similar r\^{o}les.

The similarities between (1), (2) in \ref{amendef} and (2), (3) in
\ref{equiv} are also fairly apparent; simply replace $\xi_g$ by $g\xi$.  These
relationships suggest that amenability is somewhat more dependent on the
group structure than property A; it is partly for this reason that A is
sometimes called a `non-equivariant' form of amenability.  Nontheless
(and perhaps surprisingly), amenability is a coarse invariant of 
(countable discrete) groups (see e.g. chapter 3 of \cite{roe1}). 

The following defines a weak form of amenability; the relationship between the two is analogous to that between property A and coarse embeddability.

\begin{atmen} \label{atmen}
A (countable discrete) group $G$ is \emph{a-T-menable} (or has 
\emph{the Haagerup property}) if for all $(R,\epsilon)>0$ there
exists a normalised positive type function $\phi$ on $G$ such that:
\begin{enumerate}
\item $|1-\phi(g)|<\epsilon$ for $g\in{}B(e,R)$;
\item $\lim_{S\to\infty}\sup\{|\phi(g)|:d(e,g)\geq{}S\}=0$. 
\end{enumerate} 
\end{atmen}

An alternative, more geometric, characterisation states that $G$ is a-T-menable if
and only if it acts properly affinely isometrically on Hilbert space. 

There are many groups that have property A, but are not a-T-menable.  For example,
many linear groups have property T (see section \ref{propT}) so cannot be
a-T-menable. They do have property A, however, by the result of \cite{ghw}.  The
converse is not known: for example, Thompson's group $F$ is 
a-T-menable (this is proved in \cite{far}), but not known to have property A.  This is
perhaps surprising, as a-T-menability is usually thought of as `slightly weaker' than 
amenability, while property A is `much weaker'.  It may be the case, however, that all one
can really say is that these two generalisations of amenability
that proceed in rather different directions.

The same argument as in lemma \ref{amenimpa}, combined with proposition
\ref{eh2}, part (3), shows that a-T-menability implies coarse embeddability
in Hilbert space.  This connection, i.e. the link between the formulation
of A / coarse embeddability in terms of positive definite kernels, and that of
amenability / a-T-menability in terms of positive definite functions, is part of
the reason for the slogan `property A is to coarse embeddability as
amenability is to a-T-menability'.  Another point along these lines is
looked at in section \ref{bswc}.

It is worth noting, however, that while many groups are known to be
a-T-menable that are not amenable (see the book \cite{ccjjv}), 
the only known examples of coarsely
embeddable spaces that do not have property A are those given in section
\ref{nonace}. 

Finding a \emph{group} that coarsely embeds, but does not have A would be
extremely difficult given the current state of the theory.  The only known
examples of non-A groups are the non-coarsely embeddable groups
discussed in \ref{nonagp} (see \cite{gr2}), and these are not well understood at time of writing (in fact, there is not even a complete proof of their existence).  One possible contender is
Thompson's group $F$ as mentioned above (and some similar so-called \emph{diagram groups}), but even whether or not this is amenable is a long-standing open problem.

\subsection{The relationship between property A and amenability}
\label{aamen} 

We will prove some relationships between property A and amenability
(for groups) in the hope that they clarify the relationship
between the two concepts.  The table below,
each column of which consists of equivalent conditions, summarises this.\\

\begin{tabular}{p{5cm}|p{5cm}}
$G$ has property A & $G$ is amenable \\
\hline
$G$ admits a sequence of normalised kernels of positive type and finite 
propagation that converge pointwise to $\mathbf{1}$. 
& $G$ admits a sequence of normalised functions of positive type and finite 
support that converge pointwise to $\mathbf{1}$. \\
\hline
$C^*_u(G)$ is nuclear. & $C^*_r(G)$ is nuclear. \\
\hline
$G$ admits an amenable action on some compact Hausdorff space
& $G$ admits an amenable action on a point.
\end{tabular}

The results in the first row were established earlier in this section, and the results
of \ref{c*facts1}, \ref{c*facts2}, 
help to establish the connection between
the second and third rows.

The following gives one way to prove the theorem of Lance that amenability and nuclearity of
$C^*_r(G)$ are equivalent.  

\begin{nuclamen}
A countable discrete group $G$ is amenable if and only if $C^*_r(G)$
is nuclear.
\end{nuclamen}
\begin{proof}
The proof of \ref{roenucl} can be adjusted to show that nuclearity of $C_r^*(G)$ implies
amenability of $G$; after replacing the finite set of partial isometries in that proof
with the finite set $\{\lambda_g:|g|\leq{}R\}$.\\
Conversely, the proof of \ref{nuclexact} can be rewritten, again with `in
$\mathbb{C}[G]$' replacing `finite propagation' to get a proof that amenability
implies nuclearity.
\end{proof}

This is perhaps not the easiest way to prove the above result, but it does highlight some
of the similarities we are discussing.

The final row of the table provides a good
intuitive connection between property A and amenability - the definition of an
amenable action (below) simply `spreads out' \ref{amendef} around the $G$-space. 
It is worth remarking here that it is relatively easy to prove that $C^*_r(G)$ is
exact if and only if $G$ admits an amenable action on a compact topological space
(see \cite{bo} or \cite{adr}), which gives the left-hand equivalence.  We give the
original, direct proof below \cite{hr}.

\begin{amenact} \label{amenact} 
A (countable discrete) group $G$ is said to act \emph{topologically amenably}
on the compact Hausdorff space $X$ if for all $R,\epsilon>0$ there exists
a weak-$*$ continuous map $\xi:X\to{}l^1(G)$ such that 
\begin{enumerate}
\item $\|\xi_x\|_1=1$ for all $x\in{}X$;
\item for all $g\in\bar{B}(e,R)$ and all $x\in{}X$, 
$\|g\xi_x-\xi_{gx}\|_1<\epsilon$.
\end{enumerate}
\end{amenact}

It is clear from this and \ref{amendef}, (1) that any action at all of an amenable group is
an amenable action.  On the other hand, if $G$ acts amenably on a one point space $X=\{x\}$, then $\xi_x\in{}l^1(G)$
as in the above definition satisfies \ref{amendef} (1) (technically, we also need an 
approximation argument to show that $\xi_x$ may be assumed finitely supported, but this is not hard), and thus shows amenability of $G$.
In particular, then, a group $G$ is amenable if and only if it acts amenably on a point.  
This completes the equivalences in the right hand column of the table appearing at the start of this section.

In order to be able to complete the equivalences in the left-hand column,
we first review some facts about the 
\emph{Stone-\v{C}ech compactification} of a topological space.  We will
not prove its existence here; this is dealt with in most textbooks
on general topology. 

\begin{stcech}
Let $X$ be a completely regular (recall this means that given any point $x$ and
any neighbourhood $U$ of $x$, there is a continuous function $f:X\to{}[0,1]$ 
such that $f(x)=0$ and $f(y)=1$ for all $y\in{}X\backslash{}U$) Hausdorff
topological space.\\
Then there exists a compact Hausdorff space $\beta{}X$ that contains $X$ as a dense 
subspace, and such that if $f:X\to{}Y$ is any continuous map of $X$ into a compact
Hausdorff space $Y$, then there exists a unique continuous extension of $f$ to
$\beta{}X$.
\end{stcech}

This applies to a countable discrete group $G$ to give a 
compactification $\beta{}G$.  It is also the case that $G$ has a natural
(left) translation action on $\beta{}G$, which is arrived at by extending the natural
left translation of $G$ on itself to $\beta{}G$; we can do this by the 
universal property.

The next three lemmas are basically enough to establish the result.  We start
by characterising the existence of amenable actions on the Stone-\v{C}ech
compactification of a group.

\begin{amenact1}
A (countable discrete) group $G$ admits an amenable action on $\beta{}G$ 
if and only if
for any $R,\epsilon>0$ there exists a map $\zeta:G\to{}l^1(G)$ such that:
\begin{enumerate}
\item $\|\zeta_g\|_1=1$ for all $g\in{}G$;
\item for all $g,h\in{}G$ with $|g|\leq{}R$ 
$\|g\zeta_h-\zeta_{gh}\|_1<\epsilon$;
\item the image of $\zeta$ is contained in a weak-$*$ compact subset of  $S(l^1(G))$,
the unit sphere of $l^1(G)$.
\end{enumerate}
\end{amenact1}

\begin{proof}
Say first that the (left translation) action of $G$ on $\beta{}G$ is amenable,
so there exists a map $\xi:\beta{}G\to{}l^1(G)$ satisfying the properties
in \ref{amenact}; set $\zeta=\xi|_G$.\\
$\zeta$ trivially has properties (1) and (2) above, and $\zeta(G)$ sits in 
$\xi(\beta{}G)\subseteq{}S(l^1(G))$, which is weak-$*$ compact by weak-$*$ continuity
of $\xi$.\\ 
Conversely, say $\zeta$ satisfies the properties above for some $(R,\epsilon)$.\\
Note that $\zeta$ is weak-$*$ continous (it is a function out of a discrete set!),
and Im$(\zeta)$ is contained in some weak-$*$ compact $K\subseteq{}S(l^1(G)$, whence
we can extend $\zeta$ to $\xi:\beta{}G\to{}K$ by the universal property of 
$\beta{}G$.\\
$\xi$ trivially has property (1) in \ref{amenact}, and note that
for all $g\in{}G$ with $|g|\leq{}R$, we get that
\begin{displaymath}
\sup_{h\in{}G}\|g\zeta_h-\zeta_{gh}\|_1
=\sup_{x\in{}\beta{}G}\|g\xi_x-\xi_{gx}\|_1,
\end{displaymath}
as $G$ is dense in $\beta{}G$, $x\mapsto{}g\xi_x-\xi_{gx}$ is weak-$*$
continuous, and because if $x_{\lambda}\to{}x$ in the
weak-$*$ topology, $\|x\|_1\leq{}\sup_{\lambda}\|x_{\lambda}\|_1$.
\end{proof}

We now give a reformulation of property A for countable discrete groups that
is very similar to the conditions in the above lemma.

\begin{amenact2}
A countable discrete group $G$ has property A if and only if for any
$R,\epsilon>0$ there exists a map $\theta:G\to{}l^1(G)$ such that:
\begin{enumerate}
\item $\|\theta_g\|_1=1$ for all $g\in{}G$;
\item for all $g$ with $|g|\leq{}R$, $\|g\theta_h-\theta_{gh}\|_1<\epsilon$;
\item there exists $S>0$ such that for each $g$, $\theta_g$ is supported in
$\bar{B}(e,S)$.
\end{enumerate}
\end{amenact2}

\begin{proof}
Assume that $G$ has property A (with respect to some
left-invariant, bounded geometry metric $d$), as formulated in \ref{equiv}, 
(3) for $p=1$. Take a map $\xi$ with the properties given there for any $R,\epsilon>0$.\\
$(R,\epsilon)$ variation of $\xi$ says that for all $g,h\in{}G$ such that
$d(g,h)\leq{}R$, $\|\xi_g-\xi_h\|_1<\epsilon$, whence
for any $g,h\in{}G$ with $|g|\leq{}R$, we have that $\|\xi_h-\xi_{hg}\|_1<\epsilon$
by left-invariance of $d$.\\
Set $\theta_g=g\xi_{g^{-1}}$.  Then for any $g,h\in{}G$ with $|g|\leq{}R$,
\begin{displaymath}
\|g\theta_h-\theta_{gh}\|_1=\|gh\xi_{h^-1}-
gh\xi_{h^{-1}g^{-1}}\|_1
=\|\xi_{h^{-1}}-\xi_{h^{-1}g^{-1}}\|_1<\epsilon.
\end{displaymath}
Reversing this process gives the converse.
\end{proof}

We now have that property A gives an amenable action of $G$ on
 $\beta{}G$. To get the converse, we use the following:

\begin{amenact3} \label{amenact3}
Let $X$ be a discrete space.  For any weak-$*$ compact subset $K$ of the unit
sphere of $l^1(X)$ there exists a finite subset $F_0\subseteq{}X$ such that
$\|\xi-\xi|_{F_0}\|_1<\epsilon$ for all $\xi\in{}K$.
\end{amenact3}

\begin{proof}
Let $\epsilon>0$, and consider the sets
\begin{displaymath}
U_F=\{\xi\in{}S(l^1(X)):\|\xi|_F\|_1>1-\epsilon\},
\end{displaymath}
as $F$ ranges over finite subsets of $X$.\\
They form an open cover of $S(l^1(X))$, and thus there exists a finite
collection $U_{F_1},...U_{F_n}$ that covers $K$.  Take $F_0=\cup_{i=1}^nF_i$. 
\end{proof}

\begin{amenacta}[\cite{hr}, 1.1] \label{amenacta}
A countable discrete group $G$ has property A if and only if it admits an
amenable action on some compact Hausdorff space.
\end{amenacta}
\begin{proof}
The three lemmas give that $G$ has property A if and only if it admits an
amenable action on $\beta{}G$.
If there exists an amenable action on $\beta{}G$, there is certainly an
amenable action on some compact Hausdorff space.\\
Conversely, say $G$ admits an amenable action on the compact Hausdorff space $X$.\\
Picking any point $x_0\in{}X$, we get a (continuous, equivariant) map $f:G\to{}X$ sending $g$ to 
$gx_0$; this
extends to a continuous map $\tilde{f}:\beta{}G\to{}X$ by the universal property.  Moreover, $\tilde{f}$ is $G$-equivariant, as $f$ is and
$G$ is dense in $\beta{}G$.\\
Given a map $\xi:X\to{}l^1(G)$ as in the definition of amenable action, we can
now compose to get $\xi\circ\tilde{f}:\beta{}G\to{}l^1(G)$; using equivariance of $\tilde{f}$, it is immediate that this
also has properties (1) and (2) from \ref{amenact}.
\end{proof}

To conclude this section we state without proof some results about groupoids. 
A general discussion is contained in \cite{adr}, but is beyond the scope of this
piece.

The theory of \emph{amenable groupoids} provides a framework which ties together many of the ideas
in this section (recall that a \emph{groupoid} is (roughly!) a group whose
multiplication is not everywhere defined).  In general
a groupoid $\mathcal{G}$ is amenable if and only if the associated \emph{reduced groupoid $C^*$-algebras}
$C^*_r(\mathcal{G})$ is nuclear (modulo
concerns about second countability).   One can define the following
groupoids, with associated algebras:
 
\begin{itemize}
\item A group is a special case of a groupoid.  In this case amenability of the 
groupoid is the same thing as amenability of the group, and the reduced group and
groupoid $C^*$-algebras are equal.

\item The \emph{partial translation groupoid} of a coarse space (for us, coarse
equivalence class of metric spaces), $\mathcal{G}(X)$. This was introduced
in \cite{sty}; see also \cite{roe1}, chapter 10. It is amenable if and only if
$X$ has property A, and $C^*_r(\mathcal{G}(X))\cong{}C^*_u(X)$.

\item The groupoid associated to an action of a group $G$ on a compact Hausdorff space $X$.  
It is amenable if and only if the action is.  
In general, $C^*_r(\mathcal{G})\cong{}C(X)\rtimes_rG$ (the right hand side is the
\emph{reduced crossed product $C^*$-algebra}).
In the special case of a countable discrete group acting on its Stone-\v{C}ech 
compactification
one also has $C(\beta{}G)\rtimes_rG\cong{}l^{\infty}(G)\rtimes_rG\cong{}C^*_u(G)$.
\end{itemize}

\subsection{Box spaces and warped cones} \label{bswc}

As mentioned in the introduction to this chapter, the box space and warped cone constructions 
provide further links between property A and amenability.  In particular, they allow the construction of 
non-property A spaces (see section \ref{propT}).
We start by looking at box spaces, which are the simpler of these
two constructions.

\begin{finquot}
Let $G$ be a (countable discrete) group equipped with a left-invariant,
bounded geometry metric $d$.  
Let $\pi:G\to{}F$ be the homomorphism onto a finite quotient of $G$.  We 
metrise $F$ by setting
\begin{displaymath}
d_F(f,f')=\min\{d(g,g'):g,g'\in{}G, \pi(g)=f, \pi(g')=f'\}
\end{displaymath}
(cf. the proof of theorem \ref{ext}), giving $F$ a left-invariant metric.
\end{finquot}

Recall that a group $G$ is \emph{residually finite} if homomorphisms from
$G$ to finite groups separate elements (equivalently, the intersection of all
normal finite index subgroups is trivial).

\begin{boax}\label{boax}
Assume now that $G$ is a residually finite countable discrete group.\\
Let $K_1\supseteq{}K_2\supseteq{}K_3\supseteq...$ be a decreasing sequence of finite index normal subgroups 
of $G$ such that $\cap_{n=1}^{\infty}K_n=\{1\}$.\\
We define the \emph{box space of $G$ with respect to the family $(K_n)_n$} to have underlying space the disjoint union of the finite 
quotients $F_n=G/K_n$. It is equipped with any metric that restricts to $d_{F_n}$ as above on each finite
quotient $F_n$, and keeps distinct quotients further apart than the larger
of their diameters (any two such metrics are coarsely equivalent).
\end{boax}

We write $\Box{}G$ for the box space of $G$ with respect to some sequence $(K_n)$ as in the above (this sequence is kept implicit).

Our interest in this construction stems mainly from the following result.

\begin{boxamen} \label{boxamen}
Let $G$ be a residually finite countable discrete group, and $(K_n)_{n=1}^{\infty}$ any family of finite index subgroups as in \ref{boax}.
Let $\Box{}G$ be the box space with respect to this family.\\
Then $G$ is amenable if and only if $\Box{}G$ has property A.
\end{boxamen}

\begin{proof}
Suppose first that $\Box{}G$ has property A.\\ 
Then for any $R,\epsilon>0$ there exists a normalised, symmetric kernel 
$k:\Box{}G\times\Box{}G\to\mathbb{R}$ of 
positive type, $(R,\epsilon)$ variation and finite propagation $S$.\\
For each $n$, let $F_n=G/K_n$ and define $\phi_n:F_n\to\mathbb{R}$ by setting
\begin{displaymath}
\psi_n(f)=\frac{1}{|F_n|}\sum_{f'\in{}F_n}k(f',f'f)
\end{displaymath}
Note that $\psi_n$ is normalised, and of positive type as if $f_1,...,f_k\in{}F_n$,
$\lambda_1,...,\lambda_k\in\mathbb{C}$,
\begin{align*}
\sum_{i,j}\lambda_i\overline{\lambda_j}\psi_n(f_i^{-1}f_j)
& =\sum_{i,j}\lambda_i\overline{\lambda_j}\frac{1}{|F_n|}
\sum_{f\in{}F_n}k(f,ff_i^{-1}f_j) \\
& =\frac{1}{|F_n|}\sum_{f\in{}F_n}
\sum_{i,j}\lambda_i\overline{\lambda_j}k(ff_i,ff_j)\geq{}0,
\end{align*}
using that $k$ is of positive type, and lemma \ref{krecom}.
We define $\phi_n:G\to\mathbb{R}$ to be the composition of $\psi_n$ and the
quotient $\pi_n:G\to{}F_n$; each $\phi_n$ is
normalised and of positive type, so bounded by one.\\
As $G$ is countable, we can thus use a diagonal argument to find a subsequence
of $(\phi_n)$ that converges pointwise to some $\phi$.\\
$\phi$ is then of positive type, normalised and symmetric.  Moreover:
\begin{itemize}
\item If $d(g,h)>S$, then $d(\pi_n(g),\pi_n(h))>S$ for all but finitely many
$n$ by choice of the family $(K_n)$ and the metric on $F_n$; if
$n$ is suitably large, $K_n$ will contain none of the (finitely many) non-identity
elements in
$B(e,2S)$, and it follows that for any $g\in{}B(e,S)$, the only element
of the coset $gK_n$ of length less than $S$ is $g$.\\  
Hence $\psi_n(\pi_n(g)^{-1}\pi_n(h))=0$ for all but finitely
many $n$, as $k$ is of propagation $S$, whence $\phi$ is supported on
$B(e,S)$, so in particular is finitely supported.
\item If $d(g,h)<R$, then $d(\pi_n(g), \pi_n(h))<R$ for all $n$, as
each $\pi_n$ is contractive.  Hence 
$|1-\psi_n(\pi_n(g)^{-1}\pi_n(h))|<\epsilon$, so
for $g\in{}B(e,R)$, $|1-\phi(g)|<\epsilon$.
\end{itemize}
The existence of such a $\phi$ shows $G$ amenable. 

Conversely, say that $G$ is amenable.\\
Let $R,\epsilon>0$ and choose a positive type, finitely supported (say in $B(e,S)$ for some $S>R$) function
$\phi$ on $G$ as in the definition of amenability.\\
Again, write $F_n=G/K_n$.\\
As in the first half, there are only finitely many $n$ such that $\pi_n:G\to{}F_n$ is 
not an isometry on
$B(e,S)\subseteq{}G$; we can thus choose $N$ so that $\pi_n|_{B(e,S)}$ is an 
isometry for all $n\geq{}N$.\\
Define a kernel $k$ on $\Box{}G$ by setting:
\begin{displaymath}
k(x,y)=\left\{\begin{array}{ll}
1, & x\in{}F_n, y\in{}F_m, n,m<N \\
\phi(g^{-1}h), & x,y\in{}F_n, n\geq{}N, \pi_n(g)=x, \pi_n(h)=y, d(g,h)\leq{}S \\
0, & \textrm{otherwise.} 
\end{array} \right.
\end{displaymath}
Call the three possibilities (i), (ii), (iii).
Note that the assumption on $N$ implies that if $g_1,h_1$ and $g_2,h_2$ both satisfy
(ii), then $g_1^{-1}h_1,g_2^{-1}h_2\in{}B(e,S)$; as
$\pi_n(g_1^{-1}h_1)=x^{-1}y=\pi_n(g_2^{-1}h_2)$ and $\pi_n|_{B(e,S)}$ is an
isometry, $g_1^{-1}h_1=g_2^{-1}h_2$, and $k$ is thus well-defined.\\
$k$ is normalised as $\phi$ is.\\
It is of positive type as if $x_1,...,x_n\in\Box{}G$,
$\lambda_1,...,\lambda_n\in\mathbb{C}$,
\begin{align*}
\sum_{i,j}\lambda_i\overline{\lambda_j}k(x_i,x_j)
& =\sum_{(x_i,x_j) ~\textrm{satisfies (i)}}\lambda_i\overline{\lambda_j}
+\sum_{(x_i,x_j) ~\textrm{satisfies (ii)}}\lambda_i\overline{\lambda_j}
\phi(g_i^{-1}g_j) \\
& =\left|\sum_{x_i\in{}F_n, n<N}\lambda_i\right|^2
+\sum_{(x_i,x_j) ~\textrm{satisfies (ii)}}\lambda_i\overline{\lambda_j}\phi(g_i^{-1}g_j)
\geq{}0,
\end{align*}
as $\phi$ is of positive type.\\
$k$ is of finite propagation, as it is non-zero only on the finite set 
$\cup_{n<N}F_n\times\cup_{n<N}F_n$, and elsewhere on elements at most $S$ apart.\\
$k$ is self-adjoint by the argument in \ref{amenimpa}.\\
Note finally that for $n>N$, the assumption on $N$ implies that $F_n$ has diameter
at least $S$, whence $d(F_n,F_m)>S$ for all $m$ by assumption on the metric on
$\Box{}G$.  It follows that if $d(x,y)<R<S$, then $(x,y)$ satisfies either
(i) or (ii); in either case $|1-k(x,y)|<\epsilon$ by assumption on $\phi$, i.e
$k$ has $(R,\epsilon)$ variation.\\
It follows now that $\hat{k}=$Re$k$ has the properties in \ref{equiv}, (8).
\end{proof}

We also have the following related result, which gives more justification
to the slogan `A is to coarse embeddable as amenable is to a-T-menable'.

\begin{boxatmen} \label{boxatmen}
If $G$ is a residually finite countable discrete group and $\Box{}G$ (taken with respect to
any family of finite index normal subgroups $(K_n)$ as in \ref{boax}) is
coarsely embeddable, then $G$ is a-T-menable.
\end{boxatmen}

\begin{proof}
We skip some details in this proof, as the last one was done quite thoroughly, 
and this is similar.\\
Assume that $\Box{}G$ is coarsely embeddable, so by \ref{eh2}, (4), there 
exists a negative kernel $k$ on $\Box{}G$ and non decreasing
functions $\rho_1,\rho_2:\mathbb{R}^+\to\mathbb{R}^+$ such that
$\rho_1(t)\to+\infty$ as $t\to+\infty$ and
\begin{displaymath}
\rho_1(d(x,y))\leq{}k(x,y)\leq\rho_2(d(x,y))
\end{displaymath}
for all $x,y\in\Box{}G$.
By using a similar `averaging-and-limiting' argument to that in the above proposition to get a negative type
function $\psi$ on $G$ that satisfies
\begin{displaymath}
\rho_1(d(g,h))\leq{}\psi(g^{-1}h)\leq\rho_2(d(g,h))
\end{displaymath}
for all $g,h\in{}G$.\\
We can now use an analogue of Schoenberg's lemma (\ref{expneg}) to
get that $\phi_t:g\mapsto{}e^{-t\psi(g)}$ defines a positive type
function on $G$ for all $t>0$. For any $(R,\epsilon)>0$, then,
taking $t$ suitably large now gives us a positive type function on $G$ with
the properties in the definition of a-T-menability.
\end{proof}

Unlike the previous result, the converse to this is false.\\
This follows as if $H$ is a quotient of $G$, $\Box{}H$ coarsely embeds in
(is actually a subspace of, with suitably well chosen metrics) $\Box{}G$.
As the free group on $n$ generators is a-T-menable for all $n$, the converse would imply that
any finitely generated group coarsely embeds, which is false due to the
result discussed in \ref{nonagp}.

The idea (due to E. Guentner and J. Roe) behind both of these results is that if the box 
space constructed out of a group has property A (coarsely embeds), then we can
`average out' over it to show that the original group is
amenable (a-T-menable).  Our next construction, that of warped cones, also due to Roe, uses a 
similar idea to get an analogous result.

\begin{transl}
A map $\tau:X\to{}X$ from a metric space to itself is a \emph{translation} if
it is a bijection and there exists $S$ such that $d(x,\tau(x))<R$ for all
$x\in{}X$.
\end{transl}

Now, let $G$ act on a space $X$.  In general the maps $x\mapsto{}gx$ will not
be translations; the basic idea of warping is to alter the metric on
$X$ so that they are, i.e. so that the distance from $x$ to $gx$ is
uniformly bounded in $x$ for each $g\in{}G$.

\begin{warpm} \label{warpm}
Let $(X,d)$ be a proper metric space equipped with an action by homeomorphisms of
a countable discrete group $G$.\\
Assume $G$ is equipped with a bounded geometry left invariant metric, and associated
length function $|.|$.  We assume throughout this section that $|.|$ is
integer valued; the proof of \ref{wmce} shows there is no loss of generality
in doing this.\\
The \emph{warped metric on $X$ (with respect to the action of $G$)} is
defined to be the largest metric $d_G$ on $X$ satisfying
\begin{displaymath}
d_G(x,y)\leq{}d(x,y) ~\textrm{and} ~d_G(x,gx)\leq{}|g|
\end{displaymath}
for all $x,y\in{}X$ and all $g\in{}G$.
\end{warpm}

It is easy to show that $d_G(x,y)$ exists for all $x,y$ and defines a
metric on $X$, so is well-defined (cf. the proof of \ref{metrictree}).  
Clearly, $d_G$ captures the aim
of making $G$ act by translations. Moreover, the distance from $x$ to
$gx$ depends on the size of $g$ in $G$, so the group structure does have some
effect on the space $(X,d_G)$.  We will need the following more concrete
characterisation of $d_G$:

\begin{wmet}[\cite{roe2}, 1.6] \label{wmet}
Let $G$ act on $X$ as above.  Let $|.|$ be the length function on $G$ 
associated to the generating set $S$.  Then $d_G(x,y)$ is the
infimum over all finite sums
\begin{displaymath}
\sum_{i=1}^n(d(g_ix_i,x_{i+1})+|g_i|)
\end{displaymath}
where $x_i\in{}X$ for $i=1,,,n+1$, $g_i\in{}G$ for $i=1,...n$, 
$x_1=x$ and $x_{n+1}=y$.\\
Moreover, if $d_G(x,y)\leq{}k$, and we assume that $|g|\geq{}1$ for all $g\not=e$
then there exists a sequence as above such that $n=k$ and the infimum is attained.
\end{wmet}

\begin{proof}
The above does define a metric, say $\delta$, and it satisfies the conditions in 
\ref{warpm}, whence $\delta(x,y)\leq{}d_G(x,y)$ for all $x,y\in{}X$.\\
Conversely, $d_G\leq{}\delta$ by assumptions on $d_G$ and the triangle inequality.

Looking now at the second part, note that if $d_G(x,y)\leq{}k$, then there must
exist sequences of the above type where at most $k$ $g_i$ are non-identity elements;
as removing any terms where $g_i=e$ can only decrease the size of the sum, we can 
assume that any such sequence has at most $k$ terms.\\
For fixed endpoints $x,y$, the space of such chains can be assumed to be a closed 
subspace of $B_X(x,d(x,y)+1)^k\times{}B_G(e,k)^k$, which is compact by
properness of $X$, $G$.\\
Hence the infimum is attained.  
\end{proof}

The following proposition is the final preliminary we need before proving the main
results of this section.  It implies in particular that it makes sense to talk
about property A as defined in \ref{nonda} for $(X,d_G)$.

\begin{wcbgp}
$(X,d_G)$ is proper.\\
Moreover, if $X$ has bounded geometry and $G$ acts on $X$ by coarse maps, 
$(X,d_G)$ is of bounded geometry.
\end{wcbgp}

\begin{proof}
Let $x\in{}X$ and consider $\bar{B}_{d_G}(x,R)$ for any $R>0$.\\
It is $d$-closed, as any $d_G$ convergent sequence is $d$ convergent (as $|.|$ is
assumed integer valued).  Moreover, by the proposition above it is contained in
\begin{displaymath}
\bigcup_{g\in{}\bar{B}_G(e,R)}g(\bar{B}_d(x,R)),
\end{displaymath}
a compact set as $G$ acts by homeomorphisms and $(X,d)$ is assumed proper.\\
Hence $\bar{B}_{d_G}(x,R)$ is compact, so $d_G$ is proper.

For the second half, by bounded geometry there exists $\epsilon\geq{}0$ such that for
all $r>0$, all $x\in{}X$, $\bar{B}(x,r)$ can be covered by at most $N_r$ balls of radius
$\epsilon$.\\
Now, as each $g$ is acting as a coarse map, and $\bar{B}_G(e,r)$ is finite, the quantity 
\begin{displaymath}
R(r):=\sup_{g\in\bar{B}(e,r)}\{d(g(x),g(y)):d(x,y)<r\}
\end{displaymath}
exists.  Hence for any $x\in{}X$ we have that
\begin{displaymath}
\bigcup_{g\in\bar{B}_G(e,r)}g(\bar{B}_d(x,r))\subseteq{}\bigcup_{g\in\bar{B}_G(e,r)}
\bar{B}_d(g(x),R(r)).
\end{displaymath}
The right hand side can be covered by at most $N_{R(r)}|\bar{B}_G(e,r)|$ $d$-balls of
radius $\epsilon$, whence at most $N_{R(r)}|\bar{B}_G(e,r)|$ $d_G$-balls of radius
$\epsilon$ as the latter are at least as large as the former.  As
\begin{displaymath}
\bar{B}_{d_G}(x,r)\subseteq{}\bigcup_{g\in\bar{B}_G(e,r)}g(\bar{B}_d(x,r)),
\end{displaymath} 
$\bar{B}_{d_G}(x,r)$ can be covered by at most $N_{R(r)}|\bar{B}_G(e,r)|$ $d_G$-balls 
of radius $\epsilon$; this number is independent of $x$, so this completes the result.
\end{proof}

We have now established enough background to prove our first main result connecting
warped spaces, amenability and property A.

\begin{wamenimpa}\label{wamenimpa}
If $(X,d)$ is a bounded geometry proper metric space with property A and $G$ is a
countable discrete amenable group acting on $X$ by coarse maps, then
$(X,d_G)$ also has property A.
\end{wamenimpa}

\begin{proof}
Let $R,\epsilon>0$, and let $N>\max\{1/\epsilon,R\}$.\\
First, use definition \ref{amendef} to choose a finitely supported
map $\xi\in{}l^1(G)$ such that $\|\xi\|_1=1$ and
\begin{displaymath}
\|g\xi-\xi\|<\frac{1}{2N^2}~\textrm{whenever}~g\in\bar{B}_G(e,N).
\end{displaymath}
We use the right action of $G$ here, however, so set $g\xi(h)=\xi(hg^{-1})$; this makes
no difference to the definition of amenability.\\
Let $R'=\sup\{d(gx,gy):d(x,y)<N,g\in$Supp$(\xi)\}$, which exists as $\xi$ is assumed
to have finite support, and each element of $G$ acts as a coarse map.\\
Use property A for $(X,d)$, then, to get a map $\mu:X\to(C_0(X))^*$ such that $\|\mu_x\|=1$ for all
$x\in{}X$, there exists $S>0$ such that $\mu_x$ is supported in $\bar{B}_d(x,S)$ for all 
$x\in{}X$, and
\begin{displaymath}
\|\mu_x-\mu_y\|<\frac{1}{2N^2}~\textrm{whenever}~d(x,y)\leq{}R'
\end{displaymath}
We can now define a map $\nu:X\to(C_0(X))^*$ by setting
\begin{displaymath}
\nu_x=\sum_{g\in{}G}\xi(g)\mu_{gx}.
\end{displaymath}
Note first that $\nu_x$ has norm one for all $x\in{}X$.\\
Moreover, by definition of $d_G$, if $S'=\max\{|g|:\xi(g)\not=0\}$, which exists by finite
support, we have that $\nu_x$ is supported in $\bar{B}_{d_G}(x,S+S')$ for all $x\in{}X$.\\
It remains to show that $\nu$ has $(R,\epsilon)$ variation; we start by proving
the following:
\begin{itemize}
\item Say $d(x,y)<R\leq{}N$.  Then
\begin{displaymath}
\|\nu_x-\nu_y\|=\left\|\sum_{g\in{}G}\xi(g)\mu_{gx}-\sum_{g\in{}G}\xi(g)\mu_{gy}\right\|
\leq\sum_{g\in{}G}|\xi(g)|\|\mu_{gx}-\mu_{gy}\|
\end{displaymath}
By choice of $R'$ and $\mu$, and as $\|\xi\|=1$, this is bounded by $1/2N^2$.
\item Say $|g_0|<N\leq{}R$ and $y=gx$.  Then
\begin{align*}
\|\nu_x-\nu_y\|
& =\left\|\sum_{g\in{}G}\xi(g)\mu_{gx}-\sum_{g\in{}G}\xi(g)\mu_{gg_0x}\right\|\\
& =\left\|\sum_{g\in{}G}(\xi(g)-\xi(gg_0^{-1}))\mu_{gx}\right\| \\
& =\sum_{g\in{}G}|\xi(g)-\xi(gg_0^{-1})|\|\mu_{gx}\|
=\|g_0\xi-\xi\|_1.
\end{align*}
By choice of $\xi$, this is also bounded by $1/2N^2$.
\end{itemize}
Finally, we note that by the second part of \ref{wmet}, if $d_G(x,y)<R\leq{}N$, then
there exist sequences $x_1,...,x_{N+1}$ and $g_1,...,g_N$ such that 
$x=x_1$, $y=x_{N+1}$ and
\begin{displaymath}
d_G(x,y)=\sum_{i=1}^N(d(g_ix_i,x_{i+1})+|g_i|).
\end{displaymath}
By the two comments above, this is at most $2N/2N^2=1/N<\epsilon$. 
\end{proof}

In order to get a partial converse to this, we move now to the special case of
warped cones. 

\begin{cone}
Let $Y$ be a compact metric space.\\
Embed $Y$ into the unit sphere of some (separable) Hilbert space (this is always possible).\\
The \emph{open cone} over $Y$, $\mathcal{O}Y$, is the union of all rays through the 
origin and some point of the embedded copy of $Y$.  It is equipped with the induced metric
from the Hilbert space.
\end{cone}

Note that we can parametrise points on $\mathcal{O}Y$ as pairs $(y,t)$,
where $t\in\mathbb{R}^+$, and that this is unique away from $t=0$.  Note also that while all compact metric spaces are coarsely equivalent, they can
give rise to coarsely distinct cones.

\begin{warpc}
Let $G$ be a countable discrete group acting on a
compact metric space $Y$.  It acts on $\mathcal{O}Y$ by 
$g:(y,t)\to(gy,t)$.\\
The \emph{warped cone of $Y$}, $\mathcal{O}_GY$ is the space 
$\mathcal{O}Y$ equipped with the metric $d_G$ arising from this action.
\end{warpc}

\begin{wcamenimpa}
If $Y$ is a compact metric space, $G$ is a countable discrete amenable group acting on 
$Y$ by Lipschitz homeomorphisms, and $\mathcal{O}Y$ has property A, then $\mathcal{O}_GY$ has property A.
\end{wcamenimpa}

As an example, note that if $Y$ embeds into the unit sphere of some finite dimensional Hilbert space (e.g. if $Y$ is a compact
manifold or finite simplicial complex), then $\mathcal{O}Y$ is a subset of property A $R^\mathbb{n}$, whence has property A itself.

\begin{proof}
If $G$ acts on $Y$ by Lipschitz homeomorphisms, it also acts on $\mathcal{O}Y$ by
Lipschitz homeomorphisms, which are coarse maps.\\
Moreover, $\mathcal{O}Y$ is proper, of bounded
geometry and has property A.  We can now apply \ref{wamenimpa}.
\end{proof}

We get the following partial converse to the previous result.

\begin{wcaimpamen} \label{wcaimpamen}
Let $Y$ be a compact metric space, and $G$ a countable, 
subgroup.  Assume that $Y$ is equipped with a (finite) $G$-invariant Borel measure.
If $\mathcal{O}_GY$ (constructed with respect to the left
translation action of $G$) has property A, then $G$ is amenable.
\end{wcaimpamen}

As an example, say that $Y$ is a compact Lie group (equipped with normalised Haar measure), and that $G$ any 
countable subgroup. 

\begin{proof}
By property A as formulated in
\ref{nonda2} for any $R,\epsilon>0$ there exists a continuous, positive type,
finite propagation normalised symmetric kernel $k$ on $\mathcal{O}_GY$ of
$(R,\epsilon)$ variation.\\
We now use basically the same argument as in \ref{boxamen} to average out
$k$ over the entire warped cone and get a positive type function showing
$G$ amenable.\\
Let $\mu$ be a (normalised) $G$-invariant measure on $Y$, parametrise points on
$\mathcal{O}_GY$ as pairs $(y,t)$, and average to get
\begin{displaymath}
\hat{k}((y_1,t_1),(y_2,t_2))=\int_Yk((y_1y,t_1),(y_2y,t_2))d\mu(y).
\end{displaymath}
$\hat{k}$ is thus right invariant under multiplication by elements of $G$.  It inherits the properties of being
positive type, finite propagation, normalised, symmetric and of $(R,\epsilon)$ variation from $k$.
We define for all $t\in\mathbb{R}^+$, all $g\in{}G$,
\begin{displaymath}
\phi_t(g)=\hat{k}((y,t),(gy,t)),
\end{displaymath}
which is well defined by right invariance of $\hat{k}$.\\
$\phi_t$ is then a normalised, positive type kernel on $G$ for all $t$.\\
It follows that $|\phi_t(g)|\leq{}1$ for all $t,g$, whence we can use a diagonal
argument (and countability of $G$) to find an increasing sequence $t_n\to\infty$ of real numbers
such that $\phi_{t_n}$ converges pointwise to some $\phi$, which is thus a
normalised, positive type function on $G$.\\
We finish the proof with the following observations:
\begin{itemize}
\item For any $g\in{}G$ and all $t$ suitably large, $d_G((y,t),(gy,t))=|g|$; as
$\hat{k}$ is of finite propagation $S$, say, if $|g|>S$,
$\hat{k}((y,t),(gy,t))=0$ for all $t$ suitably large.  Hence $\phi(g)=0$.
\item If $d_G(g,h)<R$ for some $g,h\in{}G$, $\hat{k}$ of $(R,\epsilon)$ variation
forces $|1-\hat{k}((e,t),(g^{-1}h,t))|<\epsilon$ for all $t$, as
$d_G(e,g^{-1}h)=d_G(g,h)<R$.\\
Hence $|1-\phi(g^{-1}h)|<\epsilon$.
\end{itemize}
The existence of such a $\phi$ shows $G$ amenable.
\end{proof}

It is worth noting (in analogy with \ref{boxatmen})
that a similar proof shows that the
same hypotheses, with `coarse embeddability' replacing `property A', imply
that $G$ is a-T-menable \cite{roe2}, 4.4.  The proof is essentially the same.

\pagebreak

\section {Spaces without property A}
\label{nona}

\subsection{Introduction} \label{nonaintro}
One of the most interesting and active areas of research centred around property
A is the attempt to find non-A spaces, motivated partly by the
connections to the coarse Baum-Connes and Novikov conjectures.  We restrict to the (harder) bounded geometry case, as this class of
spaces provides much of our motivation
(chapters 8 and 9 of \cite{bl} provide some non-bounded geometry examples).
The exception is \ref{nonace}, where
an explicit example of a locally finite space that coarsely embeds but does not have property
A is interesting enough to merit its inclusion.

This chapter is split up into four sections, each of which discusses a
particular class of non-property A spaces: expanding graphs; box spaces and warped cones arising from property T groups;
spaces without A that coarsely embed; groups without property A.

\subsection{Expanding graphs} \label{expander}
This section concerns infinite families of graphs, that in some sense `expand
too quickly' to coarsely embed into a Hilbert space, and so
cannot have property A.  They
are perhaps our simplest example of non-property A spaces.

The crucial property
is the size of the eigenvalues of the Laplacian (almost equivalently, the
adjacency matrix) of a finite graph.  For a general discussion of this
topic, including the relationship of the graph Laplacian to the
Laplace-Beltrami operator on a manifold, see \cite{lub}, chapter 4.

For simplicity, we will restrict our attention to 
undirected graphs $\mathcal{G}$ of constant degree $D$, by which we mean that there are
exactly $D$ edges incident at each vertex.  We will also
say that all edges are between distinct vertices, and that any two vertices have at most
one edge between them.

Write $\mathcal{G}=(V,E)$, where
$V$ and $E$ are the vertex and edge sets respectively of $\mathcal{G}$; 
if $\mathcal{G}$ has
$n$ vertices, we shall assume that they are indexed as $v_1,...,v_n$.
Edges will be written as pairs $(v_i,v_j)$ of endpoints.  Note that $(v_i,v_j)=(v_j,v_i)$, 
as we assume edges are undirected.
We define a metric by setting $d(v,w)$ to be the smallest number of edges forming a path from
$v$ to $w$.

The following defines the adjacency matrix, and Laplacian.

\begin{adjlap}
Let $\mathcal{G}$ be a finite graph of constant degree $D$.\\
The \emph{adjacency matrix} of $\mathcal{G}$ is the $n\times{}n$ matrix $A=(a_{i,j})$ 
such that $a_{i,j}=1$ if there is an edge between $v_i$ and $v_j$, and zero 
otherwise.\\
The \emph{Laplacian} of $\mathcal{G}$ is $\Delta=DI-A$, where $I$ is the $n\times{}n$ identity
matrix.
\end{adjlap}

Both of these matrices are symmetric, and in fact positive (we prove this for
the Laplacian as part of \ref{laplac} below).  The largest
eigenvalue of $A$ is always $D$, the smallest of $\Delta$, $0$, which have eigenspace
the constant vectors for $D$-regular connected graphs.
We shall be interested in the second largest eigenvalue of $A$; equivalently,
the smallest positive eigenvalue of $\Delta$, which we denote by $\lambda(\mathcal{G})$.

The size of this eigenvalue is connected to `connectivity' or expansion 
properties of the graph, such as the number of neighbours of a set, the number
of edges across a cut, or the speed of convergence of the random walk on $\mathcal{G}$.
These properties lead to applications of expanders in computer science and
combinatorics.

We shall concentrate on the Laplacian, starting by proving the following:

\begin{laplac} \label{laplac}
Let $\mathcal{G}=(V,E)$ be a finite connected D-regular graph, with Laplacian $\Delta$, let
$\lambda=\lambda(\mathcal{G})$ and let
$f:V\to\mathbb{R}$.  Then 
\begin{displaymath}
\sum_{v\in{}V}(f(v)-M)^2\leq\frac{1}{\lambda}\sum_{(v,w)\in{}E}(f(v)-f(w))^2
\end{displaymath}
where $M=\frac{1}{|V|}\sum_{v\in{}V}f(v)$.
\end{laplac}

\begin{proof}
Note first that by definition of the Laplacian,
\begin{displaymath}
(\Delta{}f)(v)=Df(v)-\sum_{(v,w)\in{}E}f(w)=\sum_{(v,w)\in{}E}(f(v)-f(w)),
\end{displaymath}
whence
\begin{displaymath}
\langle\Delta{}f,f\rangle=\sum_{v\in{}V}f(v)\left(\sum_{(v,w)\in{}E}f(v)-f(w)\right)
=\sum_{(v,w)\in{}E}(f(v)-f(w))^2,
\end{displaymath}
so $\Delta$ is positive (on $L^2(V,\mathbb{R})$), and the eigenvalue $0$ has one-dimensional eigenspace consisting
of the constant functions.\\
Now as $\Delta$ is symmetric and $f-M$ is orthogonal to $Ker(\Delta)$, it follows that
\begin{displaymath}
\langle\Delta{}(f-M),f-M\rangle\geq\lambda\|f-M\|^2.
\end{displaymath} 
Expanding both sides, and using that $\Delta{}M$=0=$\langle\Delta{}f,M\rangle$:
\begin{displaymath}
\sum_{v\in{}V}(f(v)-M)^2\leq\langle\Delta{}f,f\rangle=\sum_{(v,w)\in{}E}(f(v)-f(w))^2.
\end{displaymath}
\end{proof}

This says that the total variance of $f$ from its mean value is bounded by
$1/\lambda$ multiplied by the sum of variances over each edge.  If $\lambda$ is
large, then, the variance over the whole graph must be small compared to the
sum of the variances over each edge, which suggests good connectivity.

Having introduced the necessary ideas, we now make the main definition.

\begin{expander}
An infinite family $(\mathcal{G}_n)_{n=1}^{\infty}=(V_n,E_n)_{n=1}^{\infty}$ of finite $D$-regular graphs
and $\lambda_n:=\lambda(\mathcal{G}_n)$ is called an \emph{expander} if:
\begin{itemize}
\item the size of the graphs, $|V_n|$, tends to infinity as $n\to\infty$;
\item there exists $\lambda>0$ such that $\lambda_n>\lambda$ for all $n$.
\end{itemize}
\end{expander}

We equip the disjoint union of the vertex sets of such an infinite family with a metric that
restricts to the metric defined above on each $V_n$, and keeps distinct
$V_n$ at least as far apart as their diameters.  Any two such metrics are coarsely equivalent.

For intuition, and as we will use it in the next section, we state the following 
equivalent definition, which has perhaps a more obvious connection to 
the `connectivity' properties mentioned above.  A proof that the two definitions are
equivalent can be found in \cite{lub} 4.2.4-5 and 1.1.4.

\begin{exp2} \label{exp2}
An infinite family of graphs as above (but that does not
necessarily have bounded Laplacian eigenvalues) is an \emph{expander} if
and only if there exists $c>0$ such that for each $n$ and any subset $A$ of $V_n$
\begin{displaymath}
|\partial{}A|\geq{}c\left(1-\frac{|A|}{|V_n|}\right)|A|,
\end{displaymath}
where $\partial{}A$ (the \emph{boundary} of $A$) is the collection of elements of 
$V_n$ that are connected to an element of $A$ (but not in $A$).
\end{exp2}

This formalises the idea that `any set of vertices has a large number of neighbours'
mentioned above.

The graphs in the family thus have good connectivity properties in 
a `uniform' fashion.  The proof of the next result is due to N. Higson; the idea is 
that the graphs in an expanding family are so `well-connected' (i.e. vertices
have a lot of other vertices close to them) that for any map from the
expander to $l^2$ there exists a sequence of 
$l^2$ balls of fixed radius having arbitrarily large numbers of vertices mapped into them.
As the subgraphs of an expander all have the same
constant degree, such a map cannot be a coarse embedding.

\begin{expembed}
No expander coarsely embeds into Hilbert space.
\end{expembed}

\begin{proof}
Consider any expander $\mathcal{E}=(V_n,E_n)_{n=1}^{\infty}$ of
constant degree $D$ and
lower bound on the first positive eigenvalue of its Laplacian $\lambda$.\\
Say for contradiction that $f:\mathcal{E}\to\mathcal{H}$ be a coarse embedding.\\
As $f$ is effectively proper, for each $r>0$
there exists $M_r$ such that no more than $M_r$ vertices are mapped into
any ball of radius $r$.\\
As $f$ is bornologous, there exists $c>0$ such that any edge $(v_i,v_j)$ in 
any subgraph of $\mathcal{E}$, $d(f(v_i),f(v_j))<c$.\\
Let $f_n$ be the restriction of $f$ to $V_n$, and by
translating the origin, assume that $\sum_{v\in{}V_n}f_n(v)=0$.\\
Then by summing over all coordinates, the result of \ref{laplac} give that
\begin{displaymath}
\sum_{v\in{}V_n}\|f_n(v)\|^2
\leq\frac{1}{\lambda}\sum_{(v,w)\in{}E_n}\|f_n(v)-f_n(w)\|^2,
\end{displaymath}
as the reciprocals of the first positive eigenvalues
are bounded above by $1/\lambda$.
Now, the right hand side is at most $|V_n|c^2/\lambda$, whence no more 
than $|V_n|/2$ terms on the left hand side can have absolute 
value more than $2{}c^2/\lambda$.\\
This says that $f_n$ maps at least $|V_n|/2$ points into 
$B_{l^2}(0,2c^2/\lambda)$, whence $M_{2c^2/\lambda}$ is unbounded.  This is
a contradiction, so no such $f$ can exist.
\end{proof}

It follows, of course, that no expander can have property A.  Our next question, 
then, is: do expanders actually exist?  The answer is yes; there are several
ways to prove this.

Firstly, probabilistic methods can be used to show that for large $n$, `most'
$D$-regular graphs are expanders (for some fixed bound on the size of the
smallest positive eigenvalue of the Laplacian).  This provides the easiest proof;
concretely constructing an expanding family is somewhat more difficult.

One way to do this is given in the paper \cite{rvw}.  Here an operation on
graphs called the `zig-zag product' is described; it takes a large graph $\mathcal{G}_1$
and a small graph $\mathcal{G}_2$ and produces a larger graph $\mathcal{G}$ that inherits its 
degree from $\mathcal{G}_2$ and its expansion properties from both graphs.  For a
sufficiently well-chosen starting point, and with a few additional (basic)
techniques, and infinite expanding family can be (recursively) constructed.
This construction has the advantage of using only basic graph theory,
unlike other known constructions (\cite{mar} and \cite{lps}), which make fairly 
heavy use of algebraic
techniques.

We will, however, use a construction based on that in \cite{mar} due to Margulis, as it uses some group theory that we are interested in for other reasons
(this is the content of the next section).  
The construction is given in \ref{proptexp}.

\subsection{Property T} \label{propT}
Property T (or sometimes `the Kazhdan property', or `Kazhdan's property T') is
an attribute of locally compact (so in particular, discrete) groups that is in some sense `opposite' to 
amenability.  Specifically, a discrete group that is both amenable and has property
T must be finite (finite groups trivially have both properties).  
In fact, this is also true of groups that are both a-T-menable and property T groups.  The name `a-T-menable' 
comes from its being a weak form of amenability that is still strong
enough to be `anti-T'.  For an introduction (and more) to property T see \cite{bhv}.

We restrict (as usual) to the case of discrete groups.

Prior to the definition we recall that a \emph{(unitary) representation} of a discrete
group $G$ consists of a pair $(\pi,\mathcal{H})$ where $\mathcal{H}$ is a Hilbert space,
and $\pi:G\to\mathcal{B}(\mathcal{H})$ is a homomorphism such that $\pi(g)$ is unitary
for all $g\in{}G$.\\
It is \emph{reducible} if there exists a proper closed subspace
$\mathcal{H'}\subset\mathcal{H}$ such that $\pi(g)(\mathcal{H'})\subseteq{}\mathcal{H'}$ for
all $g\in{}G$.  If this occurs, we say that $\pi':G\to\mathcal{B}(\mathcal{H}')$,
$g\mapsto\pi(g)|_{\mathcal{H'}}$ is a
\emph{subrepresentation} of $\pi$.\\
If no such $\mathcal{H'}$ exists, $\pi$ is \emph{irreducible}. 

\begin{propT}
A finitely generated discrete group $G$ with some finite generating set $S$ has \emph{property T} if there exists $\epsilon>0$
such that for all non-trivial irreducible
representations $(\pi,\mathcal{H})$
of $G$ there exists $s\in{}S$ such that $\|\pi(s)v-v\|_{\mathcal{H}}\geq\epsilon$ for all unit vectors $v\in\mathcal{H}$.
\end{propT}

The idea is that a representation of a property T group that does not contain the
trivial representation  $\mathbf{1}$ (i.e. has no non-trivial fixed vectors) cannot get `close' to the trivial representation. 
(there is sense in which this can be made precise using the \emph{Fell topology} on the unitary dual of the group).
This is a rigidity property that restricts possible $G$-actions.

Note that, contrary to our usual habits, we have demanded that property T groups be finitely generated.  
This is no real restriction: any discrete group with property T has to be finitely generated (see \cite{bhv}, 1.3.1).  Assuming finite
generation from the start makes the definition slightly easier to state, however, so we did so.

We basically know the next lemma already; it connects unitary representations of a group
with positive definite functions on it.

\begin{replem} \label{replem}
A positive type function $\phi$ on a discrete group $G$ defines a
unitary representation $(\pi_{\phi},\mathcal{H}_{\phi})$.
\end{replem} 

\begin{proof}
The Hilbert space $\mathcal{H}_{\phi}$ is that arising from the positive type kernel
$k(g,h)=\phi(g^{-1}h)$ and theorem \ref{kerthe}.\\
On the finitely supported functions in $\mathcal{H}_{\phi}$, $\pi_{\phi}(g)$ is defined
by $(\pi_{\phi}(g)f)(h)=f(g^{-1}h)$.\\
Note then that for two such functions $f,f'$, and by construction of $\mathcal{H}_{\phi}$,
\begin{align*}
\langle\pi_{\phi}(g)f,\pi_{\phi}(g)f'\rangle
& =\sum_{g_1,g_2\in{}G}f(g^{-1}g_1)\overline{f(g^{-1}g_2)}\phi(g_1^{-1}g_2) \\
& =\sum_{g_1,g_2\in{}G}f(g_1)\overline{f(g_2)}\phi(g_1^{-1}gg^{-1}g_2)
=\langle{}f,f'\rangle.
\end{align*}
It follows that $\pi_{\phi}(g)$ defines a unitary operator in 
$\mathcal{B}(\mathcal{H}_{\phi})$ as required.
\end{proof}

This leads to the following idea: a-T-menable groups have positive type functions that get close
to to (but not uniformly close, as long as the
group is infinite) the trivial function $\mathbf{1}$.  This can be used to show that they have non-trivial irreducible
representations that are `close to being trivial', so cannot have property T.  The following makes this
precise:

\begin{atment}
An infinite a-T-menable group $G$ does not have property T.
\end{atment}

\begin{proof}
As property T groups must be finitely generated, we may assume that there exists a finite generating set $S\subset{}G$.  

Let now $\epsilon>0$. 
We will first show that there exists a representation $(\pi,\mathcal{H})$ and a 
unit vector $v\in\mathcal{H}$ such that $\|\pi(s)v-v\|<\epsilon$ for all
$s\in{}S$.

Equip $G$ with the word metric coming form $S$ (so all elements of $S$ have length one), 
and choose a normalised positive type function
$\phi$ as in the definition of a-T-menability (\ref{atmen}) 
for parameters $R=1,\epsilon^2/2$.\\
Let $(\pi,\mathcal{H})$ be the representation associated to $\phi$ by \ref{replem}.\\
Note first that for all $g\in{}G$, there is an equivalence class of functions
corresponding to $\delta_g$ in $\mathcal{H}$ (which we also denote $\delta_g$) and moreover that
$\langle\delta_g,\delta_g\rangle=\phi(e)=1$.\\
Now, for any $s\in{}S$,
\begin{align*}
\|\pi(s)\delta_e-\delta_e\|^2 & =\|\delta_s-\delta_e\|^2
=2-2\textrm{Re}\langle\delta_s,\delta_e\rangle \\
& = 2-2\textrm{Re}\phi(s)<2\epsilon^2/2.
\end{align*}
Hence (for $v=\delta_e$) we have that $\|\pi(s)v-v\|<\epsilon$ for all $s\in{}S$.

On the other hand, note that for any $x,y\in{}G$,
\begin{displaymath} 
\langle{}\pi(g)\delta_x,\delta_y\rangle=\phi(x^{-1}g^{-1}y)\to{}0
\end{displaymath}
as $g\to\infty$ by definition of the positive type functions given by a-T-menability, and the fact that $g\to\infty$ forces $x^{-1}g^{-1}y\to\infty$.  It follows that if $f$ is
any (equivalence class of a) finitely supported function in $\mathcal{H}$, then $\langle{}\pi(g)f,f\rangle\to{}0$ as $g\to\infty$.  As such elements are dense in $\mathcal{H}$, 
$\langle{}\pi(g)w,w\rangle\to0$ as $g\to\infty$ for any $w\in\mathcal{H}$.  In particular, then, $(\mathcal{H},\pi)$ has no fixed vectors.  The existence of such a representation now shows that $G$ cannot have property T. 
\end{proof}

From the results of section \ref{bswc}, we therefore get: 

\begin{itemize}
\item from \ref{boxatmen},if $G$ is a residually finite, countably infinite discrete group with property T, then 
$\Box{}G$ is not coarsely embeddable;
\item from \ref{wcaimpamen}, if $G$ has property T and is a countably infinite discrete subgroup of a
Lie group $Y$, then the warped cone $\mathcal{O}_GY$ does not have property A.
\end{itemize}

In the second case, this can be extended to show that $\mathcal{O}_GY$ 
does not coarsely embed, as is done in \cite{roe2}, 4.4; the proof uses the same ideas as those in \ref{wcaimpamen}.

Groups with these properties do exist in both cases:
in the box space case, a good example is $SL(n,\mathbb{Z})$ for $n\geq{}3$, while for
warped cones, $SO(5)$ contains subgroups of the required type; see e.g \cite{lub}, chapter 3 or
\cite{bhv}.  This gives us more examples of non-property A spaces.

They are not all new, however: box spaces of property T
groups are actually expanders.  This is made precise by the following result (which also provides a proof that expanders exist,
at least assuming that property T groups do).  

\begin{proptexp} \label{proptexp}
Say $G$ is a (finitely generated), infinite, residually finite property T group.  Equip it with the word metric associated to some symmetric finite generating set $S$ (not containing the identity).
Say $(K_n)_{n=1}^{\infty}$ is a family of finite index subgroups as in \ref{boax}, subject to the further condition that $S^2\cap{}K_1=\{e\}$.\\
Then $\Box{}G$ built with respect to this family is an expander.
\end{proptexp}

The condition `$S^2\cap{}K_1=\{e\}$' is not really important.  It is included as we have demanded that
the graphs making up our expanders have no loops and only one edge between any given pair of vertices; this condition guarantees
that this will be the case.  In a similar vein, equipping $G$ with the word metric coming from $S$ is not terribly important: any two bounded geometry,
left-invariant metrics on a group give rise to coarsely equivalent box spaces (as long as they are taken with respect to the same chain of subgroups).  Hence
if one was interested only in producing spaces
coarsely equivalent to an expander, any metric on $G$ would do.  As we demand that our expanders are formed from graphs, however, 
we are forced to use a word metric.

The proof uses the definition of an expander from \ref{exp2}, and is
taken from \cite{lub}, 3.3.1.
  
\begin{proof}
Let $K$ be any finite-index, proper normal subgroup of $G$ (in particular, any of the $K_n$), and consider
$l^2(G/K)$.

Letting $\bar{g}$ denote the image of $g\in{}G$ in $G/K$, we get a
representation of $G$ on $l^2(G/K)$ by setting $\pi(g)f=\bar{g}f$ (cf.
\ref{amendef} for notation).\\
Further, let $E$ be the subspace of $l^2(G/K)$ of functions with mean value $0$, and
note that $l^2(G/K)=E\oplus\mathbb{C}\mathbf{1}$, the second summand consisting of
the constant functions.  The action of $G$ fixes $E$ (and $\mathbb{C}\mathbf{1}$),
so we get a subrepresentation $\pi'$ of $G$ on $E$.

Now, $\pi'$ has no fixed (non-trivial) subspace, as the action of $G$ on $G/K$ is transitive;
in particular $\pi'$ does not contain the trivial representation.\\
Hence by property T, there exists $\epsilon>0$ and
$s\in{}S$ with $\|\bar{s}f-f\|\geq\epsilon$ for all $f\in{}l^2(G/K)$.  Crucially, note that
$\epsilon>0$ is \emph{independent} $K$.

We claim now that $\Box{}G$ as in the statement is an
expander.  Note first that  $|G/K_n|\to\infty$ as $n\to\infty$ as $\cap{}K_n=\{1\}$.  They are also of constant 
valency $|S|$.\\
To show the `expansion' part of the definition, let $A$ be any subset of $G/K_n$, let $a=|A|$, $B=(G/K_n)\backslash{}A$, $b=|B|=|G/K_n|-a$.\\
Define $f\in{}E$, $E$ the subspace of $l^2(G/K_n)$ consisting of elements of mean 
zero as above by setting
\begin{displaymath}
f(x)=\left\{\begin{array}{ll}
b, & x\in{}A \\ -a, & x\in{}B.
\end{array}\right.
\end{displaymath}
Write $m=|G/K_n|$.\\
Then $\|f\|=ab^2+ba^2=mab$ and for any $s\in{}S$, we have that
\begin{displaymath}
\|\bar{s}f-f\|^2=
\sum_{x\in{}A,\bar{s}^{-1}x\in{}B}|-a-b|^2
+\sum_{x\in{}B,\bar{s}^{-1}x\in{}A}|b+a|^2=(a+b)^2E_s(A,B),
\end{displaymath}
where $E_s(A,B)$ is the number of `$s$-edges' between $A,B$, i.e.
\begin{displaymath}
E_s(A,B)=|\{x\in{}G/K_n:x\in{}A ~\textrm{and}~ \bar{s}^{-1}x\in{}B, ~\textrm{or}~ 
x\in{}A ~\textrm{and}~\bar{s}^{-1}x\in{}B\}|.
\end{displaymath}
Using the result above, then, there exists $\epsilon>0$ independent of $K_n$ and
$s\in{}S$ such that
\begin{displaymath}
|\partial{}A|\geq\frac{1}{2}E_s(A,B)=\frac{\|\bar{s}f-f\|^2}{2(a+b)^2}
\geq\frac{\epsilon^2\|f\|^2}{2m^2}=\epsilon^2\frac{ab}{2m}
=\frac{\epsilon^2}{2}.(1-\frac{a}{m})a.
\end{displaymath}
This is exactly the condition in \ref{exp2}.
\end{proof}

At time of writing (and to the best of the author's knowledge)
there is no known general group-theoretic property that is strong enough to deny 
the resulting box space or warped cone property A, but not so strong that it
denies coarse embeddability.

There are some interesting extensions in this direction.  The first is detailed in \cite{lub}, where a 
condition called \emph{property $\tau$} is formulated.  This gives a precise criterion for when the box space
corresponding to an infinite family of subgroups is an expander.  The second is due to V. Lafforgue, \cite{laf}.
Here the author formulates a stronger version of property T; in the last section it is shown that box spaces arising
from such a group cannot coarsely embed in \emph{any} uniformly convex space.

The next section gives a different construction, which uses a group (an amenable group, in
fact) to build a space that does coarsely embed, but does not have property A. 

\subsection{Spaces without property A that coarsely embed in Hilbert space}
\label{nonace}
There is only one known construction of spaces that coarsely embed, but do not have
property A, which is due to P. Nowak, \cite{no3}.  The spaces constructed are surprisingly
simple; showing that they cannot have property A, however, requires a careful
analysis of `how well' property A is satisfied for direct products of amenable groups. 
Nowak has expanded on
these ideas in \cite{no3}, where he relates `how well' a space has property A to amenability,
growth rates, isoperimetric profiles and finite asymptotic dimension.

Nowak's spaces are locally finite, but are not, however, bounded geometry.  The search for a
bounded geometry, coarsely embeddable space without property A (it is generally assumed one 
exists, but...) is an active area of current research.  

Throughout, we shall be interested in functions $\xi:X\to{}l^1(X)$ satisfying \ref{equiv},
condition (2) (for $p=1$), and similarly functions $f\in{}l^1(X)$ satisfying \ref{amendef}, (1).
We start with the following definitions, which are the first step in
quantifying `how well' a space has property A, or `how amenable' a group is.

We assume throughout this section (as 
\ref{wmce} entitles us to) that any left-invariant bounded geometry
metric on a countable discrete group is integer valued and has
at least three elements of norm at most one. 

\begin{quanta}
Let $X$ be a discrete space, and $G$ a countable discrete group.
\begin{enumerate}
\item We define $\textrm{diam}^{\mathcal{A}}_X(R,\epsilon)$
to be
\begin{displaymath}
\inf\{S:\textrm{supp}\xi_x\subseteq\bar{B}(e,S)
~\textrm{for all}~x\in{}X~\textrm{where}~ \xi~\textrm{satisfies \ref{equiv}, (2)}\},
\end{displaymath}
if this exists, and infinity otherwise.
\item We define
\begin{displaymath}
\textrm{diam}^{\mathcal{F}}_G(R,\epsilon)=\inf\{S:\textrm{supp}f\subseteq\bar{B}(e,S),
~\textrm{where}~ f~\textrm{satisfies \ref{amendef}}\},
\end{displaymath}
if this exists, and infinity otherwise.
\end{enumerate}
\end{quanta}

Note that in the case of a group, our assumption that the metric on $G$ is integer valued
implies that the infimums in the above are actually obtained by some functions 
$\xi$, $f$; we will use this.

The idea is that $\textrm{diam}^{\mathcal{A}}_X$ measures how small we can keep the
support of a collection of elements in $l^1(X)$ that show $X$ has property A.
Similarly $\textrm{diam}^{\mathcal{F}}_G$ measures how small the support of 
an element of $l^1(G)$ showing that $G$ is amenable can be.

The following proposition sums this up.  It is immediate once we recall that 
\ref{equiv} (2) is equivalent to
property A in the bounded geometry case, and implied by it otherwise.

\begin{qaprops}
Let $X$ be a discrete space, $G$ a countable discrete group.
\begin{enumerate}
\item If $X$ has property A then 
\emph{diam}$^{\mathcal{A}}_X(R,\epsilon)<\infty$ for all $R,\epsilon>0$.
\item If $X$ is of bounded geometry, then it has A if and only if
\emph{diam}$^{\mathcal{A}}_X(R,\epsilon)<\infty$ for all $R,\epsilon>0$.
\item $G$ is amenable if and only if \emph{diam}$^{\mathcal{F}}_G(R,\epsilon)<\infty$
for all $R,\epsilon>0$.
\end{enumerate}
\end{qaprops}

Note that the numerical values of $\textrm{diam}^{\mathcal{A}}_X$ and 
$\textrm{diam}^{\mathcal{F}}_G$ will depend on the metrics on $X$, $G$ respectively;
whether they are infinite or not will only depend on the coarse equivalence class
of metric.

The following theorem can be construed as saying
that an amenable group has amenability `as well' as it has property A; given that the
latter is quite a lot weaker than the former, this is somewhat surprising.  The basic
idea is to use amenability to average out some $\xi$ showing property A to get an $f$
such that the support of $f$ is bounded by that of $\xi_g$ for any $g$.

\begin{aeqamen}[\cite{no3}, 3.3] \label{aeqamen}
Let $G$ be a countable discrete amenable group, and fix $R\geq{}1$, $\epsilon>0$. Then
\begin{displaymath}
\textrm{\emph{diam}}^{\mathcal{A}}_G(R,\epsilon)=\textrm{\emph{diam}}^{\mathcal{F}}_G(R,\epsilon)
\end{displaymath}
\end{aeqamen}

\begin{proof}
To show that 
$\textrm{diam}^{\mathcal{A}}_G(R,\epsilon)\leq\textrm{diam}^{\mathcal{F}}_G(R,\epsilon)$ for
any given $R,\epsilon$, it suffices to note that given $f$ attaining the infimum on the
right hand side, $\xi$ defined by $\xi_g=gf$ satisfies the right properties, with the
same support bound as $f$.\\
Conversely, let $R,\epsilon>0$, and let $\xi:G\to{}l^1(G)$ attain the infimum on the
left hand side, say $S$.  As usual, by replacing $\xi_g$ by $h\mapsto|\xi_g(h)|$, we
may assume $\xi_g$ positive for all such $h$ and all $g\in{}G$.\\
Let $M(f)=\int_Gf(g)dg$ be the invariant mean on $l^{\infty}(G)$ provided by 
amenability (see \ref{amendef}), and define a function $f:G\to\mathbb{C}$ by
\begin{displaymath}
f(h)=\int_G\xi_g(gh)dg,
\end{displaymath}
which is well defined as $\xi_{.}(\cdot{}h):g\mapsto\xi_g(gh)$ is bounded by
one, so in $l^{\infty}$.\\
Note first that $f(h)=0$ for all $|h|>S$, as then $d(g,gh)=|h|>S$ for all
$g\in{}G$, whence $\xi_g(gh)=0$ for all $g\in{}G$.\\
In particular, then
\begin{displaymath}
\|f\|_1=\sum_{h\in{}\bar{B}(e,S)}\int_G\xi_g(gh)dg
=\int_G\left(\sum_{h\in{}\bar{B}(e,S)}\xi_g(gh)\right)dg=\int_G\mathbf{1}dg=1.
\end{displaymath}
Moreover, if $k\in{}G$, $|k|<R$, we get that
\begin{align*}
\|f-kf\|_1 & =\sum_{h\in{}G}|f(h)-f(k^{-1}h)| \\
& =\sum_{h\in{}\bar{B}(e,S)\cup\bar{B}(k,S)}
\left|\int_G\xi_g(gh)dg-\int_G\xi_g(gk^{-1}h)dg\right| \\
& =\sum_{h\in{}\bar{B}(e,S)\cup\bar{B}(k,S)}
\left|\int_G\xi_g(gh)dg-\int_G\xi_{gk}(gh)dg\right| \\
& \leq\sum_{h\in{}\bar{B}(e,S)\cup\bar{B}(k,S)}
\int_G|(\xi_g-\xi_{gk})(gh)|dg \\
& =\int_G\left(\sum_{h\in{}\bar{B}(e,S)\cup\bar{B}(k,S)}
|(\xi_g-\xi_{gk})(gh)|\right)dg
<\int_G\epsilon{}dg=\epsilon.
\end{align*}
The existence of an $f$ with these properties shows that
$\textrm{diam}^{\mathcal{A}}_G(R,\epsilon)\geq\textrm{diam}^{\mathcal{F}}_G(R,\epsilon)$.  
\end{proof}

We now restrict our attention to countable discrete amenable groups, and their
direct products (which preserve amenability).\\
Given such a group $G$, we give $G\times{}G$ the $l^1$ metric, i.e.
\begin{displaymath}
d_{G\times{}G}((g_1,g_2),(h_1,h_2))=d_G(g_1,h_1)+d_G(g_2,h_2).
\end{displaymath} 
This construction clearly preserves bounded geometry and left-invariance of metrics, and
extends naturally to the direct product of $n$ copies of $G$, written $G^n$.

Note that by our assumption that $d_G$ is integer valued, 
\begin{displaymath}
|B_{d_G^n}(e,1)|=n|B_{d_G}(e,1)|;
\end{displaymath}
$B_{d_G^n}(e,1)$ consists of all $n$-tuples where only one element is not the identity, 
and that element has norm one in $G$.

The next theorem is the crucial result; it says roughly that if $G$ is (non-trivial),
amenable, countable discrete, then the direct products $G^n$ have amenability 
arbitrarily badly.  By \ref{aeqamen}, this implies that they have property A
arbitrarily badly too, and this will allow us to construct spaces that are
coarsely embeddable, but do not have property A.

\begin{suppinf}[\cite{no3}, 4.3]
Let $G$ be a non-trivial countable discrete group.\\
Then for any $0<\epsilon<2$, 
\begin{displaymath}
\lim\inf_{n\to\infty}\textrm{\emph{diam}}^{\mathcal{F}}_{G^n}(1,\epsilon)=\infty.
\end{displaymath}
\end{suppinf}

\begin{proof}
Assume otherwise for contradiction, so there exists $S\in\mathbb{N}$ such that for
infinitely many $n\in\mathbb{N}$ there is a function $f^{(n)}\in{}l^1(G^n)$ supported in
$\bar{B}(e,S)$ and that satisfies $\|f^{(n)}-gf^{(n)}\|_1<\epsilon$ whenever
$|g|=1$.\\
In the usual way, we assume that $f^{(n)}$ is positive for all $n$.

The first part of the proof constructs for arbitrary $m$, some $f_m$ with
`good variation properties', but `small' support.\\
More specifically, fix $\delta\leq\frac{2-\epsilon}{2S}$, and any $m\in\mathbb{N}$.  
We will construct $f_m$ with $(1,\epsilon+\delta)$ variation and supported
in $B(e,S-1)$.\\
For any $n\geq{m}$ consider the decomposition
\begin{displaymath}
G^n=G^m\times{}G^m\times...\times{}G^m\times{}G^r,
\end{displaymath}
where $0\leq{}r<m$.\\
for each $k=1,...,(n-r)/m$, we define the `$k^{th}$ boundary of $f^{(n)}$', 
$\partial_kf^{(n)}$, by first restricting $f^{(n)}$ to the set
\begin{displaymath}
\{g\in\textrm{supp}f^{(n)}:|g|=S,g_i\not=e ~\textrm{if and only if} 
~(k-1)m<i\leq{}mk\},
\end{displaymath}  
where $g_i$, $i=1,...,n$, is the $i^{th}$ coordinate of $g$.\\
$\partial_kf^{(n)}$ is then this restriction extended by zero to all of $G^n$.\\
Note that if $k\not=l$, we get that 
supp$\partial_kf^{(n)}\cup$supp$\partial_lf^{(n)}=\emptyset$, whence
\begin{displaymath}
\sum_{k=1}^{(n-r)/m}\|\partial_kf^{(n)}\|_1\leq\|f^{(n)}\|_1=1.
\end{displaymath}
It follows that for any $\epsilon'$ satisfying 
$\frac{\epsilon+2\epsilon'}{1-\epsilon'}\leq\epsilon+\delta$, there exists (large) 
$n\in\mathbb{N}$ and some $k$ with $\|\partial_kf^{(n)}\|_1<\epsilon'$.\\
Define positive $\phi\in{}S(l^1(G))$ by setting
\begin{displaymath}
\phi=\frac{f^{(n)}-\partial_kf^{(n)}}{\|f^{(n)}-\partial_kf^{(n)}\|_1}.
\end{displaymath}
Note then that if $|g|=1$, we get that
\begin{displaymath}
\|\phi-\lambda_g\phi\|_1=
\frac{\|(f^{(n)}-\lambda_gf^{(n)})+(\lambda_g\partial_kf^{(n)}-\partial_kf^{(n)})\|_1}
{\|f^{(n)}-\partial_kf^{(n)}\|_1}
<\frac{\epsilon+2\epsilon'}{1-\epsilon'}\leq\epsilon+\delta.
\end{displaymath}
Keeping the same $n$ and $k^{th}$ factor $G^m$ of $G^n$ where all this was performed,
and considering $G^n=G^m\times{}G^{n-m}$, we define $f_m:G^m\to\mathbb{C}$ by setting
\begin{displaymath}
f_m(g)=\sum_{h\in{}G^{n-m}}\phi(gh).
\end{displaymath}
This concludes the first part of the proof.

In the second part we show that $f_m$ has the required properties.\\
Note first that
\begin{displaymath}
\|f_m(g)\|_1=\sum_{g\in{}G^m}\sum_{h\in{}G^{n-m}}\phi(gh)=\sum_{g\in{}G^n}\phi(g)=1.
\end{displaymath}
Moreover, supp$f_m\subseteq\bar{B}(e,S-1)$, as $\phi$ must be zero on any element
of $G^m$ of length $S$.\\
Finally, note that if $k\in{}G^m$, $|k|=1$, we get that
\begin{align*}
\|f_m-\lambda_kf_m\|_1
& =\sum_{g\in{}G^m}|f_m(g)-f_m(k^{-1}g)|
=\sum_{g\in{}G^m}\left|\sum_{h\in{}G^{n-m}}(\phi(gh)-\phi(k^{-1}gh))\right| \\
& \leq\sum_{g\in{}G^n}|\phi(g)-\phi(k^{-1}g)|=\|\phi-\lambda_k\phi\|_1<\epsilon+\delta.
\end{align*}
As $m$ was arbitrary, we have a family $\{f_m\}_{m\in\mathbb{N}}$ with these properties.

In the third and final part of the proof, we use this construction to arrive at a
contradiction.\\
The construction can be repeated to get a family
$\{f_m\}_{m\in\mathbb{N}}$ such that each $f_m\in{}l^1(G^m)$ has $(1,\epsilon+n\delta)$ 
variation, and is supported in $\bar{B}(e,S-n)$ for any $1\leq{n}\leq{S}$.\\
Taking this to its logical conclusion, then, we get a family $\{f_m\}_{m\in\mathbb{N}}$,
where $f_m$ is supported on $\bar{B}(e,0)=\{e\}$ and has 
$(1,\epsilon+S\delta)$ variation for all $m$.\\
Recall, however, that by choice of $\delta$, 
\begin{displaymath}
\epsilon+S\delta=\epsilon+S\frac{2-\epsilon}{2S}<2.
\end{displaymath}
On the other hand, for any $g\in{}G\backslash\{e\}$ with $|g|=1$, we conclude that
$\|f_1-\lambda_gf_1\|=2$.  This is a contradiction.
\end{proof}

After the last theorem, it is now a relatively simple matter to construct coarsely
embeddable spaces without property A; our sole ingredient is a finite group (even the
two element group will do!).
So let $G$ be a finite (hence amenable), non-trivial group, and define
$\chi_G$ to be the disjoint union $\sqcup_{n=1}^{\infty}G^n$ equipped with a metric
$d$ that restricts to the $l^1$ metric on $G^n$ and satisfies
\begin{itemize}
\item $d(X_n,X_{n+1})\geq{}n+1$,
\item for $n\leq{}m$, $d(G^m,G^n)=\sum_{k=n}^{m-1}d(G^k,G^{k+1})$.
\end{itemize}

\begin{cenota}
Let $G$ be a finite non-trivial group.
\begin{itemize}
\item $\chi_G$ does not have property A.
\item $\chi_G$ coarsely embeds in $l^2$.
\end{itemize}
\end{cenota}

\begin{proof}
To prove the first part, assume for contradiction that $\chi_G$ did have property A,
whence $\textrm{diam}^{\mathcal{A}}_{\chi_g}(1,\epsilon)<\infty$ for all $0<\epsilon<2$.\\
Let $\xi$ be a map realising this property for some such $\epsilon$.\\
Then for all sufficiently large $n$, (specifically, 
$n>\textrm{diam}^{\mathcal{A}}_{\chi_g}(1,\epsilon)$), the restriction of $\xi$ to
$G^n\subseteq{}\chi_G$ shows that 
$\textrm{diam}^{\mathcal{A}}_{G^n}(1,\epsilon)$ is bounded away from infinity by
$\textrm{diam}^{\mathcal{A}}_{\chi_g}(1,\epsilon)$.\\
As $\textrm{diam}^{\mathcal{A}}_{G^n}(1,\epsilon)=\textrm{diam}^{\mathcal{F}}_{G^n}(R,\epsilon)$
by amenability of $G^n$, and the right hand side tends to infinity, this is a contradiction.

To prove the second part, note that any finite metric space, so in particular $G$, bi-Lipschitz
embeds in $l^1$, i.e. there exists $f:G\to{}l^1$ and $C>0$ such that for all $g,h\in{}G$,
$C^{-1}d(g,h)\leq{}d(f(g),f(h))\leq{}Cd(x,y)$.\\
It follows that $G^n$ bi-Lipschitz embeds in the direct sum of $n$ copies of $l^1$, with the
$l^1$ norm on the sum; write $(\sum_{i=1}^nl^1)_1$ for this.  Moreover, it does so with the same
constant $C$.\\
Hence $\chi_G$ coarsely embeds in
\begin{displaymath}
\left(\sum_{k=1}^{\infty}\left(\sum_{i=1}^kl^1\right)_1\right)_1.
\end{displaymath}
As this is isometrically isomorphic to $l^1$, and $l^1$ coarsely embeds in $l^2$ by the 
result of \ref{pto2}, this completes the proof.
\end{proof}

The simplest example of this construction is given by $\mathbb{Z}_2$, the 
two element group.  In this case, $\mathbb{Z}_2^n$ with the $l^1$ metric is (isometric to) the
vertices of the edge graph of an n-dimensional cube, with metric setting $d(v,w)$ equal to the minimal
number of edges joining $v,w$.

\subsection{Groups without A}
\label{nonagp}
We sketch a construction of groups that do not coarsely embed in Hilbert
space; such a group cannot, of course, have property A.
This is a surprising result, partly because it implies the existence of non-exact
groups; this was part of the motivation for the Guentner-Kaminker-Ozawa theorem we 
looked at in chapter \ref{analysis}.  The author does not pretend for a moment
to fully understand the construction, but the basic ideas are picturesque and reasonably easy to get a handle on.

Gromov has written a much
longer exposition of some of the ideas \cite{gr3}, which was expanded on in
\cite{sil} and \cite{oll} (and see these papers for further references). Despite this,
aspects of the theory remain controversial for some; no complete proof exists in the literature and
experts are divided on exactly how much remains to be done.  Nonetheless, Gromov's result is now widely believed to be true.

One starts with the idea of a \emph{random group}.  Given a group $G$ (for
concreteness, say $G=F_k$), we choose a probability measure $\mu$ on 
$\mathcal{P}(F_k)$.  
A random group is then a quotient $F_k/\langle\langle{}R\rangle\rangle$, where 
$\langle\langle{}R\rangle\rangle$ is the normal
subgroup generated by $R\in\mathcal{P}(F_k)$, a random subset chosen with respect to $\mu$.

Obviously, the characteristics of random subgroups will depend on $\mu$, and for many
$\mu$, random subgroups are trivial with probability one (the set $R$ is large enough to
generate all of $F_k$ with probability one).  This is not, however, always true, and
one can study the properties that a random group is likely to have.

Gromov then modifies this idea in the following way: take $F_N$ for some large $N$ and a graph 
$\mathcal{G}$.  Generators of $F_N$ are assigned randomly to edges in $\mathcal{G}$ and a map
$\phi:\pi_1(\mathcal{G})\to{}F_N$ defined to send each element of $\pi_1(\mathcal{G})$ 
(loop in $\mathcal{G}$) to 
the corresponding word in the generators of $F_N$.  The random group resulting
is $F_k/\langle\langle\phi(\pi_1(\mathcal{G}))\rangle\rangle$.  Clearly, it is likely to be trivial if $\mathcal{G}$ contains a 
large number of loops; a central problem to be overcome for the construction to work is to show that this group is not always (or almost always) trivial.

The crucial claim (for us) is that \emph{if $\mathcal{G}$ contains an expander, 
then random groups as above
do not coarsely embed in Hilbert space with non-zero probability}.  Intuitively, one might think of taking
the quotient by $\langle\langle\phi(\pi_1(\mathcal{G}))\rangle\rangle$ as creating (enough of) the same loops in the Cayley graph of $F_N$ as 
appear in $\mathcal{G}$, and thus (`probably') enough loops to obstruct its coarse
embeddability in Hilbert space.
Further, Gromov claims that one can find `quasi-random' groups
with this property that embed in some finitely presented group.  The conclusion, then,
is the existence of finitely presented groups that do not coarsely embed, and thus
do not have property A (and are not exact).

It is worth commenting that Gromov's groups are constructed as inductive limits
of Hyperbolic groups (see 
\cite{gr3}, 1.4 for hints); as such they do still satisfy the (strong) Novikov conjecture, so are not even 
a possible counterexample to this.  They would, however, provide counterexamples to (surjectivity in) certain versions 
of the Baum-Connes conjecture; see \cite{hls}.

\pagebreak

\section{Comments on the bibliography} \label{bib}
We make (brief) comments on all the bibliography items (some of which were not
directly used in the discussion above, but are related to the circle of ideas this
piece is interested in).  In no way are the
comments below meant to sum up the entire contents of any of the papers or books listed!  
We simply attempt to point out their relevance to property A and coarse embeddability.\\

\cite{ad} proves a result about a Gromov hyperbolic groups that leads to the conclusion of
their having property A.\\

\cite{adr} deals with amenable groupoids.  $X$ has property A if and only if the
translation groupoid $G(X)$ is amenable (cf. \cite{sty}), and amenability of an action is
the same thing as amenability of the appropriate groupoid, so the equivalences discussed
in this book are very useful and relevant.\\

\cite{bhv} is an introduction to property T groups, also containing a lot of
general information on e.g. kernels, representations and amenability in the
appendices.\\

\cite{bell2} proves that if $(X,x_0)$ is metric space with a basepoint of finite
asymptotic dimension on which
finitely generated $G$ acts by isometries in such a way that the \emph{$r$-stabiliser},
$\{g\in{}G:d(x_o,gx_0)<r\}$, has property A, then $G$ has property A.\\

\cite{bell} proves that if $G$ is the fundamental group of a developable complex
over a small category without loops $\mathcal{Y}$, the geometric realisation of the
development $|\mathcal{X}|$ has finite asymptotic dimension, and the vertex groups
of $\mathcal{Y}$ have property A, then $G$ has property A.\\

\cite{bd} covers some ideas relating to finite asymptotic dimension of discrete groups.

\cite{bl} describes aspects of the geometry of Banach space, much of which has relevance to 
coarse-embeddability-type problems\\

\cite{bdk} relates (in a somewhat different context to ours) coarse embeddings into
different $L^p$ spaces.\\

\cite{bg} proves that any bounded geometry space coarsely embeds into a reflexive, strictly
convex Banach space. \\

\cite{bo} is an introduction to the use of finite dimensional approximation
techniques in $C^*$-algebra theory.  The chapter on exact groups contains a lot
of material relevant to property A.\\

\cite{bnw} constructs so-called `partial translation structures' out of 
coarse spaces, and (with some assumptions) relates them to properties of
$C^*$-algebras.  It also contains a wealth of equivalent formulations of A.\\

\cite{cn} proves property A for any group acting properly and cocompactly on a CAT(0)
cube complex.\\

\cite{cdgy} gives two proofs that coarse embeddability is preserved under free
products. The second also works for property A.\\

\cite{cw} proves that if $X$ has property A, then any so-called \emph{ghost
operator} (see \cite{roe1}, p. 170) in $C^*_u(X)$ is compact. \\

\cite{cw2} proves that if $X$ has property A, then the controlled (in our context,
finite) propagation operators are dense in any ideal of $C^*_u(X)$, and some further
equivalent conditions.\\

\cite{ccjjv} gives a comprehensive survey of a-T-menable (or Haaagerup)
groups.  while containing little directly related to property A it is very
much one of the same circle of ideas.\\

\cite{dg} discusses closure properties of the class of (countable discrete)
groups that coarsely embed in Hilbert space.  A lot of the proofs go
over to the property A case.\\

\cite{dg2} has as its main result that if $G$ is hyperbolic relative
to subgroups $H_1,...,H_n$ that coarsely embed, then $G$ coarsely embeds.
It also contains the formulation of property A in terms of partitions
of unity, and various `gluing'-type lemmas that it can be used to prove.\\

\cite{dr} proves that a group whose asymptotic dimension function grows at
a polynomial rate has property A, and equivalence of property A with
admitting an \emph{$l^p$-spherical anti-\v{C}ech approximation} for
$1\leq{}p<\infty$ (for bounded geometry geodesic spaces).\\

\cite{dr2} proves that the countably generated free group does not coarsely
embed when equipped with some word metric.  This is despite its having property
A when equipped with any bounded geometry, left-invariant metric.\\

\cite{dgly} discusses some aspects of coarse embeddability, and a constructs
a non-coarsely embeddable, separable metric space.\\

\cite{far} proves that Thompson's groups are a-T-menable.  It is a major 
unsolved problem to decide
whether Thompson's group $F$ (which is finitely presented) is amenable or has
property A.\\

\cite{frr} covers a lot of information about various versions of the Novikov
conjecture, including a detailed survey and history by the editors.\\

\cite{gr1} was one of the main impetuses to the study of large-scale properties
of groups.  It contains the seed of a lot of the ideas discussed above.\\

\cite{gr2} is mainly about geometric properties of manifolds.  Its last few
pages, however, sketch a construction of groups that do not coarsely
embed in Hilbert space.\\

\cite{gr3} is a much longer exposition of the ideas mentioned at the end of the
previous article.\\

\cite{gu} proves property A for countable discrete groups with one relation.\\

\cite{ghw} proves property A for any countable linear group over any field.\\

\cite{gk1} and \cite{gk2} relate exactness and property A, specifically
proving that exactness implies coarse embeddability.\\

\cite{gk3} is a survey of the relationships between (coarse) geometric
(A, amenability, coarse embeddability) and
analytic (exactness, nuclearity) properties of groups, 
including notes on the Baum-Connes and Novikov conjectures.\\

\cite{gk4} defines the `Hilbert space compression' of a space and relates
it to property A.\\

\cite{hig} proves (amongst other things) split injectivity of the Baum-Connes assembly map for countable property A groups using
KK-theoretic techniques.\\

\cite{hls} provides counterexamples to certain versions of the Baum-Connes and coarse Baum-Connes
conjectures.  These use ideas relating to exactness from our chapter \ref{analysis} and 
Gromov's groups, from section \ref{nonagp}.\\

\cite{hr} proves that a group has property A if and only if it admits an 
amenable action on a compact Hausdorff space.  It also proves that finite
asymptotic dimension implies A.\\

\cite{hr2} is included for its notes on the construction of the Roe algebra
(which it does rather more generally than we did above).\\

\cite{jr} proves that $l^p$ does not coarsely embed in $l^2$ for $p>2$.\\

\cite{ky} proves the injectivity of the coarse assembly map for spaces which admit a 
coarse embedding into a uniformly convex Banach space.\\

\cite{laf} formulates a stronger version of property T.  Of particular relevance to this piece,
it also shows that box spaces constructed from such a group cannot coarsely embed in any uniformly
convex Banach space.\\

\cite{lub} discusses expanding graphs, and their relationship to property T and other
areas of mathematics.\\

\cite{lps} constructs expanding graphs using the (proved) Ramanujan conjecture.\\

\cite{mar} constructs expanding graphs using a property-T-like condition.\\

\cite{no1} proves that $l^p$ coarsely embeds in $l^2$ for $0<p<2$.\\

\cite{no2} proves that $l^2$ coarsely embeds in $l^p$ for $1\leq{}p<2$.\\

\cite{no4} expands on the function from \cite{no3} that measures `how well' a space has 
property A, and relates it to other functions associated to amenability, growth rate
and finite asymptotic dimension. \\

\cite{no3} constructs (locally finite) coarsely embeddable spaces that do not 
have property A.\\

\cite{oll} fills in some of the details about random groups, following the ideas
given in \cite{gr2} and \cite{gr3}.\\

\cite{oz} proves that a groups is exact if and only if it has property A, via
amenable actions.\\

\cite{oz2} proves that a group that is hyperbolic relative to a collection
of property A subgroups has property A.\\

\cite{rvw} introduces the zig-zag graph product and uses it to construct
expanders.\\

\cite{roe3} discusses the
coarse Baum-Connes conjecture and its relation to index theory, as well as 
being the source for other coarse geometric ideas.\\

\cite{roe1} is a general introduction to coarse geometry.  It also has a section
covering property A and coarse embeddability.\\

\cite{roe2} constructs warped cones, and discusses conditions under which they
have property A.\\

\cite{roe4} proves that hyperbolic groups have finite asymptotic dimension.\\

\cite{sch1} and \cite{sch2} are the historical basis for the theory of kernels of
positive type.\\

\cite{ser} was used for its discussions of groups acting trees, which relate to
amalgamated free products and HNN extensions.\\

\cite{sil} fills in some of the details from \cite{gr2} and \cite{gr3}, focusing in
particular for property T for random groups.\\

\cite{sty} relates coarse geometry, the translation groupoid $G(X)$ of a space $X$ and the
coarse Baum-Connes conjecture.  This treatment unifies a lot of the ideas in our
discussion.  The paper also removes the necessity of the finiteness
assumption on $BG$ and the finite generation assumption on $G$ for the main result of \cite{y1}, using an idea from \cite{hig}.\\

\cite{tes} contains some applications of property A to the distortion of embeddings of certain metric spaces into $L^p$-spaces.\\

\cite{tu} proves certain closure properties of A, in particular spaces of subexponential
growth and excisive unions, which we don't mention above.  It also proves that if $G$ acts
on a tree in such a way that each vertex stabiliser has A, then $G$ has A.\\

\cite{wass} is a general discussion of (and the original source for many ideas about)
exact $C^*$-algebras.\\

\cite{y1} introduces property A, and proves the strong Novikov
conjecture for a finitely generated group $G$ that coarsely embeds and is 
such that the
classifying space $BG$ has the homotopy type of a finite CW complex (the
finiteness assumption was later removed in \cite{sty}).

\pagebreak


\begin{thebibliography}{120}

\bibitem[Ada94]{ad} S. Adams, \emph{Boundary amenability for word hyperbolic groups and an application
to smooth dynamics of simple groups}, Topology \textbf{33} (1994), pp. 765-783.

\bibitem[A-DR00]{adr} C. Anantharaman-Delarcohe and J. Renault, \emph{Amenable groupoids},
L'Enseignement Mathematique \textbf{36} (2000), with a foreword by George Skandalis and Appendix
B by E. Germain.

\bibitem[BHV]{bhv} B. Bekka, P. De La Harpe and A. Valette, \emph{Kazhdan's property
(T)}, manuscript available at 
http://perso.univ-rennes1.fr/bachir.bekka/KazhdanTotal.pdf.

\bibitem[Bell03]{bell2} G. Bell,  \emph{Property A for groups acting on metric spaces},
Topology Appl. \textbf{130} (2003), no. 3, pp. 239-251.

\bibitem[Bell05]{bell} G. Bell, \emph{Asymptotic properties of groups acting on complexes},
Proc. Amer. Math. Soc. \textbf{133} (2005), no. 2, pp. 387-396 (electronic).

\bibitem[BD01]{bd} G. Bell and A. Dranishnikov, \emph{On asymptotic dimension of discrete
groups}, Algr. Geom. Topl. \textbf{1} (2001), pp. 57-71.

\bibitem[BL00]{bl} Y. Benyamini and J Lindenstrauss, \emph{Geometric nonlinear
functional analysis}, Vol. 1, AMS Colloquium publications, Vol. 48, 2000.

\bibitem[BDK66]{bdk} J. Bretagnolle, D. Dacunha-Castelle and J. Krivine, \emph{Lois stables
et espaces $L^p$}, Ann. Inst. H. Poincar\'{e} \textbf{2} (1966), pp. 231-259.

\bibitem[BG05]{bg} N. Brown and E. Guentner, \emph{Uniform embeddings of bounded geometry
spaces into reflexive Banach spaces}, Proc. Amer. Math. Soc. \textbf{133} (2005), No. 7, 
pp. 2045-2050.

\bibitem[BO06]{bo} N. Brown and N. Ozawa, \emph{$C^*$-algebras and finite
dimensional approximations}, preprint, 2006.

\bibitem[BNW06]{bnw} J. Brodzki, G. A. Niblo and N. J. Wright, \emph{Property A, 
partial translation structures and uniform embeddings in groups}, preprint 2006.

\bibitem[CN05]{cn} S. Campbell and G. Niblo, \emph{Hilbert space compression 
and exactness of discrete groups},
J. Funct. Anal. \textbf{222} (2005), no. 2, pp. 292-305.

\bibitem[CDGY03]{cdgy} X. Chen, M. Dadarlat, E. Guentner and G. Yu, \emph{Uniform
embeddability and exactness of free products}, J. Funct. Anal. \textbf{205} (2003), pp. 293-314.

\bibitem[CW04]{cw2} X. Chen and Q. Wang,  \emph{Ideal structure of uniform Roe algebras 
of coarse spaces}, J. Funct. Anal. \textbf{216} (2004), no. 1, pp. 191-211.

\bibitem[CW05]{cw} X. Chen and Q. Wang, \emph{Ghost ideals in uniform Roe algebras of
coarse spaces}, Arch. Math. (Basel), \textbf{84} (2005), no. 6, pp. 519-526.

\bibitem[CCJJV01]{ccjjv} P.-A. Cherix, M. Cowling, P. Jolissaint, P. Julg and A. Valette,
\emph{Groups with the Haagerup property (Gromov's a-T-menability)}, 
Birkh\"{a}user, Boston, 2001.

\bibitem[DG03]{dg} M. Dadarlat and E. Guentner, \emph{Constructions preserving Hilbert
space uniform embeddability of discrete groups}, Vol. 355, Trans. Amer. Math. Soc. 
No. 8 (2003), pp. 3253-3275.

\bibitem[DG06]{dg2} M. Dadarlat and E. Guentner, \emph{Uniform embeddability of relatively
hyperbolic groups}, preprint 2006.

\bibitem[Dra02]{dr} A. Dranishnikov, \emph{Anti-\v{C}ech approximation in coarse geometry},
IHES, preprint, 2002.

\bibitem[DGLY]{dgly} A. Dranishnikov, G. Gong, V. Lafforgue and G. Yu, 
\emph{Uniform embeddings
into Hilbert space and a question of Gromov}, to appear in Canad. Math. Bull.

\bibitem[Dra]{dr2} A. Dranishnikov, \emph{On Generalized amenability}, preprint.

\bibitem[Far03]{far} D. Farley, \emph{Proper isometric actions of Thompson's groups on Hilbert
space}, Int. Math. Res. Notes, 2003, pp. 2409-2414.

\bibitem[FRR95]{frr} S. Ferry, A. Ranicki and J. Rosenberg (eds.), \emph{Novikov
conjectures, index theorems and rigidity}, London Math. Soc. Lecture Notes,
no.s 226, 227 Cambridge University Press, 1995.
 
\bibitem[Gro93]{gr1} M. Gromov, \emph{Asymptotic invariants of infinite groups}, 
London Math. Soc., Lecture Notes Series 182, Geometric group theory, Vol. 2,
CUP 1993.

\bibitem[Gro00]{gr2} M. Gromov, \emph{Spaces and questions}, GAFA 2000, Special Volume,
Part I, Birkh\"{a}user Verlag, Basel (2000), pp. 118-161.

\bibitem[Gro01]{gr3} M. Gromov, \emph{Random walks in random groups}, Geom. Funct. Anal., 
\textbf{13} (2003), pp. 73-146.

\bibitem[Gue02]{gu} E. Guentner, \emph{Exactness of the one relator groups},
Proc. Amer. Math. Soc. \textbf{130} (2002), no. 4, pp. 1087-1093 (electronic).

\bibitem[GHW04]{ghw} E. Guentner, N. Higson and S. Weinberger, \emph{The Novikov conjecture for
linear groups}, Noncommutative Geometry, 137-251, Lecture notes in Math., \textbf{1831}, Springer,
Berlin, 2004. 

\bibitem[GK02a]{gk1} E. Guentner and J. Kaminker, \emph{Exactness and the Novikov
conjecture}, Topology \textbf{41} (2002), no. 2, pp. 411-18.

\bibitem[GK02b]{gk2} E. Guentner and J. Kaminker, \emph{Addendum to
`Exactness and the Novikov conjecture'}, Topology \textbf{41} (2002), pp. 411-418.

\bibitem[GK03]{gk3} E. Guentner and J. Kaminker, \emph{Analytic and geometric
properties of groups}, in \emph{Noncommutative Geometry: Lectures at the C.I.M.E.
summer school held in Martina Franca, Italy, September 3-9, 2000}, Lecture notes in
Mathematics vol. 1831, Springer-Verlag Heidelberg, 2003, pp. 253-262.

\bibitem[GK04]{gk4} E. Guentner and J. Kaminker, \emph{Exactness and uniform embeddability
of discrete groups}, Journal Lon. Math. Soc. \textbf{70} (2004) pp. 703-718.

\bibitem[Hig00]{hig} N. Higson, \emph{Bivariant K-theory and the Novikov conjecture}, Geom. Funct. Anal., \textbf{10} (2000), pp. 563-581.

\bibitem[HLS02]{hls} N. Higson, V. Lafforgue, G. Skandalis, \emph{Counterexamples to the Baum-Connes
conjecture}, Geom. Funct. Anal., \textbf{12} (2002), 330-354.

\bibitem[HR00a]{hr} N. Higson and J. Roe, \emph{Amenable group actions and the
Novikov conjecture}, J. Reine Agnew. Math. \textbf{519} (2000).

\bibitem[HR00b]{hr2} N. Higson and J. Roe, \emph{Analytic K-homology}, Oxford Mathematical
Monographs, Oxford Science Publication, OUP, Oxford, 2000.

\bibitem[JR04]{jr} W. Johnson and N. Randrianarivony, \emph{$l^p$ $(p>2)$ does not coarsely 
embed into a Hilbert space}, Proc. Amer. Math. Soc., \textbf{134} (2005), No. 4, pp. 1045-1050.

\bibitem[KYu06]{ky} G. Kasparov and G. Yu, \emph{The coarse geometric Novikov conjecture 
and uniform convexity}, Adv. in Math., \textbf{206} (2006), Issue 1, pp. 1-56. 

\bibitem[Laf06]{laf} V. Lafforgue, \emph{Un renforcement de la propri\'et\'e (T)}, preprint.

\bibitem[Lub94]{lub} A. Lubotzky, \emph{Discrete groups, expanding graphs and
invariant measures}, Birk\"a{}user, Boston, 1994.

\bibitem[LPS86]{lps} A. Lubotzky, R. Philips and P. Sarnak, \emph{Ramanujan conjectures
and explicit constructions of expanders}, Proc. Symp on Theo. of Comp. Sci.
(STOC), \textbf{86} (1986), pp. 240-246.

\bibitem[Mar75]{mar} G. A. Margulis, \emph{Explicit construction of concentrators},
Probl. of Inform. Transm. \textbf{10} (1975), pp. 325-332.

\bibitem[Now05]{no1} P. Nowak, \emph{Coarse embeddings of metric spaces into Banach spaces},
Proc. Amer. Math. Soc. \textbf{133} (2005), no. 9, pp. 2589-2596.

\bibitem[Now06]{no2} P. Nowak, \emph{On coarse embeddability into $l^p$ spaces and a conjecture
of Dranishnikov}, Fund. Math. \textbf{189} (2006), pp. 111-6.

\bibitem[Now07a]{no4} P. Nowak, \emph{On exactness and isoperimetric properties of discrete
groups}, J. Funct. Anal., \textbf{243} (2007), no. 1, pp. 323-344.

\bibitem[Now07b]{no3} P. Nowak, \emph{Coarsely embeddable metric spaces without property A}, J. Funct. Anal.,
\textbf{252} (2007), no. 1, pp. 126-136.

\bibitem[Oll04]{oll} Y. Ollivier,\emph{Sharp phase transition theorems for hyperbolicity of random groups},
Geom. Funct. Anal. \textbf{14} (2004), no. 3, pp. 595-679.

\bibitem[Oz00]{oz} N. Ozawa, \emph{Amenable actions and exactness for discrete groups},
C. R. Acad. Sci. Paris Ser. I Math., \textbf{330} (2000), pp. 691-5.

\bibitem[Oz06]{oz2} N. Ozawa, \emph{Boundary amenability of relatively hyperbolic groups}, 
Topology Appl., \textbf{153} (2006), 2624-2630.

\bibitem[RVW02]{rvw} O. Reingold, S. Vadhan and A. Wigderson, \emph{Entropy waves,
the zig-zag product, and new constant degree expanders}, Ann. of
Math. \textbf{155} (2002), pp. 157-187.

\bibitem[Roe96]{roe3} J. Roe, \emph{Index theory, coarse geometry and the topology of
manifolds}, CBMS Regional Conference Series in Mathematics, AMS 1996.

\bibitem[Roe03]{roe1} J. Roe, \emph{Lectures on coarse geometry}, University Lecture Series,
31. AMS Providence, Rhode Island, 2003.

\bibitem[Roe05a]{roe2} J. Roe, \emph{Warped cones and property A}, Geometry and Topology, Vol.
\textbf{9} (2005), pp. 163-78.

\bibitem[Roe05b]{roe4} J. Roe, \emph{Hyperbolic groups have finite asymptotic dimension}, 
Proc. AMS, Vol. \textbf{133} no. 9, pp. 2489-2490.

\bibitem[Sch37]{sch1} I. Schoenberg, \emph{On certain metric spaces arising from Euclidean
spaces by a change of metric and their imbedding in Hilbert space}, Ann. Math. \textbf{38} 
(1937), pp. 787-93.

\bibitem[Sch38]{sch2} I. Schoenberg, \emph{Metric spaces and positive definite functions},
Trans. Am. Math. Soc. \textbf{44} (1938), pp. 522-36. 

\bibitem[Ser80]{ser} J-P. Serre, \emph{Trees}, Springer-Verlag, 
New York-Heidelberg-Berlin, 1980.

\bibitem[Sil03]{sil} L. Silberman, \emph{Addendum to: `Random walks in random groups'},
Geom. Funct. Anal. \textbf{13} (2003), no. 1, pp. 147-177.

\bibitem[STY02]{sty} G. Skandalis, J-L. Tu and G. Yu, \emph{The coarse Baum-Connes conjecture
and groupoids}, Topology \textbf{41} (2002), pp. 807-34.

\bibitem[Tes07]{tes} R. Tessera, \emph{Quantitative property A, Poincare inequalities, $L^p$-compression and $L^p$-distortion for metric measure spaces},
preprint, 2007.

\bibitem[Tu01]{tu} J-L. Tu, \emph{Remarks on Yu's property A for discrete metric spaces
and groups}, Bull. Soc. Math. France \textbf{129} (2001), pp. 115-39.

\bibitem[Wass94]{wass} S. Wasserman, \emph{Exact $C^*$-algebras and related topics},
Seoul National University Research Institute of Mathematics Global
Analysis Research centre, Seoul, 1994.

\bibitem[Yu00]{y1} G. Yu, \emph{The coarse Baum-Connes conjecture for spaces which
admit a uniform embedding into Hilbert space}, Inventiones \textbf{139} (2000), pp. 201-240.
 



\end{thebibliography}
\end{document}